\pdfoutput=1
\RequirePackage{ifpdf}
\ifpdf 
\documentclass[pdftex]{sigma}
\else
\documentclass{sigma}
\fi

\numberwithin{equation}{section}

\newtheorem{Theorem}{Theorem}[section]
\newtheorem{Corollary}[Theorem]{Corollary}
\newtheorem{Lemma}[Theorem]{Lemma}
\newtheorem{Proposition}[Theorem]{Proposition}
{\theoremstyle{definition}
\newtheorem{Definition}[Theorem]{Definition}
\newtheorem{Example}[Theorem]{Example}
\newtheorem{Remark}[Theorem]{Remark}
\newtheorem{Notation}[Theorem]{Notation}
}
\DeclareMathOperator{\sign}{sign}

\begin{document}

\newcommand{\arXivNumber}{1304.2284}

\allowdisplaybreaks

\renewcommand{\thefootnote}{$\star$}

\renewcommand{\PaperNumber}{055}

\FirstPageHeading

\ShortArticleName{Ordered $*$-Semigroups and a~$C^*$-Correspondence for a~Partial Isometry}

\ArticleName{Ordered $\boldsymbol{*}$-Semigroups and a~$\boldsymbol{C^*}$-Correspondence\\
for a~Partial Isometry\footnote{This paper is a~contribution to the Special Issue on Noncommutative Geometry and Quantum
Groups in honor of Marc A.~Rief\/fel.
The full collection is available at
\href{http://www.emis.de/journals/SIGMA/Rieffel.html}{http://www.emis.de/journals/SIGMA/Rieffel.html}}} 

\Author{Berndt BRENKEN}

\AuthorNameForHeading{B.~Brenken}

\Address{Department of Mathematics and Statistics, University of Calgary,\\
Calgary, Canada T2N 1N4} \Email{\href{mailto:bbrenken@ucalgary.ca}{bbrenken@ucalgary.ca}}

\ArticleDates{Received August 30, 2013, in f\/inal form May 22, 2014; Published online May 31, 2014}

\Abstract{Certain $*$-semigroups are associated with the universal $C^*$-algebra generated by a~partial isometry, which
is itself the universal $C^*$-algebra of a~$*$-semigroup.
A~fundamental role for a~$*$-structure on a~semigroup is emphasized, and ordered and matricially ordered \mbox{$*$-semigroups}
are introduced, along with their universal $C^*$-algebras.
The universal $C^*$-algebra generated by a~partial isometry is isomorphic to a~relative Cuntz--Pimsner $C^*$-algebra of
a~$C^*$-correspondence over the $C^*$-algebra of a~matricially ordered $*$-semigroup.
One may view the $C^*$-algebra of a~partial isometry as the crossed product algebra associated with a~dynamical system
def\/ined by a~complete order map modelled by a~partial isometry acting on a~matricially ordered $*$-semigroup.}


\Classification{46L05; 46L08; 20M30; 06F05; 46L55}

\vspace{-3mm}

{\small \tableofcontents}

\renewcommand{\thefootnote}{\arabic{footnote}} 
\setcounter{footnote}{0}

\section{Introduction}

There are at least two motivations for this project.
One is to emphasize a~particular direction in the study of $C^*$-algebras associated with semigroups, illustrating roles
for universal $C^*$-algebras, particularly through a~consideration of star and order structures on semigroups.
The second is the delineation of an underlying algebraic structure, basically via dynamical considerations, for the
universal $C^*$-algebra of a~partial isometry.
This structure intimately involves the order structure of a~$*$-semigroup, so provides some evidence for, and
justif\/ication of, this direction of inquiry.

A~structure for the universal \mbox{$C^*$-algebra} of a~partial isometry is that of a~Cuntz--Pimsner \mbox{$C^*$-algebra} associated
with a~$C^*$-correspondence, so a~representation of a~$C^*$-algebra as adjointable operators on a~Hilbert module.
In our context the Cuntz--Pimsner $C^*$-algebra can be thought of as type of crossed product $C^*$-algebra arising from
an action by a~completely positive map, here one modelled through a~partial isometry, on a~particular nonunital
$C^*$-sub\-algebra.
The $C^*$-algebras involved occur as universal $C^*$-algebras associated with particular families of contractive
$*$-representations of certain $*$-semigroups.
These families of contractive $*$-representations arise through a~consideration of various order structures on
$*$-semigroups.

The universal $C^*$-algebras considered arise through contractive representations of a~semigroup.
They explicitly involve a~star structure, and an overlying order, and complete order structure.
The semigroups arising from partial isometries are developed in this context.
The particular operator algebras introduced here are further examples of Banach algebras associated with certain
semigroups (cf.~\cite{bd, dls}).
In the case of a~discrete group~$G$, viewed as a~$*$-semigroup with a~natural order and complete order, the various
universal $C^*$-algebras associated with~$G$ all coincide with the usual group $C^*$-algebra of~$G$.
For one of the $*$-semigroups considered below these universal $C^*$-algebras also coincide, while for other
$*$-semigroups they dif\/fer.
It is hoped that further interesting structure results for other $C^*$-algebras will be developed, as well as
investigations involving $*$-semigroup and ordered $*$-semigroup universal $C^*$-algebras.
Noting that investigations of actions on algebras of operators via isometries and their induced endomorphisms has led to
large areas of ongoing inquiry~\cite{a},
there may be potential for investigations involving actions by partial isometries.
We anticipate that these universal $C^*$-algebras of ordered $*$-semigroups will point to potentially interesting
directions for further investigations.

A~brief opening Section~\ref{section1} introduces a~universal $C^*$-algebra for contractive $*$-representations of a~$*$-semigroup.
As a~basic illustration of how this universal $C^*$-algebra may capture useful structure of a~semigroup,
a~$*$-semigroup~$A$ is introduced and its universal $C^*$-algebra is the universal $C^*$-algebra of a~partial isometry.
Section~\ref{section2} introduces particular $*$-subsemigroups of~$A$ and the notion of an irreducible element in these
$*$-semigroups.
For example, the $*$-subsemigroup $A^{0}$ is the kernel of a~tracial map.
The main $*$-semigroups of our current interest, the proper subsemigroup $D_{1}$ and its unital $*$-subsemigroup $D_{0}$
of $A^{0}$, are def\/ined in terms of irreducible generators.
There is a~$*$-semigroup homomorphism $\omega$ from the unital $*$-semigroup~$D_{0}$ to the nonunital $*$-semigroup
$D_{1}$ along with a~$*$-map $\alpha:D_{1}\rightarrow D_{0}$ with $\alpha\circ\omega$ the identity homomorphism on
$D_{0}$.
The map~$\alpha$ is modelled on the completely positive map induced by a~single partial isometry, and is our main focus.

Section~\ref{section3} concerns partially ordered $*$-semigroups and their order maps, as well as order maps satisfying
the Schwarz inequality and order representations in a~$C^*$-algebra.
A~partial order for the $*$-semigroup $D_{1}$ is described.
The universal $C^*$-algebra of a~partially ordered $*$-semigroup exists for the semigroups we consider, however, this is
generally not suf\/f\/icient to obtain a~bounded representation of a~$*$-semigroup.
For the semigroups considered here the order representations are automatically via contractions, so the universal
$C^*$-algebra is a~quotient of the initial universal $C^*$-algebra.
A~further extension is described in Section~\ref{section4}, namely matricial partial orders and the corresponding
complete order maps on $*$-semigroups.
The map $\alpha$ is a~complete order map, and the universal $C^*$-algebra for complete order representations of
a~matricially ordered $*$-semigroup is def\/ined.
Forming the Hilbert module over the universal matricially ordered \mbox{$C^*$-algebra} for the semigroup $D_{1}$ then yields
a~complete order representation of $D_{1}$, resulting in a~$C^*$-correspondence $\mathcal{E}$ over this universal $C^*$-algebra.
The section concludes with a~closer look at a~matrix decomposition for elements in the subset of matrices used to def\/ine
the matricial order structure on $D_{1}$.

The f\/inal section shows that a~relative Cuntz--Pimsner $C^*$-algebra~\cite{ms}  associated with this
$C^*$-correspondence is isomorphic to the universal $C^*$-algebra
generated by a~partial isometry.
Note that a~general theorem~\cite{aee} already  implies that the $C^*$-algebra generated by a~partial isometry is the
Cuntz--Pimsner algebra over the $C^*$-subalgebra generated by the image of the larger semigroup~$A^{0}$.
The ideal used to describe this relative Cuntz--Pimsner algebra is contained in the standard ideal $J_{\mathcal{E}}$,
the largest ideal on which the restriction of the left action for the correspondence $\mathcal{E}$ is an injection into
compact adjointable maps of the Hilbert module.
A~consequence is that the complete order universal $C^*$-algebra of $D_{1}$ is isomorphic to its image in the universal
$C^*$-algebra generated by a~partial isometry.
It follows that each of the various universal $*$-semigroup $C^*$-algebras of $D_{1}$ are not exact.

\subsection*{Notation}

A~semigroup is a~set with a~binary associative operation.
For present purposes they are assumed to have the discrete topology.
Denote by $\mathbb{N}$ the usual additive semigroup of nonzero natural numbers, while $\mathbb{N}_{0}$ is the unital
semigroup with adjoined unit, denoted~$0$ for this abelian case.
To distinguish certain copies of the semigroup $\mathbb{N}$ use $\mathbb{N}^{+}$ and $\mathbb{N}^{-}$ to denote the
abelian semigroups of strictly positive and strictly negative integers respectively.
If a~unit is adjoined to these semigroups it is potentially useful to distinguish these, in which case they are~$0^{+}$
and~$0^{-}$ respectively.

For $\mathcal{H}$ a~Hilbert space, $\mathcal{B}\left(\mathcal{H}\right) $ denotes the $C^*$-algebra of bounded operators
on $\mathcal{H}$, while $\mathcal{K}$ denotes the norm closed $*$-subalgebra of compact operators, an ideal in
$\mathcal{B}\left(\mathcal{H}\right)$,
where ideal in the $C^*$-context means a~closed two-sided ideal.
For a~set~$X$ and $n,m\in \mathbb{N}$, $\mathrm{M}_{m,n}(X)$ is the set of $m\times n$ matrices with entries in~$X$,
while $\mathrm{M}_{n}(X)$ is the set of $n\times n$ matrices.
If~$X$ is a~$C^*$-algebra then so is $\mathrm{M}_{n}(X)$.
For~$X$ a~set equipped with an involution, the set of selfadjoint elements of~$\mathrm{M}_{n}(X)$ is
$\mathrm{M}_{n}(X)^{\rm sa}$.
If $\varphi:X\rightarrow Y$ is a~map, then $\varphi_{n}:\mathrm{M}_{n}(X)\rightarrow\mathrm{M}_{n}(Y)$ denotes the map
$\varphi_{n}([x_{i,j}])=[\varphi(x_{i,j})]$.

For results and conventions on $C^*$-modules we follow Lance~\cite{l}; so if~$B$ is a~$C^*$-algebra a~Hilbert~$B$-module
$\mathcal{E}_{B}$ is a~Banach space $\mathcal{E}$ which is a~right~$B$-module with a~$B$-valued inner product
$\langle\,,\,\rangle_{B}$, denoted $\langle\,,\,\rangle $ if the context is clear.
The norm on $\mathcal{E}$ is given by $\Vert x\Vert^{2}=\Vert \langle x,x\rangle \Vert
$, $(x\in\mathcal{E})$; $\mathcal{L}(\mathcal{E})$ denotes the $C^*$-algebra of adjointable linear operators def\/ined on
$\mathcal{E}$.
In analogy with the case when~$A$ is the complex numbers, $\mathcal{K}(\mathcal{E})$ is the closed two-sided ideal of `compact'
operators $\overline{\rm span}\{\theta_{x,y}^{\mathcal{E}}|x,y\in\mathcal{E}\}$,
where $\theta_{x,y}^{\mathcal{E}}(z)=x\langle y,z\rangle $, $(z\in\mathcal{E})$.
If $\mathcal{E}$ is a~Hilbert~$B$ - module the linear span of $\{\langle x,y\rangle |x,y\in\mathcal{E}\}$,
denoted $\langle \mathcal{E},\mathcal{E} \rangle $, has closure a~two-sided ideal of~$B$.
Note that $\mathcal{E} \langle \mathcal{E},\mathcal{E}\rangle $ is dense in $\mathcal{E}$~\cite{l}.
The Hilbert module $\mathcal{E}$ is called full if $\langle \mathcal{E},\mathcal{E}\rangle $ is dense in~$B$.
If~$B$ is a~$C^*$-algebra then $B_{B}$ refers to the Hilbert module~$B$ over itself, where $\langle
a,b\rangle =a^{\ast}b$ for $a,b\in B$.
A~$C^*$-correspondence over a~$C^*$-algebra~$B$, denoted $_{B}\mathcal{E}_{B}$, is a~Hilbert~$B$-module
$\mathcal{E}_{B}$ along with a~specif\/ied $*$-homomorphism $\phi:B\rightarrow\mathcal{L}(\mathcal{E}_{B})$.

\section[$C^*$-algebras of a~$*$-semigroup]{$\boldsymbol{C^*}$-algebras of a~$\boldsymbol{*}$-semigroup}
\label{section1}

A~semigroup~$S$ is a~set with an associative binary operation.
An element~$x$ in a~semigroup~$S$ with at least two elements is a~zero element, or absorbing element, if $xs=sx=x$ for
all $s\in S$.
Such an absorbing element is clearly unique.
The stipulation that~$S$ have at least two elements follows the convention that the element in the trivial semigroup of
one element is not referred to as a~zero~\cite{h}.
There is at most one such element, and we use the convention that a~homomorphism of semigroups with absorbing element
$0$ maps the zero element to the zero element.
For a~semigroup~$S$ without an absorbing element zero, the semigroup $S^{0}$ is the semigroup~$S$ with a~zero adjoined.
Most of the examples of semigroups considered below do not have an absorbing element zero.
A~homomorphism of semigroups which is one to one is called a~monomorphism.
If~$T$ is a~semigroup with absorbing element $0$ and $\beta:S\rightarrow T$ is a~semigroup monomorphism then if there is
an $s\in S$ with $\beta(s)=0$ it follows that~$s$ must be an absorbing element for~$S$, so must be a~(so only) zero
for~$S$.

A~$*$-semigroup is a~semigroup equipped with an involutive antihomomorphism, denoted~$*$.
There may be several involutive maps on a~given semigroup, as examples below indicate.
\mbox{A~$*$-ho\-mo\-morphism}
$\beta:S\rightarrow T$ of $*$-semigroups is a~semigroup homomorphism with $\beta(s^{\ast})=(\beta(s))^{\ast}$, $s\in S$.
If~$S$ has an identity it follows that the $*$-operation must leave it f\/ixed, and if~$S$ has a~left or right identity
the $*$-operation implies~$S$ has an identity.
Similarly if the semigroup~$S$ has a~left or right absorbing element
then it has an absorbing element, and the $*$-operation must leave it f\/ixed.
An element~$s$ of a~$*$-semigroup~$S$ is called selfadjoint if $s^{\ast}=s$, and a~$*$-homomorphism necessarily maps
selfadjoint elements to selfadjoints.

The $*$-operation on a~$*$-semigroup~$S$ yields a~$*$-algebra structure on the vector space $\mathbb{C}[S]$ of complex
valued functions of f\/inite support on~$S$, with $f^{\ast}(s)=\overline{f(s^{\ast})}$ for~$f$ such a~function.
The same $*$-operation is isometric on the Banach $*$-algebra $l_{1}(S)$, with the algebra structure on $l_{1}(S)$
arising from the convolution operation.

Note that if~$S$ is a~group or an inverse semigroup then the inverse operation is a~natural $*$-operation on~$S$ since
it satisf\/ies $(ab)^{-1}=b^{-1}a^{-1}$.
For an inverse semigroup~$S$ this identity follows from knowing that idempotents in such a~semigroup~$S$ must commute,
and that $a^{-1}a$ are idempotents in~$S$ whenever $a\in S$~\cite{h}.
If~$S$ is an abelian semigroup then the identity map on~$S$ naturally provides a~$*$-operation, so, for example, the
semigroup of natural numbers may be considered as a~$*$-semigroup.
The semigroup of linear contractions on a~Hilbert space form a~$*$-semigroup (with the $0$ operator an absorbing element
for this semigroup), where the Hilbert space adjoint of a~bounded linear map is the $*$-operation.
This $*$-semigroup provides the primary focus for representations of $*$-semigroups considered below.
It is not a~priori clear that there are easily available nontrivial $*$-representations of~$S$ or $l_{1}(S)$ as bounded
operators on a~Hilbert space.
For~$S$ without an absorbing zero element, and nonempty there is always at least one nonzero $*$-semigroup homomorphism
of~$S$ to the contractions on a~Hilbert space $\mathcal{H}$, namely the $*$-homomorphism mapping~$S$ to the identity
element of $B\left(\mathcal{H} \right)$.

For a~$*$-semigroup~$S$ def\/ine a~universal $C^*$-algebra $C_{b}^{\ast}(S)$ as a~$C^*$-algebra $C_{b}^{\ast}(S)$ along
with a~$*$-semigroup homomorphism $\iota:S\rightarrow C_{b}^{\ast}(S)$ satisfying a~universal property: if
$\gamma:S\rightarrow B$ is a~$*$-semigroup homomorphism to a~$C^*$-algebra~$B$ then there is a~unique $*$-homomorphism
$\pi_{\gamma}=\pi:C_{b}^{\ast}(S)\rightarrow B$ such that $\pi_{\gamma}\circ\iota=\gamma$.
That such a~universal $C^*$-algebra of a~given semigroup exists would follow if there are $*$-semigroup homomorphisms
from~$S$ to a~$C^*$-algebra, and if
\begin{gather*}
\sup\{\Vert \beta(s)\Vert \,|\,\beta~\text{is a $*$-representation of}~S~\text{in a~$C^*$-algebra}\}
\end{gather*}
is f\/inite for each $s\in S$.
If~$S$ is described by generators and relations and this $\sup$ is bounded for each generator then existence follows.
That the map $\iota$ is a~monomorphism would follow if, for example, there is a~$*$-monomorphism of~$S$ into some
$C^*$-algebra.
In Def\/inition~\ref{universal} below we consider a~more restrictive universal $C^*$-algebra associated with~$S$.
For the main $*$-semigroups we work with the two def\/initions coincide.

The contractions on a~Hilbert space form a~natural $*$-semigroup of bounded operators under composition.
It is of some interest, using the group situation as a~guide, to restrict attention to the situation where elements
of~$S$ are all mapped to (bounded) linear operators that are uniformly bounded in norm by~$1$, i.e., to contractions.
Call these ($*$-)semigroup representations contractive ($*$-)representations of a~($*$-)semigroup.
These are the representations of a~$*$-semigroup that are our main concern below.
Since there are subsemigroups of $B\left(\mathcal{H}\right)$ which may have a~unit other than the identity operator, one
may want to also direct attention to those representations of~$S$ that are nondegenerate.
If the semigroup~$S$ is unital, then a~nondegenerate $*$-representation would necessarily map the unit to the unit of
$B\left(\mathcal{H}\right)$.
A~$*$-representation of a~$*$-semigroup~$S$ in the contractions of a~$C^*$-algebra necessarily gives rise to a~bounded
$*$-representation of the Banach $*$-algebra $l_{1}(S)$ in this $C^*$-algebra.
Note that in an inverse semigroup~$S$, $a^{-1}a$ is a~selfadjoint idempotent for any~$a$ in~$S$.
Consequently any $*$-representation of an inverse semigroup~$S$ must always occur via partial isometries, so in
particular is always a~representation by contractions.
In fact, for any $*$-semigroup~$S$ generated by a~set of elements~$a$ satisfying $aa^{\ast}a=a$ (or
$a^{\ast}aa^{\ast}=a^{\ast}$) any $*$-representation of~$S$ is via contractions.

More generally, for a~specif\/ied $*$-semigroup~$S$ and family~$F$ of $*$-homomorphisms from~$S$ to $C^*$-algebras, one
may consider a~universal $C^*$-algebra $C_{F}^{\ast}(S)$, paired with a~map $\iota:S\rightarrow C_{F}^{\ast}(S)$ in~$F$,
associated with~$S$.
For example, $F$ could consist of $*$-representations, or contractive \mbox{$*$-representations}, or as described in subsequent
sections, certain order $*$-representations.
Note that the universal property implies that the image of the map $\iota:S\rightarrow C_{F}^{\ast}(S)$ generates the
$C^*$-algebra $C_{F}^{\ast}(S)$, or in other words the $*$-subalgebra $\mathbb{C}[\iota(S)]$ of $C_{F}^{\ast}(S)$ is
a~dense $*$-subalgebra of $C_{F}^{\ast}(S)$.
Although a~standard argument we include it for completeness.

\begin{Proposition}
\label{basic}
For~$S$ a~$*$-semigroup, assume the universal $C^*$-algebra $C_{F}^{\ast}(S)$ exists.
The $*$-subalgebra $\mathbb{C}[\iota(S)]$ of $C_{F}^{\ast}(S)$ is a~dense $*$-subalgebra of $C_{F}^{\ast}(S)$.
\end{Proposition}

\begin{proof}
Let~$C$ be the closure of the $*$-algebra $\mathbb{C}[\iota(S)]$ and $j:C\rightarrow C_{F}^{\ast}(S)$ the canonical
inclusion, so $j\circ \underline{\iota}=\iota$, where $\underline{\iota}$ is the restriction of $\iota$ to codomain~$C$.
We show that $(C,\underline{\iota})$ satisf\/ies the universal property.
If $\varphi:S\rightarrow B$ is a~$*$-semigroup homomorphism in~$F$ to a~$C^*$-algebra~$B$ then set $\pi_{C}=\pi\circ j$,
where $\pi:C_{F}^{\ast}(S)\rightarrow B$ is the unique $*$-homomorphism such that $\pi\circ\iota=\varphi$.
Clearly $\pi_{C}\circ\underline{\iota}=\pi\circ j\circ\underline{\iota}=\varphi$, so it only remains to show that
$\pi_{C}$ is the unique such map.
The universal property for $C_{F}^{\ast}(S)$ implies there is a~$*$-homomorphism $\rho:C_{F}^{\ast}(S)\rightarrow C$
with $\rho \circ\iota=\underline{\iota}$, so $\rho\circ j\circ\underline{\iota}=\underline{\iota}$.
Since $\underline{\iota}(S)$ generates~$C$ as a~$C^*$-algebra, $\rho\circ j={\rm Id}_{C}$.
If $\pi_{C}^{\prime}$ is another map with $\pi_{C}^{\prime}\circ\iota=\varphi$ then $\pi_{C}^{\prime}\circ\rho\circ
\iota=\pi_{C}^{\prime}\circ\underline{\iota}=\varphi=\pi\circ\iota$, so the universal property for $C_{F}^{\ast}(S)$
yields $\pi_{C}^{\prime}\circ\rho =\pi$.
Thus $\pi_{C}^{\prime}=\pi_{C}^{\prime}\circ\rho\circ j=\pi\circ j=\pi_{C}$.
\end{proof}

The basic universal $C^*$-algebra for a~$*$-semigroup under consideration is that associated with contractive
$*$-representations.

\begin{Definition}
\label{universal}
The universal $C^*$-algebra $C^{\ast}(S)$ for a~$*$-semigroup~$S$ is a~$C^*$-algebra $C^{\ast}(S)$ along with
a~$*$-semigroup homomorphism $\iota:S\rightarrow C^{\ast}(S)$ into the contractions of $C^{\ast}(S)$ satisfying the
following universal property:

if $\gamma:S\rightarrow B$ is a~contractive $*$-semigroup homomorphism to
a~$C^*$-algebra~$B$ then there is a~unique $*$-homomorphism $\pi_{\gamma}=\pi:C^{\ast}(S)\rightarrow B$ such that
$\pi_{\gamma}\circ\iota=\gamma$.

The homomorphism $\iota$ is denoted $\iota_{S}$ when necessary.
\end{Definition}

\subsection{Examples}

\begin{Example}[groups]If a~group~$G$ is viewed as a~unital $*$-semigroup, where the $*$ of an element is def\/ined to be its inverse, then a~unital $*$-representation in a~$C^*$-algebra must be a~unitary representation, so a~unital
$*$-representation.
Conversely, a~unitary representation of~$G$ is clearly a~unital $*$-representation.
Thus the universal $C^*$-algebra $C^{\ast}(G)$ for the group~$G$ as a~$*$-semigroup exists and is the usual group
$C^*$-algebra of~$G$.
\end{Example}

\begin{Example}[abelian semigroups]For a~given abelian $*$-semigroup~$S$ with $*$ equal to the identity map on~$S$, if the universal
$C^*$-algebra $C^{\ast}(S)$ exists, then each element of~$S$ is viewed as a~selfadjoint contraction in the abelian
algebra $C^{\ast}(S)$.
For example $\mathbb{N}$, the nonunital abelian semigroup generated by a~single element, may also be viewed as the
$*$-semigroup generated by a~single selfadjoint element.
The $C^*$-algebra $C^{\ast} (\mathbb{N})$ is therefore the universal $C^*$-algebra generated by a~selfadjoint
contraction, namely the abelian nonunital $C^*$-algebra of continuous $\mathbb{C}$-valued functions $\{f\in
C([-1,1])\,|\, f(0)=0\}$.
For the unital $*$-semigroup $\mathbb{N}_{0}$ (here $0$ denotes the additive unit of $\mathbb{N}_{0}$),
$C^{\ast}(\mathbb{N}_{0})$ is the abelian unital $C^*$-algebra $C([-1,1])$.
\end{Example}

Let~$B$ be a~set with an involution $\sigma$.
The free $*$-semigroup $F_{B}$ on~$B$ is a~$*$-semigroup, with a~$*$-map $\iota:B\rightarrow F_{B}$, satisfying
the following: if $\gamma:B\rightarrow T$ is a~$*$-map to a~$*$-semigroup~$T$ then there is a~unique $*$-semigroup
homomorphism $\pi_{\gamma}:F_{B}\rightarrow T$ with $\pi_{\gamma}\circ\iota=\gamma$.
The universal $*$-semigroup $F_{B}$ exists and is unique up to a~bijective $*$-homomorphism.
Given a~set~$C$ and a~copy $\widetilde{C}$ of~$C$ there is an involution $\sigma$ on the disjoint union
$B=C\sqcup\widetilde{C}$; we may view $F_{B}$ as the universal $*$-semigroup generated by~$C$.
The free inverse semigroup $I_{B}$ on a~set~$B$ with an involution can be similarly def\/ined.
It exists and is a~quotient $*$-semigroup of $F_{B}$.

\begin{Example}[inverse semigroups]
The symmetric inverse semigroup $\mathcal{J}(X)$ on a~set~$X$ consists of the partial bijections on
a~set~$X$.
Assume a~$*$-semigroup~$S$ can be embedded via a~$*$-mono\-morphism $\rho$ as a~$*$-subsemigroup of the inverse semigroup
$\mathcal{J}(X)$.
Then there is a~$*$-representation $\pi$ of~$S$ by bounded linear operators on the Hilbert space $l^{2}(X)$ given by
\begin{gather*}
(\pi(s)f)(x)=f\big(\rho(s)^{-1}x\big)\qquad\text{for}\qquad f\in l^{2}(X),\quad x\in X
\end{gather*}
(where~$f$ is interpreted to have the value zero whenever $\rho(s)^{-1}x$ is not def\/ined).
The Vagner--Preston representation theorem~\cite{h}
ensures that any inverse semigroup~$S$ has such an embed\-ding~$\rho$.
Thus the universal $C^*$-algebras $C_{b}^{\ast}(S)$ and $C^{\ast}(S)$ always exist for a~general inverse semigroup~$S$,
and are isomorphic by our previous remarks (cf.~\cite{dp}).
We remark that the selfadjoint idempotents $a^{-1}a$  $(a\in S)$ in an inverse semigroup form a~commutative subsemigroup
of idempotents.
If this subsemigroup is totally ordered, we are led to aspects of a~structure needed to consider nest algebras
associated with~$S$.
\end{Example}

\begin{Example}[free $*$-semigroup with one generator]\label{free sgrp}
Recall that $\mathbb{N}^{+}$ and $\mathbb{N}^{-}$ denote the abelian semigroups of strictly positive and strictly
negative integers respectively.
Def\/ine $A_{c}$ to be the free semigroup product of $\mathbb{N}^{+}$ and $\mathbb{N}^{-}$.
Form the set of nonempty reduced words, or strings
\begin{gather*}
\left\{(n_{0},n_{1},\dots,n_{k})\,|\, n_{i}\in\mathbb{N}^{+}  \Leftrightarrow n_{i+1}\in\mathbb{N}^{-}\right\}
\end{gather*}
of f\/inite alternating sequences of elements in $\mathbb{N}^{+}$ and $\mathbb{N}^{-}$, with multiplication given by
concatenation of words followed by reduction to a~reduced word.
A~$*$-operation on $A_{c}$ is given~by
\begin{gather*}
(n_{0},n_{1},\dots,n_{k})^{\ast}=(-n_{k},-n_{k-1},\dots,-n_{0}).
\end{gather*}
It is clear that $A_{c}$ is the free $*$-semigroup on the two element set $B=\{1,-1\}$ with involution $1^{\ast}=-1$, so
the free $*$-semigroup generated by a~single element.
Its universal $C^*$-algebra $C^{\ast}(A_{c})$ is the universal $C^*$-algebra generated by a~single contraction.
\end{Example}

Consider the free inverse semigroup generated by a~single element~$a$.
The universal $C^*$-algebra of this inverse semigroup is the universal $C^*$-algebra of a~partial isometry~$v$ with the
property that all words in~$v$ and $v^{\ast}$ are partial isometries.
This is equivalent to~$v$ being a~power partial isometry~\cite{cdp}.
Note that the universal $C^*$-algebra generated by an isometry is $C^{\ast}(B)$ where~$B$ is the bicyclic semigroup,
a~unital inverse semigroup~\cite{cdp, h}.

\begin{Example}[universal abelian $*$-semigroup with one generator] The universal abelian $*$-semi\-group
generated by a~single element is
$(\mathbb{N}_{0}^{+}\times \mathbb{N}_{0}^{-})\backslash\{(0^{+},0^{-})\}$, where $(n,-m)(\widetilde
{n},-\widetilde{m})=(n+\widetilde{n},-m-\widetilde{m})$ and $(n,-m)^{\ast}=(m,-n)$.
The universal $C^*$-algebra of this $*$-semigroup, for contractive $*$-representations, is the universal $C^*$-algebra
generated by a~normal contraction, so $\{f\in C(\overline{\mathbb{D}})\,|\, f(0,0)=0\}$, where $\overline{\mathbb{D}}$ is
the closed unit disc in the plane.
The universal unital abelian $*$-semigroup generated by a~single element is $\mathbb{N}_{0}^{+}\times \mathbb{N}_{0}^{-}$,
with unit $(0^{+},0^{-})$, and its universal $C^*$-algebra
$C^{\ast}(\mathbb{N}_{0}^{+}\times\mathbb{N}_{0}^{-})$ is $C(\overline{\mathbb{D}})$.
\end{Example}

Let $\beta:S\rightarrow T$ be a~$*$-homomorphism of $*$-semigroups, where the universal $C^*$-algebras $C^{\ast}(S)$ and
$C^{\ast}(T)$ for~$S$ and~$T$ exist.
The universal property applied to the map $\beta\circ\iota_{T}$ implies that there is a~$*$-homomorphism of
$C^*$-algebras $\beta_{\ast}:C^{\ast}(S)\rightarrow C^{\ast}(T)$.
If $\gamma:R\rightarrow S$ is another such $*$-homomorphism then a~standard universal property argument implies
$(\beta\gamma)_{\ast} =\beta_{\ast}\gamma_{\ast}$.
If $\beta:S\rightarrow T$ is a~surjection of $*$-semigroups and the universal $C^*$-algebras for~$S$ and~$T$ exist then
Proposition~\ref{basic} implies $\beta_{\ast}:C^{\ast}(S)\rightarrow C^{\ast}(T)$ is a~surjection of $C^*$-algebras.
There are also appropriate versions of these types of properties for various families~$F$ of $*$-homomorphisms.

For example the group $\mathbb{Z}$ viewed as a~$*$-semigroup is a~unital abelian $*$-semigroup generated by the element
$1$.
The universal property applied to the universal unital abelian $*$-semigroup $\mathbb{N}_{0}^{+}\times \mathbb{N}_{0}^{-}$
yields a~unital $*$-homomorphism $\tau: \mathbb{N}_{0}^{+}\times \mathbb{N}_{0}^{-}\rightarrow\mathbb{Z}$ mapping $(m,-n)\rightarrow m-n$.
Using the other generator $-1$ of $\mathbb{Z}$ yields the $*$-homomorphism $-\tau$.
Then $\tau_{\ast}$ is the restriction map, the surjective $*$-homomorphism $i^{\#}:C(\overline{\mathbb{D}})\rightarrow
C(\mathbb{T})$, where $i:\mathbb{T\rightarrow}\overline{\mathbb{D}}$ is the inclusion map.
With the $*$-operation
on $\mathbb{N}_{0}$ def\/ined by the identity map there is a~unital $*$-homomorphism
$j:\mathbb{N}_{0}^{+} \rightarrow\mathbb{N}_{0}^{+}\times \mathbb{N}_{0}^{-}$ mapping~$n$ to $(n,-n)$ (and $0$ to
$(0^{+},0^{-})$) in the $*$-subsemigroup $\tau^{-1}(0)$ of $\mathbb{N}_{0}^{+}\times \mathbb{N}_{0}^{-}$.
The resulting $*$-homomorphism of universal contractive $C^*$-algebras is $j_{\ast}:C([-1,1])\rightarrow
C(\overline{\mathbb{D}})$ which maps the function $f(t)=t$ to the function $g(z)=\left\vert z\right\vert^{2}$.

Using the alternative
$*$-structure on $\mathbb{Z}$ given by the identity map the $*$-semigroup $\mathbb{N}_{0}^{+}$
embeds into the $*$-semigroup $\mathbb{Z}$ via $n\rightarrow n$.
With this $*$-operation $C^{\ast}(\mathbb{Z})$ is $C(\{-1,1\})$, the universal $C^*$-algebra generated by an invertible
selfadjoint contraction (with contractive inverse).
The induced homomorphism $C([-1,1])\rightarrow C(\{-1,1\})$ of the continuous function $C^*$-algebras is the restriction
$*$-homomorphism.

\subsection[The $*$-semigroup~$A$ and the universal $C^*$-algebra generated by a~partial isometry]
{The $\boldsymbol{*}$-semigroup~$\boldsymbol{A}$ and the universal $\boldsymbol{C^*}$-algebra generated\\ by a~partial isometry}

Both past and ongoing approaches for associating $C^*$-algebras with semigroups often involve representing semigroups as
isometries on a~Hilbert space.
This approach implicitly involves viewing the semigroup as a~subsemigroup of a~$*$-semigroup, since any semigroup of
isometries on a~Hilbert space naturally sits inside the $*$-subsemigroup generated by the representing isometries and
their adjoints.
More generally any representation of a~semigroup as operators on a~Hilbert space is necessarily a~map of the semigroup
into a~$*$-semigroup.
Since embeddings of semigroups in other semigroups may negate semigroup structure properties, the naturality properties
of such embeddings could prof\/itably be acknowledged and initially explored at the semigroup level.

An illustration of this implicit approach occurs in~\cite{li}, where a~full $C^*$-algebra for a~left cancellative unital
semigroup~$P$ is def\/ined by f\/irst embedding~$P$ in a~particular inverse semigroup, so a~$*$-semigroup, $S_{P}$.
Left translation by the elements of~$P$ on the set~$P$ allows~$P$ to be embedded as a~subsemigroup of the symmetric
inverse semigroup $\mathcal{J}(P)$ of partial bijections on the set~$P$.
The full $C^*$-algebra of~$P$ is, in ef\/fect, def\/ined as the universal $C^*$-algebra of the $*$-subsemi\-group~$S_{P}$ of
$\mathcal{J}(P)$generated by~$P$.
The structure of~$S_{P}$, unital and with an absorbing zero element, then ensures that the elements of~$P$, which must
be partial isometries since~$\mathcal{J}(P)$ is an inverse semigroup, are actually isometries in this $C^*$-algebra.

We introduce a~particular $*$-semigroup~$A$ which is naturally associated with a~partial isometry.
Describe an equivalence relation $\sim$ on $A_{c}$ compatible with the $*$-semigroup operations, so if~$x$,~$y$ and~$z$ are
elements of $A_{c}$ with $x\sim y$ then $x^{\ast}\sim y^{\ast}$, $xz\sim yz$ and $zx\sim zy$, by stipulating that
\begin{gather*}
(n_{0},n_{1},\dots,n_{k})\sim(n_{0},n_{1},\dots,n_{i}\pm1+n_{i+2},\dots, n_{k})
\end{gather*}
whenever $0\leq i\leq k-2$, and $n_{i+1}=\pm1$.
For example if $n_{i+1}=1$ then $n_{i}$ and $n_{i+2}$ must both be in $\mathbb{N}^{-}$, and the sum $n_{i}+1+n_{i+2}$
(performed in $\mathbb{Z}$) is again in $\mathbb{N}^{-}$.

\begin{Definition}[semigroup~$A$] Let~$A$ be the $*$-semigroup of congruence classes $(A_{c})/\sim$.
A~reduced word in~$A$ is described using the same notation $n=(n_{0},n_{1},\dots,n_{k})$, where $n_{i+1}\neq\pm1$
whenever $0\leq i\leq k-2$.
Def\/ine the length of a~reduced word~$n$, denoted $l(n)$, to be $k+1$.
\end{Definition}

The $*$-semigroup~$A$ is not unital, does not have a~zero, and is not an inverse semigroup.
It is also not a~left (so also not right) cancellative semigroup; namely if $mn=mr$ for $m$, $n$, $r$ in~$A$,
then it does not follow that $n=r$.
For example, choose $m=(m_{0},\dots,m_{k})$ to be any (reduced) word with $m_{k}>1$,~$n$ to be a~(reduced) word
$(-1,2,n_{2},\dots,n_{k})$, and $r=(1,n_{2},\dots,n_{k})$.
In particular one could choose $n=(-1,2,-1)$ and $r=(1,-1)$.

Since the elements $(-1,1)$ and $(1,-1)$ are selfadjoint idempotents in~$A$ they must be mapped to projections under
any $*$-homomorphism into a~$C^*$-algebra, and the generator $(1)$ of~$A$ as a
$*$-semigroup is therefore always mapped
to a~contraction, in fact a~partial isometry.
Thus the semigroup~$A$ has the property that any $*$-representation is also a~contractive representation, so
$C_{b}^{\ast}(A)$ exists and is $*$-isomorphic to $C^{\ast}(A)$.

Denote by $\mathcal{P}$ the universal $C^*$-algebra~\cite{bn} generated by a~partial isometry $\upsilon$.

\begin{Theorem}
\label{first iso}
There is a~contractive $*$-semigroup homomorphism $\sigma:A\rightarrow\mathcal{P}$ defined by mapping the element $(1)$
of~$A$ to the partial isometry~$v$ of $\mathcal{P}$.
This yields a~$*$-isomorphism $\pi_{\sigma}:C^{\ast}(A)\rightarrow\mathcal{P}$ of $C^*$-algebras.
\end{Theorem}

\begin{proof}
Since $\upsilon^{\ast}\upsilon\upsilon^{\ast}=\upsilon^{\ast}$ and $\upsilon\upsilon^{\ast}\upsilon=\upsilon$ and $(1)$
is a~generator of~$A$ as a~$*$-semigroup, mapping the element $(1)$ of~$A$ to the partial isometry~$v$ of $\mathcal{P}$
def\/ines a~contractive $*$-semigroup homomorphism $\sigma$ of~$A$ into the $C^*$-algebra $\mathcal{P}$.
By the universal property there is a~$*$-homomorphism $\pi_{\sigma}:C^{\ast}(A)\rightarrow\mathcal{P}$ with
$\pi_{\sigma}\circ\iota=\sigma$.
Since $\mathcal{P}$ is generated as a~$C^*$-algebra by the partial isometry $\upsilon=\pi_{\sigma}\iota(1)$,
$\pi_{\sigma}$ is surjective.
Since $\iota(1)$ is a~partial isometry in the $C^*$-algebra $C^{\ast}(A)$, the universal property of the $C^*$-algebra
$\mathcal{P}$ implies there is a~$*$-homomorphism $\lambda:\mathcal{P} \rightarrow C^{\ast}(A)$ with
$\lambda(\upsilon)=\iota(1)$.
Since $\lambda \pi_{\sigma}\iota(1)=\iota(1)$ and $\iota(1)$ generates $\iota(A)$ as a~$*$-semigroup we have
$\lambda\pi_{\sigma}={\rm Id}_{\iota(A)}$, so by Proposition~\ref{basic} $\lambda\pi_{\sigma}={\rm Id}_{C^{\ast}(A)}$.
\end{proof}

Given that free products were used to describe~$A$ it is not surprising that~$A$ provides an example of a~discrete
$*$-semigroup whose universal $C^*$-algebra is not exact.
In fact $C^{\ast}(A)$ is Morita equivalent to the universal unital $C^*$-algebra generated by a~contraction, so to the
universal $C^*$-algebra for the (unital) $*$-semigroup generated by a~single element~\cite{bn}.

\section[$*$-Subsemigroups of~$A$]{$\boldsymbol{*}$-Subsemigroups of~$\boldsymbol{A}$}
\label{section2}

The $*$-semigroup~$A$ has many $*$-subsemigroups of interest.
For example, the words $n=(n_{0},n_{1}$, $\ldots,n_{k})$
in~$A$ with $k\geq1$, $n_{0}\in\mathbb{N}^{-}$ and
$n_{k}\in\mathbb{N}^{+}$ form a~$*$-subsemigroup $A_{+}$.
Similarly def\/ine the semigroup $A_{-}$ to be the words $n=(n_{0},n_{1}, \dots,n_{k})$ in~$A$ which have $k\geq1$,
$n_{0}\in\mathbb{N}^{+}$ and $n_{k}\in\mathbb{N}^{-}$.
Note that $A_{+}$ and $A_{-}$ are disjoint unital $*$-semigroups, with $(-1,1)$ the unit of $A_{+}$ and $(1,-1)$ the
unit of $A_{-}$, and the length of any element of $A_{+}$ or $A_{-}$ must be even.

View $\mathbb{Z}$ as a~$*$-semigroup where the $*$-operation arises from the group structure.
Since $A_{c}$ is the free $*$-semigroup generated by a~single element $(1)$ (Example~\ref{free sgrp}) there is
a~$*$-semigroup homomorphism $\tau:A_{c}\rightarrow\mathbb{Z}$ uniquely determined by $\tau((1))=1$.
We have
\begin{gather*}
\tau(n)=\sum\limits_{i=0}^{k}n_{i}\qquad\text{for}\quad n=(n_{0},n_{1},\dots ,n_{k})\in A_{c}.
\end{gather*}
Since $\tau$ is constant on the appropriate congruence classes, $\tau$ determines a~$*$-homomorphism, also denoted
$\tau$, on~$A$.
Thus $\tau (nm)=\tau(n)+\tau(m)$, $\tau(n^{\ast})=-\tau(n)$, and $\tau(nm)=\tau(mn)$.
Set
\begin{gather*}
A^{0}:=\tau^{-1}(0)
\end{gather*}
a~$*$-subsemigroup of~$A$.
It is clear that reduced words in $A^{0}$ must have length at least $2$.
Using~$\tau$ it follows that any left or right unit of an element of~$A$, so in particular any idempotent, must be in $A^{0}$.
The selfadjoint elements of~$A$ are all contained in $A^{0}$.
Referring to Theorem~\ref{first iso} note that the canonical action $\gamma$ of the group $\mathbb{T}$ by
$*$-automorphisms of $\mathcal{P}$, where $\gamma_{t} (\upsilon)=t\upsilon$ for $t\in\mathbb{T}$, determines
a~conditional expectation of $\mathcal{P}$ to a~$C^*$-subalgebra $\mathcal{P}^{0}$.
It is clear that the contractive representation $\sigma$ maps the $*$-semigroup $A^{0}$ to $\mathcal{P}^{0}$.

Consideration of the length function yields the following properties for multiplication of reduced words
$m=(m_{0},m_{1},\dots,m_{k})$ and $n=(n_{0},n_{1},\dots,n_{j})$ in~$A$: $l(mn)=l(m)+l(n)-1$ if\/f $m_{k}$ and $n_{0}$ have
the same sign; if $m_{k}$ and $n_{0}$ have opposite signs then $l(mn)$ is either $l(m)+l(n)$ or $l(m)+l(n)-2$.
The later occurs if\/f at least one of $\left\vert m_{k}\right\vert =1$ and $l(m)\geq2$, or $\left\vert n_{0}\right\vert
=1$ and $l(n)\geq2$.
It is also straightforward to see that if $mn=pq$ with $l(m)=l(p)$ and $l(n)=l(q)$ for $m,n,p,q\in A$ then $m=p$ and $n=q$.
Using these observations, and $\tau$, it is straightforward to show that the only idempotents of~$A$ are the elements $(-1,1)$ and $(1,-1)$.

Other $*$-semigroups of $A^{0}$ may be identif\/ied by considering further structure.
First we def\/ine and identify certain irreducible elements of the $*$-semigroup $A^{0}$.
Two $*$-subsemigroups of $A^{0}$ are then introduced, which, although not satisfying unique factorization into
irreducibles, satisfy a~certain unique factorization property.

If~$p$ is an element in a~general semigroup~$S$ with $p=pn$ for some~$n$ then $p=pn=pn^{2}=\dots =pn^{k}$ for all
$k\in\mathbb{N}^{+}$.
It follows, using the properties of $\tau$ and the length~$l$ under multiplication, that an element~$n$ of~$A$ is
a~product~$nm$ (or~$mn$) in~$A$ if\/f~$m$ is one of two idempotents $(-1,1)$ or $(1,-1)$.
Although a~useful approach to irreducible elements for more general semigroups may dif\/fer from the one below, for our
specif\/ic context the following approach is suf\/f\/icient.

\begin{Definition}[irreducible]
An element~$p$ of a~semigroup~$S$ is irreducible if~$p$ cannot be written as a~product of (at least) two elements
in~$S$, or whenever $p=mn$ with $m,n\in S$ then $p=m$, or $p=n$.
Denote the set of irreducible elements of~$S$ by
$\operatorname*{Irr}(S)$.
Say~$p$ factors nontrivially in~$S$ if it is not irreducible, so there exist $m,n\in S$ with $p=mn$, and neither $p=m$
nor $p=n$.
\end{Definition}

If~$S$ is a~subsemigroup of a~semigroup~$T$, and $d\in S$, then~$d$ is irreducible in~$T$ implies~$d$ is irreducible
in~$S$, so $\operatorname*{Irr} (S)\supseteq\operatorname*{Irr}(T)\cap S$.
If~$S$ is a~$*$-semigroup then an element~$p$ of~$S$ is irreducible if and only if $p^{\ast}$ is irreducible in~$S$.

In our situation if~$p$ is irreducible in any subsemigroup of $A^{0}$ containing either of the idempotents $(-1,1)$ and
$(1,-1)$ then $p=mn$ implies that either $p=m$ (so necessarily~$n$ is one of these idempotents) or $p=n$ (with~$m$ one
of the same two idempotents).
It also follows that the idempotents of~$A$ are irreducible in any subsemigroup of $A^{0}$ containing either of them.

To investigate the irreducible elements of $A^{0}$ we introduce a~sequence of maps
\begin{gather*}
\sigma_{r}: \ A\rightarrow\mathbb{Z},
\qquad
\sigma_{r}(n)=\sum\limits_{i=0}^{r}n_{i},
\qquad
\text{and}
\qquad
\sigma_{r}(n)=\tau(n)
\quad
\text{if}
\quad
k\leq r,
\end{gather*}
for $n=(n_{0},n_{1},\dots,n_{k})\in A$ and $r\in\mathbb{N}_{0}$.
The next lemma shows that irreducibility is characte\-ri\-zed by a~lack of sign changes in the sequence~$\sigma_{r}(n)$
as~$r$ increases from~$0$ through $k-1$.

\begin{Lemma}\label{fund lemma}
Let $p=(p_{0},\dots,p_{l})\in A^{0}$.
If $\sigma_{0} (p)\sigma_{r}(p)>0$ for all $r\in\left\{1,2,\dots,l-1\right\} $ then~$p$ is irreducible in the
$*$-semigroup $A^{0}$.
If $p\in\operatorname*{Irr}(A^{0})$ and~$p$ is in reduced form, so $p_{i}\neq\pm1$ whenever $1\leq i\leq l-1$, then
$\sigma_{0}(p)\sigma_{r}(p)>0$ for all $r\in\left\{1,2,\dots,l-1\right\}$.
\end{Lemma}

\begin{proof}
We assume for the sake of def\/initeness that $p_{0}\in\mathbb{N}^{-}$.
The situation with $p_{0}\in\mathbb{N}^{+}$ is similar.

Suppose~$p$ factors nontrivially in $A^{0}$, so there are $m,n\in A^{0}$ with $p=mn$ and neither $p=m$ nor $p=n$.
We will show that $\sigma_{0}(p)\sigma_{r}(p)\leq0$ for some $r\in\left\{1,2,\dots,l-1\right\}$.
Let $m=(m_{0},\dots,m_{j})$ and $n=(n_{0},\dots,n_{k})$, with $m_{0}=p_{0}<0$.
We consider various cases.

First suppose $\sign(m_{j})=\sign(n_{0})$.
Then $m_{j}+n_{0}$ cannot be $\pm1$ and
\begin{gather*}
p=(m_{0},\dots,m_{j}+n_{0},\dots,n_{k}).
\end{gather*}
If $\sigma_{0}(p)\sigma_{r}(p)>0$ for all $r\in\left\{1,2,\dots,j-1\right\} $ then (since $p_{0}\in\mathbb{N}^{-}$) we
must have that $\sigma_{r}(p)<0$ for all $r\in\left\{1,2,\dots,j-1\right\}$.
Since $0=\tau(m)=\sigma_{j-1}(p)+m_{j}$ it follows that $m_{j}>0$.
Thus $n_{0}$ is also strictly larger than $0$, and we have $\sigma_{j}(p)=\sigma_{j}(m)+n_{0}>\tau(m)=0$.
Therefore $\sigma_{0}(p)\sigma_{j}(p)<0$, and therefore either $\sigma_{0}(p)\sigma_{r}(p)\leq0$ for some
$r\in\left\{1,2,\dots,j-1\right\} $ or $\sigma_{0}(p)\sigma_{j}(p)<0$.

Next suppose $\sign(m_{j})\neq \sign(n_{0})$.
If neither $m_{j}$ nor $n_{0} =\pm1$, then $p=(m_{0},\dots,m_{j},n_{0},\dots$, $n_{k})$ so $\sigma_{j} (p)=\tau(m)=0$ and
$\sigma_{0}(p)\sigma_{j}(p)\leq0$.
Next consider if $m_{j}=\pm1$.
If for example $m_{j}=1$, then $m_{j-1}$ and $n_{0}$ are both less than $0$, and $p_{j-1}=m_{j-1}+1+n_{0}$ (or
$p_{j-1}=m_{j-1}$).
If $m_{j-1}=-1$ then either $m=(-1,1)$ and the factorization of~$p$ is not nontrivial, or $m=(m_{0},\dots,m_{j-2},-1,1)$
so $0=\tau(m)=\sigma_{r-2}(p)$ and $\sigma_{0}(p)\sigma_{j-2}(p)\leq0$.
Therefore we may assume $m_{j-1}\leq-2$.
Since $\tau(m)=0$ and $m_{0}=p_{0}<0$, necessarily $j\geq3$.
Since
\begin{gather*}
0=\tau(m)=m_{0}+\dots +m_{j-2}+(m_{j-1}+1)<m_{0}+\dots +m_{j-2}=\sigma_{j-2}(p),
\end{gather*}
we have $\sigma_{j-2}(p)>0$, and the product $\sigma_{0}(p)\sigma_{j-2}(p)<0$.
On the other hand, if $m_{j}=-1$, then as before $m_{j-1}$ and $n_{0}$ are both greater than $0$ and
$p_{j-1}=m_{j-1}-1+n_{0}$.
We have
\begin{gather*}
\sigma_{j-1}(p)=m_{0}+\dots + m_{j-2}+(m_{j-1}-1+n_{0})=\tau(m)+n_{0}>0,
\end{gather*}
so $\sigma_{0}(p)\sigma_{j-1}(p)<0$.
The cases where $n_{0}=\pm1$ are dealt with similarly.

To prove the condition is necessary, suppose that $\sigma_{0}(p)\sigma_{r}(p)\leq0$ for some
$r\in\left\{1,2,\dots,l-1\right\}$.
Since $p_{0} \in\mathbb{N}^{-}$, we have $\sigma_{r}(p)\geq0$ for this $r\in\left\{1,2,\dots,l-1\right\}$.
Thus $r\geq1$ and $r\leq l-1$ so $l\geq2$, and~$p$ is a~word of length at least $3$.
Choose~$r$ to be the least such element.
Thus $\sigma_{r-1}(p)<0$, so we have that
\begin{gather*}
p_{r}>0
\qquad
\text{and}
\qquad
\sum\limits_{i=r}^{l}p_{i}=\tau(p)-\sigma_{r-1}(p)=0-\sigma_{r-1}(p)>0.
\end{gather*}
Consider the two possibilities $\sigma_{r}(p)=0$ or $0<\sigma_{r}(p)$.
If $\sigma_{r}(p)=0$ (so $r<l-1$, as otherwise $p_{l}=0$) then set
\begin{gather*}
m=(p_{0},\dots,p_{r-1},p_{r})
\qquad
\text{and}
\qquad
n=(p_{r+1},\dots,p_{l}).
\end{gather*}
If $0<\sigma_{r}(p)$  ($=p_{r}+\sigma_{r-1}(p)$) set
\begin{gather*}
m=(p_{0},\dots,p_{r-1},-\sigma_{r-1}(p))
\qquad
\text{and}
\qquad
n=(p_{r}+\sigma_{r-1}(p),p_{r+1},\dots,p_{l}).
\end{gather*}
In either case $m,n\in A^{0}$ and $p=mn$.
To f\/inish we need to ensure that this is not a~trivial factorization, namely that $p=m$ or $p=n$ is not the case.
Since~$p$ is assumed irreducible, either $p=m$ and so $n=(-1,1)$, or $p=n$ and so $m=(-1,1)$.
In the case $\sigma_{r}(p)=0$ this means either $p_{r+1}=-1$ or $p_{r}=1$, contradicting that~$p$ is in reduced form.
If $0<\sigma_{r}(p)$ then $-\sigma_{r-1}(p)$ and $p_{r}+\sigma_{r-1}(p)$ are both positive.
Since $p=mn$, and $p_{r}=m_{r}+n_{0}$, where $r\geq1$,
 neither $p=m$ nor $p=n$ is possible.
\end{proof}

It follows that the only irreducible elements $p=(p_{0},\dots,p_{l})$ of $A^{0}$ with $\left\vert p_{0}\right\vert =1$,
or $\left\vert p_{l}\right\vert =1$, are the idempotents $(-1,1)$ and $(1,-1)$.
Also all elements of $A^{0}$ of length $2$, so $(-n,n)$ and $(n,-n)$ of $A^{0}$, $n\in\mathbb{N}$ are irreducible.
The irreducible idempotents $(-1,1)$ and $(1,-1)$ are units for various elements of $A^{0}$, so $A^{0}$ does not have
unique factorization into irreducibles.
For example
\begin{gather*}
(-n,n)(m,-m)=(-n,n)(-1,1)(m,-m)=(-n,n)(1,-1)(m,-m).
\end{gather*}
However, this describes all that can fail with unique factorization.

\begin{Theorem}\label{! decomp}
Any element of $A^{0}$ may be written as a~product of irreducible elements of $A^{0}$.
If this decomposition is minimal, meaning only single powers of the idempotent irreducibles appear, and an idempotent
does not appear if it acts as either a~left or right unit for the element to its right or left, then this decomposition
is unique.
\end{Theorem}

\begin{proof}
Let $p=(p_{0},\dots,p_{l})\in A^{0}$.
Choose the least $r\in\left\{1,2,\dots,l-1\right\} $ with $\sigma_{0}(p)\sigma_{r}(p)\leq0$.
If there is no such~$r$ then~$p$ is irreducible by Lemma~\ref{fund lemma}.
If there is such an~$r$ then as in the last paragraph of Lemma~\ref{fund lemma} $p=mn$, where $m,n\in A^{0}$ and
$m=(p_{0},\dots,p_{r-1},-\sigma_{r-1}(p))$ and $n=(p_{r} +\sigma_{r-1}(p),p_{r+1},\dots,p_{l})$ (if
$p_{r}+\sigma_{r-1}(p)=0$ then $n=(p_{r+1},\dots,p_{l})$).
The element~$m$ must be irreducible, since $\sigma_{0}(m)\sigma_{s}(m)=\sigma_{0}(p)\sigma_{s}(p)>0$ for all
$s\in\left\{1,2,\dots,r-1\right\}$.
Now apply the same argument to $n\in A^{0}$.

To see uniqueness of a~minimal decomposition, assume $p=(p_{0},\dots,p_{k})\in A^{0}$ is a~f\/inite pro\-duct
$m_{1}m_{2}\cdots m_{k}$ of irreducible elements $m_{j}=(m_{j0},\dots,m_{jr_{j}})$ of $A^{0}$ in the described minimal
way.
With this product write~$p$ as a~word $(m_{00},\dots,m_{kr_{k}})\in A^{0}$ (where the signs alternate).
Note that the f\/irst and last terms must be $m_{00}$ and $m_{kr_{k}}$.
The minimality condition in the given decomposition implies that this word must be in reduced form.
Applying the process described in the f\/irst paragraph one obtains the original minimal decomposition.
\end{proof}

Intersections of $*$-semigroups are $*$-semigroups.

\begin{Definition}
Def\/ine $*$-subsemigroups $A_{+}^{0}=A^{0}\cap A_{+}$ and $A_{-}^{0}=A^{0}\cap A_{-}$.
\end{Definition}

Since the unit $(-1,1)$ of $A_{+}$ is in $A^{0}$, and the unit $(1,-1)$ of $A_{-}$ is in $A^{0}$, each of the
$*$-subsemigroups $A_{+}^{0}$ and $A_{-}^{0}$ are unital, with units $(-1,1)$ and $(1,-1)$ respectively.
Many of the results shown below for $A_{+}^{0}$ have corresponding statements holding for $A_{-}^{0}$.
Since elements of either $A_{+}$ or $A_{-}$ have even length, elements of either $A_{\pm}^{0}$ must have even length.
The following result shows that these $*$-semigroups arise naturally.

\begin{Proposition}
\label{basic irr}
If $p\in A^{0}$ is irreducible in $A^{0}$ then $p\in A_{+}^{0}$ or $p\in A_{-}^{0}$, and~$p$ has even length.
\end{Proposition}

\begin{proof}
Let $p=(p_{0},\dots,p_{l})\in A^{0}$ be irreducible.
The preceding lemma shows $\sign(p_{0})=\sign(\sigma_{l-1}(p))$.
Since $0=\tau(p)$, $p_{l}=-\sigma_{l-1}(p)$, $\sign(p_{0})$ and $\sign(p_{l})$ are dif\/ferent.
Thus $p\in A_{+}^{0}$ or $p\in A_{-}^{0}$ and the length of~$p$ is even.
\end{proof}

Note that it is not the case that an element of $A_{+}^{0}$ that is irreducible as an element of the semigroup
$A_{+}^{0}$ is also irreducible in $A^{0}$.
For example $(-2,3,-3,2)$ is not irreducible in~$A^{0}$~-- it is the product $(-2,2)(1,-1)(-2,2)$ of irreducibles in
$A^{0}$~-- but is irreducible in the semigroup~$A_{+}^{0}$.

\begin{Definition}[semigroup $D_{0}$]\label{def}
Set $\operatorname*{Irr}\nolimits_{+}(A^{0}):=\operatorname*{Irr}(A^{0})\cap A_{+}$, the set of irreducible elements of~$A^{0}$
contained in~$A_{+}^{0}$.
Similarly for $\operatorname*{Irr}\nolimits_{-}(A^{0})$.
Def\/ine $D_{0}$ to be the subsemigroup of $A_{+}^{0}$ generated
by $\operatorname*{Irr}\nolimits_{+}(A^{0})$.
\end{Definition}

Since $A_{+}^{0}$ (and $A_{-}^{0}$) is $*$-closed, the previous proposition shows $\operatorname*{Irr}\nolimits_{+}(A^{0})$ (and
$\operatorname*{Irr}\nolimits_{-}(A^{0})$) is $*$-closed.
Since $(-1,1)$ is a~unit of $A_{+}^{0}$, and $\operatorname*{Irr}\nolimits_{+}(A^{0})$ is $*$-closed, $D_{0}$ is a~unital
$*$-semigroup.

\begin{Proposition}\label{D 1}
Let $p=(p_{0},\dots,p_{l})\in A^{0}$.
Then
\begin{gather*}
\operatorname*{Irr}\nolimits_{+}(A^{0})=\left\{p\in A^{0}\,|\,\sigma_{r}(p)<0~\text{for}~r\in\left\{1,2,\dots,l-1\right\}\right\},
\\
D_{0}=\left\{n\in A^{0}\,|\,\sigma_{r}(n)\leq0~\text{for all}~r\geq0\right\}.
\end{gather*}
\end{Proposition}

\begin{proof}
Lemma~\ref{fund lemma} yields the f\/irst statement.
The second follows from this f\/irst statement, the def\/inition of $D_{0}$, and Theorem~\ref{! decomp}.
\end{proof}

\begin{Proposition}
\label{irred D}
$\operatorname*{Irr}(D_{0})=\operatorname*{Irr}\nolimits_{+}(A^{0})$ and elements of $D_{0}$ have unique minimal decomposition
into irreducibles.
\end{Proposition}

\begin{proof}
Since $D_{0}$ is a~subsemigroup of $A^{0}$, an element of $D_{0}$ which is irreducible in $A^{0}$ must be an irreducible
element of $D_{0}$, so $\operatorname*{Irr}(D_{0})\supseteq\operatorname*{Irr}(A^{0})\cap D_{0}$.
This in turn $=\operatorname*{Irr}(A^{0})\cap A_{+}^{0}\cap D_{0} =\operatorname*{Irr}\nolimits_{+}(A^{0})$.

Conversely, if $d\in D_{0}$ is not the unit then by def\/inition~$d$ is a~product of elements in
$\operatorname*{Irr}\nolimits_{+}(A^{0})\backslash\{(-1,1)\}$.
By Theorem~\ref{! decomp} any element of $A^{0}$ may be written in a~unique (minimal) manner as a~product of elements in
$\operatorname*{Irr}(A^{0})$, so the given product must also be the unique product when~$d$ is viewed as an element of
$A^{0}$.
By the def\/inition of an irreducible element, if $d\in\operatorname*{Irr}(D_{0})$ then~$d$ must equal one of these
elements in $\operatorname*{Irr}\nolimits_{+}(A^{0})\backslash\{(-1,1)\}$.
\end{proof}

It therefore follows that nonunital elements of the semigroup $D_{0}$ have unique factorization into irreducible
elements of $D_{0}\backslash\{(-1,1)\}$.
The element $(-2,3,-3,2)$ of $A_{+}^{0}$ noted before Def\/inition~\ref{def} is therefore not in $D_{0}$, showing $D_{0}$
is a~proper subsemigroup of $A_{+}^{0}$.

\begin{Remark}
\label{free}
The $*$-semigroup $D_{0}$ is the free unital $*$-semigroup over the $*$-closed set\linebreak
$\operatorname*{Irr}\nolimits_{+}(A^{0})\backslash(-1,1)$.
\end{Remark}

\begin{Definition}[semigroup $D_{1}$]Let $D_{1}$ denote the $*$-subsemigroup of $A_{+}^{0}$ generated by the idempotent $(1,-1)$ and the
irreducible elements $\operatorname*{Irr}\nolimits_{+}(A^{0})$, so by the idempotent $(1,-1)$ and the $*$-semigroup $D_{0}$.
\end{Definition}

The same type of example before Theorem~\ref{! decomp} shows $D_{1}$ is not a~cancellative semigroup.
On the other hand the unital $*$-semigroup $D_{0}$ is cancellative.

\subsection[Canonical $*$-mappings of~$A$]{Canonical $\boldsymbol{*}$-mappings of~$\boldsymbol{A}$}
\label{section2.1}

We remark that for a~general $*$-semigroup~$S$, if $a\in S$ is an element with $aa^{\ast}a=a$ then $a^{\ast}a$ and
$aa^{\ast}$ are idempotents and the map $\alpha_{a}:s\rightarrow a^{\ast}sa$ def\/ines a~$*$-isomorphism (surjective
$*$-monomorphism) of $*$-subsemigroups $aa^{\ast}Saa^{\ast} \rightarrow a^{\ast}aSa^{\ast}a$, with inverse map
$\alpha_{a^{\ast}}$ restricted to $a^{\ast}aSa^{\ast}a$.

\begin{Definition}
Def\/ine $*$-maps $\alpha:A\rightarrow A$ and $\omega:A\rightarrow A$ by
\begin{gather*}
\alpha(n)=(-1)n(1),
\qquad
\omega(n)=(1)n(-1)
\end{gather*}
for $n\in A$.
\end{Definition}

These are the $*$-maps `implemented' by the element $(1)\!\in\! A$.
We have $\omega\circ\alpha(a)\!=\!(1,-1)a(1,-1)$ while $\alpha\circ\omega (a)=(-1,1)a(-1,1)$, $a\in A$.
Since $\tau\circ\alpha=\tau$ and $\tau \circ\omega=\tau$, both maps restrict to maps of $A^{0}$ to itself.
By the preceding comment, $\alpha$ therefore restricts to a~$*$-isomorphism
\begin{gather*}
\alpha: \ (1,-1)A^{0}(1,-1)\rightarrow(-1,1)A^{0}(-1,1)
\end{gather*}
of $*$-semigroups with inverse $\omega$ restricted to the $*$-semigroup $(-1,1)A^{0}(-1,1)$.
We remark that
\begin{gather*}
\alpha(1,-1)=(-1,1)
\qquad
\text{and}
\qquad
\omega(-1,1)=(1,-1).
\end{gather*}

The map $\alpha$ also maps $A_{+}$ to itself, and therefore maps $A_{+}^{0}$ to itself, though it is no longer `inner'
in these subsemigroups.
However $\omega(-1,1)=(1,-1)$, so $\omega$ does not restrict to a~map of $A_{+}^{0}$ to itself.
Since $A_{+}^{0}\subset(-1,1)A^{0}(-1,1)$, $\omega$ is a~$*$-homomorphism on $A_{+}^{0}$, so is a~$*$-homomorphism on
$D_{0.}$.
We have $(-1,1)n(-1,1)=n$ for $n\in A_{+}^{0}$, so $\alpha(\omega(n))=n$ for $n\in A_{+}^{0}$.
In particular $\alpha\circ\omega$ is the identity map on $D_{0}$.

Note that the $*$-map $\alpha$ on the $*$-semigroup~$A$ may be viewed as the restriction to $\sigma(A)$
(Theorem~\ref{first iso}) of the completely positive linear map on $\mathcal{P}$ given by $\varphi:x\rightarrow
v^{\ast}xv$.
This is used later to construct a~representation of a~certain $C^*$-correspondence.

\begin{Definition}
For $n=(n_{0},n_{1},\dots,n_{k})\in\operatorname*{Irr}\nolimits_{+}(A^{0})\backslash (-1,1)$ def\/ine
\begin{gather*}
\beta_{\omega}(n)=(n_{0}+1,n_{1},\dots,n_{k}-1).
\end{gather*}
\end{Definition}

For $n\in\operatorname*{Irr}\nolimits_{+}(A^{0})\backslash(-1,1)$, $n_{0}<-1$ (and $n_{k}>1$) so $\beta_{\omega}(n)\in
A_{+}^{0}$.
Therefore
\begin{gather*}
\omega(n)=(1,-1)\beta_{\omega}(n)(1,-1)
\end{gather*}
for $n\in\operatorname*{Irr}\nolimits_{+}(A^{0})\backslash(-1,1)$.
This follows since the right side is
\begin{gather*}
(1,-1,n_{0}+1,n_{1},\dots,n_{k}-1,1,-1)=(1,n_{0},n_{1},\dots,n_{k},-1)
\end{gather*}
which is $\omega(n)$.
Also
\begin{gather*}
\alpha(\omega(n))=\alpha(\beta_{\omega}(n))
\end{gather*}
for $n\in\operatorname*{Irr}\nolimits_{+}(A^{0})\backslash(-1,1)$ since, by the preceding equality, the left side is
\begin{gather*}
\alpha((1,-1)\beta_{\omega}(n)(1,-1))=(-1,1,-1)\beta_{\omega}(n)(1,-1,1)
\end{gather*}
which is $\alpha(\beta_{\omega}(n))$.

\begin{Proposition}
\label{irred}
The map $\alpha$ restricts to a~$*$-map $\alpha:D_{0}\rightarrow D_{0}$ with image
$\alpha(D_{0})=\operatorname*{Irr}\nolimits_{+}(A^{0})\backslash (-1,1)$.
The map $\beta_{\omega}\circ\alpha$ is the identity map on $D_{0}$ and $\alpha\circ\beta_{\omega}$ is the identity map
on $\operatorname*{Irr}\nolimits_{+}(A^{0})\backslash(-1,1)$.
The image of $\beta_{\omega}$ is $D_{0}$.
\end{Proposition}

\begin{proof}
If $n\in\operatorname*{Irr}\nolimits_{+}(A^{0})\backslash(-1,1)$ then Proposition~\ref{D 1} implies (using a~reduced form
$(n_{0},n_{1},\dots$, $n_{k})$ for~$n$ if necessary) that $\sigma_{r}(n)<0$ for $r\in\left\{1,\dots,k-1\right\} $ and
therefore $\sigma_{r}(\beta_{\omega}(n))\leq0$ for $r\in\left\{1,\dots,k-1\right\}$.
Thus $\beta_{\omega}(n)\in D_{0}$ by Proposition~\ref{D 1}.

Let $n=(n_{0},n_{1},\dots,n_{k})\in D_{0}$.
Then $\sigma_{0}(n)<0$ and by Proposition~\ref{D 1} $\sigma_{r}(n)\leq0$ for $r\geq1$.
Therefore, for $r\in \{1,2,\dots,k-1 \}$,  $\sigma_{r}(\alpha(n))=\sigma_{r}(n)-1<0$, so Proposition~\ref{D 1}
implies $\alpha(D_{0})\subseteq \operatorname*{Irr}\nolimits_{+}(A^{0})$.
It is clear that $\alpha(n)=(n_{0} -1,n_{1},\dots,n_{k}+1)$ so $\beta_{\omega}\circ\alpha$ is the identity map on
$D_{0}$.
We have already seen $\alpha\circ\beta_{\omega}=\alpha\circ\omega$ on $\operatorname*{Irr}\nolimits_{+}(A^{0})\backslash(-1,1)$.
Now use that $\alpha\circ\omega$ is the identity map on $D_{0}$.
\end{proof}

\begin{Proposition}
The $*$-semigroup $D_{0}$ is the smallest subsemigroup of~$A$ containing the element $(-1,1)$ and invariant under the
map $\alpha$.
\end{Proposition}

\begin{Proposition}
The image $\omega(D_{0})\subseteq(1,-1)D_{1}(1,-1)$, so $\omega:D_{0} \rightarrow D_{1}$ is a~$*$-semigroup homomorphism
with $\alpha\circ\omega$ the identity on $D_{0}$; $\omega$ splits the $*$-map $\alpha$.
\end{Proposition}

\begin{proof}
For $n\in\operatorname*{Irr}\nolimits_{+}(A^{0})\backslash(-1,1)$ we saw that $\omega(n)=(1,-1)\beta_{\omega}(n)(1,-1)\in D_{1}$,
while $\omega (-1,1)$ $=(1,-1)\in D_{1}$.
\end{proof}

\begin{Remark}
\label{b}
Note $(-1,1)\alpha(a)=\alpha(a)(-1,1)=\alpha(a)$ and
\begin{gather*}
\alpha((1,-1)a)=\alpha(a(1,-1))=\alpha(a),
\qquad
\alpha(a(1,-1)b)=\alpha (a)\alpha(b)
\end{gather*}
for $a,b\in A$.
Since $\alpha(1,-1)=(-1,1)\in D_{0}$ we have that $\alpha$ may be viewed as a~surjective $*$-map
\begin{gather*}
\alpha: \ D_{1}\rightarrow D_{0}.
\end{gather*}
It is not a~semigroup homomorphism, nor is it one-to-one.
We also remark that $D_{1}$ is the smallest semigroup of~$A$ containing $(1,-1)$ and invariant under the map~$\alpha$.

There is a~useful alternative
characterization of $D_{1}$ that is similar to the corresponding
description of $D_{0}$ (Proposition~\ref{D 1}), namely
\begin{gather*}
D_{1}=\left\{n\in A^{0}\,|\,\sigma_{r}(n)\leq1~\text{for all}~r\geq0\right\}.
\end{gather*}
\end{Remark}

We may obtain the elements of $\operatorname*{Irr}(D_{0})$ recursively.
For $n=(n_{0},n_{1},\dots ,n_{k})\in A^{0}$ in reduced form consider the map $\tau
^{+}(n)=\sum\limits_{n_{i}\geq1}n_{i}\in\mathbb{N}^{+}$.
Any subset~$L$ of $A^{0}$ can be written as the disjoint union over $k\in\mathbb{N}^{+}$ of the sets $L(k)=\left\{n\in
L\,|\,\tau^{+}(n)=k\right\}$.
For $n\in A_{+}^{0}$ we have $\tau^{+}(\alpha(n))=\tau^{+}(n)+1$ (although this is not generally true on $A$;
consider $n=(1,-1)$).
The subset~$L$ we are interested in is $\operatorname*{Irr}\nolimits_{+}(A^{0})=\operatorname*{Irr}(D_{0})$.
If $n\in \operatorname*{Irr}(D_{0})\backslash(-1,1)$ then $\tau^{+}(\beta_{\omega}(n))=\tau^{+}(n)-1$.
It is straightforward to check that $\tau^{+} (mn)=\tau^{+}(m)+\tau^{+}(n)$ for~$m$, $n\in\operatorname*{Irr}\nolimits_{+}
(A^{0})\setminus(-1,1)$.

\begin{Theorem}
\label{gen}
The subset
\begin{gather*}
\operatorname*{Irr}(D_{0})(k):=\left\{m\in\operatorname*{Irr}(D_{0})\,|\, \tau^{+}(m)=k\right\}
\end{gather*}
of $D_{0}$ consists of $\alpha^{k_{0}}\left(\prod\limits_{i=1}^{r}m_{i}\right)$, where $m_{i}\in\operatorname*{Irr}(D_{0})(k_{i})$
and $k_{0}+\sum\limits_{i=1}^{r}k_{i}=k$, where $k_{0}>0$ if $k>1$, and $k_{i}>1$ if $r>1$.
\end{Theorem}

\begin{proof}
For $k=1$ the only element~$n$ of $\operatorname*{Irr}\nolimits_{+}(A^{0})(1)$ is $(-1,1)$.
Choose $n\in\operatorname*{Irr}\nolimits_{+}(A^{0})(k)$ with $k>1$.
Propositions~\ref{irred} and~\ref{irred D} show $\beta_{\omega}(n)$ is a~product of elements, say~$r$ of them,
$m_{i}\in\operatorname*{Irr}\nolimits_{+} (A^{0})$, each of which must have $\tau^{+}(m_{i})<k$.
If $r>1$ we must have (again by Proposition~\ref{irred D}) none of the $m_{i}$ equal the unit $(-1,1)$ of $A_{+}^{0}$,
so $k_{i}>1$ for all these~$i$.
However, if $r=1$ then $m_{1}=(-1,1)$ is certainly possible.
Since $\alpha\beta_{\omega}(n)=n$ we obtain the description of $\operatorname*{Irr}\nolimits_{+}(A^{0})(k)$.
The statement follows by induction.
\end{proof}

One can now iteratively list the irreducible elements of $D_{0}$.
For example, $\operatorname*{Irr}(D_{0})(k)=\{-k,k\}$
for $0\leq k\leq4$.
We obtain $\operatorname*{Irr}(D_{0})(5)=\left\{(-5,5),(-3,2,-2,3)\right\} $, while
\begin{gather*}
\operatorname*{Irr}(D_{0})(6)=\left\{(-6,6),(-4,2,-2,4),(-4,3,-2,3),(-3,2,-3,4)\right\}.
\end{gather*}

\section[Partial orders on $*$-semigroups]{Partial orders on $\boldsymbol{*}$-semigroups}
\label{section3}

There is a~partial order on the selfadjoint elements $D_{1}^{\rm sa}$ of the $*$-semigroup $D_{1}$, which
may be viewed as taking place in the more general context of the category of ordered $*$-semigroups $(S,\preceq)$ and
their order maps.
The Schwarz property for a~$*$-map from a~$*$-semigroup to an ordered $*$-semigroup is introduced, and the $*$-map
$\alpha:D_{1}\rightarrow D_{1}$, def\/ined in the previous section, is shown to be an order map with the Schwarz property.
A~universal $C^*$-algebra $C^{\ast}(S,\preceq)$ for a~general ordered $*$-semigroup $(S,\preceq)$ and its order
representations is introduced as an interim step towards considering, in the following section, the universal
$C^*$-algebra for a~more restrictive matricial order on $*$-semigroups.
There are also some technical factorization results which have later use.

\begin{Definition}[ordered $*$-semigroup]An ordered $*$-semigroup is a~$*$-semigroup~$S$ along with a~partial order $\preceq$ on the
subset of selfadjoint elements $S^{\rm sa}$ of~$S$ satisfying
\begin{gather*}
x^{\ast}ax\preceq x^{\ast}bx\qquad\text{for}\quad a\preceq b\quad\text{in}\quad S^{\rm sa},
\qquad
x\in S.
\end{gather*}
Write $a\prec b$ if $a\preceq b$ and $a\neq b$.

An order map $\beta:S\rightarrow T$ of ordered
$*$-semigroups is a~$*$-map, so mapping $S^{\rm sa}$ to $T^{\rm sa}$, with $\beta(r)\preceq\beta(s)$ if $r\preceq s$ in $S^{\rm sa}$.

An order homomorphism $\beta:S\rightarrow T$ of ordered $*$-semigroups is an order map which is
a~$*$-semigroup homomorphism.
\end{Definition}

If $(T,\preceq)$ is an ordered $*$-semigroup and $\beta:S\rightarrow T$ is a~$*$-map which is a~monomorphism then one
can pull back the partial order on $T^{\rm sa}$ to a~partial order, denoted $\preceq_{\beta}$, on $S^{\rm sa}$.
Def\/ine $a\preceq_{\beta}b$ in $S^{\rm sa}$ if and only if $\beta(a)\preceq\beta(b)$ in $T^{\rm sa}$.
With this ordering $\beta$ becomes an order map.
One may also def\/ine isomorphic ordered $*$-semigroups via semigroup isomorphisms that are order homomorphisms in both directions.

\begin{Definition}
\label{o on o}
The partial orders on a~$*$-semigroup~$S$ possess a~partial order:
\begin{gather*}
(S,\preceq_{1})\preceq(S,\preceq_{2})
\end{gather*}
if the identity homomorphism $(S,\preceq_{1})\rightarrow(S,\preceq_{2})$ is an order map.
\end{Definition}

The semigroups we consider below do not generally have an absorbing (zero) element, but if the $*$-semigroup~$S$ has
a~zero element $0$ then, since representations in $C^*$-algebras is our focus, we assume that $0\preceq a^{\ast}a$ for all $a\in S$.

\begin{Example}
[$C^*$-algebras]Any $C^*$-algebra~$C$ may be viewed as an ordered $*$-semigroup, with multiplication as the semigroup
operation, where $C^{\rm sa}$ has the usual partial ordering def\/ined via the positive cone $C^{+}$.
An order map from an ordered $*$-semigroup~$S$ to a~$C^*$-algebra~$C$ is an order map of~$S$ to the ordered
$*$-semigroup of~$C$.
A~$*$-homomorphism of $C^*$-algebras maps positive elements to positive elements, so is necessarily an order map.
\end{Example}

\begin{Example}[groups]If~$G$ is a~discrete group viewed as a~unital $*$-semigroup, then the self\-adjoint elements $G^{\rm sa}$ are those
elements~$g$ with $g^{-1}=g$, i.e., $g^{2}$ is the unit of~$G$.
A~$*$-ho\-mo\-morphism of~$G$ to a~$C^*$-algebra therefore maps every element of $G^{\rm sa}$ to a~selfadjoint unitary, and
a~partial order on $G^{\rm sa}$ may be def\/ined by using the standard ordering on the projections which def\/ine these
selfadjoint unitaries.
For example if~$U$ is a~selfadjoint unitary then $(U+1)/2$ is a~projection.
Def\/ine a~partial order on $G^{\rm sa}$ by $a\preceq b$ in $G^{\rm sa}$ if and only if $ab=ba$ and $(b-a)^{2}=2(b-a)$ in the
(commutative) polynomial ring $\mathbb{Z}[a,b]=\{n_{1}+n_{2}a+n_{3} b+n_{4}ab\,|\, n_{i}\in\mathbb{Z}\}$.
This ensures the existence of square roots for positive elements in certain operator spaces.
These conditions are equivalent to $a\preceq b$ in $G^{\rm sa}$ if and only if $b(1+a)=(1+a)b=(1+a)$ in $\mathbb{Z}[a,b]$.
One can check that this is a~partial order on $G^{\rm sa}$ satisfying the conditions for an ordered semigroup $(G,\preceq)$,
and that any $*$-homomorphism of~$G$ to a~$C^*$-algebra must be an order map on $(G,\preceq)$.
\end{Example}

\subsection[Partial orders on $*$-subsemigroups of~$A$]{Partial orders on $\boldsymbol{*}$-subsemigroups of~$\boldsymbol{A}$}

The $*$-semigroup~$A$ is an ordered $*$-semigroup.
To describe a~partial order on $A^{\rm sa}$ note that if $n=(n_{0},n_{1},\dots,n_{k})\in A^{\rm sa}$ then $k+1$ must be even,
and $n=(n_{0},n_{1},\dots,n_{l},-n_{l},-n_{l-1},\dots,-n_{0})$.
Note, for later use, that therefore $n=w^{\ast}w$ for some $w\in A$.
For $n_{l}\geq2$ def\/ine
\begin{gather*}
(n_{0},n_{1},\dots,n_{l},-n_{l},-n_{l-1},\dots,-n_{0})\preceq(n_{0},n_{1},\dots,n_{l}-1,-(n_{l}-1),-n_{l-1},\dots,-n_{0}),
\end{gather*}
while if $n_{l}\leq-2$ def\/ine
\begin{gather*}
(n_{0},n_{1},\dots,n_{l},-n_{l},-n_{l-1},\dots,-n_{0})\preceq(n_{0},n_{1},\dots,n_{l}+1,-(n_{l}+1),-n_{l-1},\dots,-n_{0}).
\end{gather*}
If $\left\vert n_{l}\right\vert =1$ then~$n$ is actually $(n_{0},n_{1},\dots,n_{l-1},-n_{l-1},\dots,-n_{0})$ as
a~reduced word, unless $n=(1,-1)$ or $(-1,1)$.
The partial order on $A^{\rm sa}$ is the partial order generated by these relations.

\begin{Definition}[$A,\preceq$]
For $n,m\in A^{\rm sa}$ def\/ine $n\preceq m$ if $n=m$ or if there is a~f\/inite chain $n=n(0)\preceq n(1)\preceq\dots \preceq
n(k)=m$ of (reduced) words $n(i)\in A^{\rm sa}$ with $n(i)\preceq n(i+1)$ as described above.
\end{Definition}

To check that $n\preceq m$ and $m\preceq n$ implies that $n=m$ note that if $n\preceq m$ then the length $l(m)\leq l(n)$.
The symmetric assumption implies $l(n)=l(m)$, and so the length must be constant along any f\/inite chain of reduced words
leading from~$n$ to~$m$.

These relations restrict to $A^{0}$.
Moreover, if $n$, $m$ are selfadjoint and $n\preceq m$ then an integer~$k$ which is an upper bound for the integers
$\sigma_{r}(n)$ for all $r\geq0$ will remain an upper bound for~$\sigma_{r}(m)$ for all $r\geq0$.
Since $D_{0}=\left\{n\in A^{0}\,|\, \sigma_{r}(n)\leq0~\text{for all}~r\geq0\right\}$ (Proposition~\ref{D 1}), and $D_{1}$
has a~similar description (Remark~\ref{b}) this partial order restricts to a~partial order on each of the $*$-semigroups
$D_{0}$ and $D_{1}$.

\begin{Definition}[$(D_{0},\preceq)$ and $(D_{1},\preceq)$]
The given partial order on~$A$ restricted to $D_{0}$ and $D_{1}$ are the partial
orders on the $*$-semigroups $D_{0}$ and $D_{1}$.
\end{Definition}

An application of Lemma~\ref{fund lemma} shows that if $n\preceq m$ in either $D_{0}$ or $D_{1}$ and
$n\in\operatorname*{Irr}\nolimits_{+}(A^{0})$, so a~generator, then $m\in\operatorname*{Irr}\nolimits_{+}(A^{0})$.

We consider an alternative
description of the partial ordered $*$-semigroup $(D_{1},\preceq)$ that is better suited for
describing the matrix ordering introduced in the next section.
It makes formal use of factorizations occurring in the larger semigroup~$A$.
For $n\in A^{\rm sa}$ we already noted that $n=w^{\ast}w$ for some reduced element $w\in A$.
If $w\in A$ is chosen with minimal length this factorization is unique; however~$n$ can be factored as $n=w^{\ast}w$ in
two dif\/ferent ways using reduced elements~$w$ in~$A$.
To see this recall that for $w\in A$ either $w=(-1,1)w$ or $w=(1,-1)w$.
If $n=w^{\ast}w$ and~$w$ is reduced with $(1,-1)w=w$ then $n=\widetilde{w}^{\ast}\widetilde{w}$, where
$\widetilde{w}=(-1)w$, and (the reduced form of) $\widetilde{w}$ provides an alternative
factorization of~$n$.
Similarly if $(-1,1)w=w$ then $n=\widetilde{w}^{\ast}\widetilde{w}$, where $\widetilde {w}=(1)w$ provides another
factorization of~$n$.

\begin{Definition}[alternative
approach to $(A,\preceq)$]
For $n\in A^{\rm sa}$ write $n=w^{\ast}w$ with~$w$ reduced in~$A$.
If $w=(-1,1)w$ def\/ine $w_{+}$ to be a~(reduced) element in~$A$ with $w=(-1)w_{+}$.
Similarly if $w=(1,-1)w$ def\/ine $w_{-}\in A$ satisfying $w=(1)w_{-}$. Then $n=w_{\pm}^{\ast}(\pm1,\mp1)w_{\pm}$.
For $n=w_{\pm}^{\ast}(\pm1,\mp1)w_{\pm}$ def\/ine $n\preceq w_{\pm}^{\ast}w_{\pm}$.
With the partial order generated by these relations~$A$ becomes a~partially ordered $*$-semigroup.
\end{Definition}

Note if $w=(-1)$ then $w_{+}=(1,-1)$ and if $w=(1)$ then $w_{-}=(-1,1)$.
It is usually possible to chose two candidates for $w_{+}$.
For example if $w=(-1,1)w$ then the element $w_{+}=(1)w$ satisf\/ies $w=(-1)w_{+}$, however $w_{+}^{\ast}w_{+}=n$ in this
case.
Similarly if $(1,-1)w=w$.

The characterizations of the $*$-semigroups $D_{0}$ and $D_{1}$ in Proposition~\ref{D 1} and Remark~\ref{b} imply
$w_{\pm}^{\ast}w_{\pm}$ is again an element of $D_{i}^{\rm sa}$ for $n=w^{\ast}w\in D_{i}^{\rm sa}$ ($i=0$ or~$1$).
The restrictions of the partial order therefore yield ordered $*$-semigroups $(D_{1},\preceq)$ and $(D_{0},\preceq)$.

Since there are two factorizations possible for~$n$ there are potentially four elements larger than or equal to~$n$.
However, since the unique minimal length factor always appears as one of these two factoring possibilities, it is
straightforward to check that at least three of these possibilities will always equal the original element~$n$ of
$D_{1}^{\rm sa}$.
The reason for considering both factorizations and not just the one of minimal length may not be apparent now, but for
partial orderings on matrices with entries in~$A$ this becomes necessary to consider (Remark~\ref{hmm}).

If~$n$ is one of the two idempotents in $D_{1}$, this process is stationary and yields the same idempotent~$n$.
If~$n$ is not one of the two idempotents in $D_{1}^{\rm sa}$ then this process always produces an element~$m$ with $n\prec
m$.
One may think of this process as a~hollowing out of a~selfadjoint element~$n$ in $D_{1}^{\rm sa}$.

The def\/inition of $\preceq$ implies that
\begin{gather*}
(-k,k)\preceq(-k+1,k-1)\preceq\dots \preceq(-1,1)
\end{gather*}
for $k\geq1$.
It follows that the maximal elements $(-1,1)$ and $(1,-1)$ are upper bounds for this partial order; namely for any $n\in
D_{1}^{\rm sa}$ we have either $n\preceq(-1,1)$ or $n\preceq(1,-1)$.
The former occurs exactly when $(-1,1)n=n$ and the later occurs exactly when $(1,-1)n=n$.
In particular these elements are the only two maximal elements.
For the $*$-semigroup $D_{0}$ the unit~$(-1,1)$ of~$D_{0}$ is the unique upper bound for this partial order.
Also any two elements of $D_{0}^{\rm sa}$ must have a~least upper bound, although the two elements~$(-1,1)$ and $(1,-1)$
show this is not the case in~$D_{1}$.

Although not needed at present, the following two results provide more detailed information on the form of selfadjoint elements of~$D_{1}$.
These will prove useful for the matricial order properties developed in the next section.

\begin{Lemma}
\label{sa lemma}
If $n\in D_{1}^{\rm sa}$ then~$n$ has one of the forms $m^{\ast}m$ or $m^{\ast}am$, where $m\in D_{1}$ if it occurs, and
$a\in\operatorname*{Irr} (D_{1})^{\rm sa}$ if it occurs.
If~$n$ is written to have a~minimal number of irreducible factors in $D_{1}$ then this form is unique.
If $a=(1,-1)$ then~$m$, if it occurs and is written with a~minimal number of irreducible factors in $D_{1}$, must
satisfy $(-1,1)m=m$.
If $a=(-1,1)$ then~$m$, if it occurs, must satisfy $(1,-1)m=m$.
\end{Lemma}

\begin{proof}
Factorization in $D_{1}$ into irreducibles implies that $n\in D_{1}^{\rm sa}$ must have the form
\begin{gather*}
(1,-1)a_{1}^{\ast}(1,-1)\cdots (1,-1)a_{l}^{\ast}(1,-1)a(1,-1)a_{l} (1,-1)\cdots (1,-1)a_{1}(1,-1)
\end{gather*}
uniquely as a~product, where $a_{i}\in D_{0}$, $a\in\operatorname*{Irr} (D_{1})^{\rm sa}$.
The initial and f\/inal $(1,-1)$ terms, the $(1,-1)$ pair bracketing~$a$, as well as~$a$ itself, may or may not occur.

If~$a$ occurs and is equal to $(1,-1)$ then the unique product must be
\begin{gather*}
(1,-1)a_{1}^{\ast}(1,-1)\cdots (1,-1)a_{l}^{\ast}(1,-1)a_{l}(1,-1)\cdots (1,-1)a_{1} (1,-1),
\end{gather*}
where the initial and f\/inal $(1,-1)$ pair either occur or not, but the $(1,-1)$ pair around~$a$ does not occur (by the
minimality condition on~$m$).
In this case $m=a_{l}(1,-1)\cdots (1,-1)a_{1}(1,-1)$, if~$m$ occurs, (perhaps without the f\/inal $(1,-1)$ term) and
$(-1,1)m=m$.

If~$a$ occurs and is equal to $(-1,1)$ then, if~$m$ occurs, the unique product must be
\begin{gather*}
(1,-1)a_{1}^{\ast}(1,-1)\cdots (1,-1)a_{l}^{\ast}(1,-1)(-1,1)(1,-1)a_{l} (1,-1)\cdots (1,-1)a_{1}(1,-1).
\end{gather*}
The pair of interior terms $(1,-1)$ must occur, since otherwise
\begin{gather*}
n=(1,-1)a_{1}^{\ast}(1,-1)\cdots (1,-1)a_{l}^{\ast}a_{l}(1,-1)\cdots (1,-1)a_{1} (1,-1)
\end{gather*}
contradicting the assumption that~$a$ occurs.
Thus~$m$, if it occurs, must be
\begin{gather*}
(1,-1)a_{l}(1,-1)\cdots (1,-1)a_{1}(1,-1)\qquad\text{and}\qquad (1,-1)m=m.\tag*{\qed}
\end{gather*}
  \renewcommand{\qed}{}
\end{proof}

It is now possible to compare the two dif\/fering approaches, each unique, to describe an element $n\in D_{1}^{\rm sa}:$ one
as $w^{\ast}w$ with $w\in A$ of minimal length, the other as in Lemma~\ref{sa lemma}, using products of elements in
$D_{1}$.
Consider two ordered sequences, indexed by increasing $r\in\mathbb{N}$, namely the sequence $\sigma_{r}(w^{\ast}w)$, and
the sequence $\sigma_{r}(m^{\ast}am)$ (or $\sigma_{r}(m^{\ast}m)$, depending on which case of Lemma~\ref{sa lemma}
holds).
These two sequences are generally not the same.
The sequence $\sigma_{r}(w^{\ast}w)$ has an initial odd number of terms that occur before the (inf\/inite) string of zero
terms, so has a~well def\/ined middle term about which the sequence is symmetric.
In particular the sequence $\sigma_{r}(w^{\ast})$ occurs as the initial half of this sequence, ending with the term
$\tau(w^{\ast})$, which is this middle term.
The sequence $\sigma_{r}(w^{\ast}w)$ is an ordered subsequence of the sequence $\sigma_{r}(m^{\ast}am)$, indeed the
sequence $\sigma_{r}(m^{\ast}am)$ is also symmetric about a~middle term, and only has some possible extra~$0$ (partial
sum) terms occurring.
In particular the sequence $\sigma_{r}(w^{\ast})$ occurs as a~subsequence of the initial half of the sequence
$\sigma_{r}(m^{\ast}am)$.
Thus $\tau(w^{\ast})$ also occurs as the middle term in the sequence $\sigma_{r}(m^{\ast}am)$.
By Remark~\ref{b}, $D_{1}=\left\{n\in A^{0} \,|\,\sigma_{r}(n)\leq1~\text{for all}~r\geq0\right\}$, so $\sigma
_{r}(w^{\ast}w)\leq1$ for all $r\geq0$, and therefore $\tau(w^{\ast})\leq1$.
Note that the f\/irst half of the sequence $\sigma_{r}(m^{\ast}am)$ begins with the sequence $\sigma_{r}(m^{\ast})$, which
ends in~$0$ if~$m$ occurs, followed by the f\/irst half of the sequence~$\sigma_{r}(a)$ if~$a$ occurs.
These types of observations yield the following.

\begin{Proposition}
\label{details}
Write $n\in D_{1}^{\rm sa}$ as $w^{\ast}w$ with $w\in A$ reduced of minimal length, and also as $m^{\ast}m$ or $m^{\ast}am$,
where $m\in D_{1}$ if it occurs, has a~minimal number of irreducible factors in $D_{1}$, and
$a\in\operatorname*{Irr}(D_{1})^{\rm sa}$ if it occurs.
Then $\tau(w^{\ast})\leq1$.
We have the following equivalences:
\begin{itemize}\itemsep=0pt
\item $\tau(w^{\ast})=1$ if and only if $n=m^{\ast}(1,-1)m$, where $m=(-1,1)m$ may not occur,

\item $\tau(w^{\ast})=0$ if and only if $n=m^{\ast}m$ and $m=(-1,1)m$, and

\item $\tau(w^{\ast})\leq-1$ if and only if $n=m^{\ast}am$ with $a\in \operatorname*{Irr}(D_{0})^{\rm sa}$, where~$m$ may not occur.
\end{itemize}
\end{Proposition}

\subsection{Schwarz property}

\begin{Definition}
A~$*$-map $\beta:S\rightarrow T$ of a~$*$-semigroup~$S$ to $(T,\preceq)$ an ordered $*$-semigroup~$T$ is said to have
the Schwarz property if it satisf\/ies the Schwarz inequality
$\beta(a)^{\ast}\beta(a)\preceq\beta(a^{\ast}a)
$
for all $a\in S$.
\end{Definition}

The composition of two order maps with the Schwarz property also has the Schwarz pro\-per\-ty.
If $\gamma:S\rightarrow T$ is a~$*$-semigroup homomorphism then $\gamma$ clearly has the Schwarz property.
Therefore if $\beta:R\rightarrow S$ an order map with the Schwarz property and $\gamma:S\rightarrow T$ a~$*$-semigroup
homomorphism as well as an order map, then $\gamma\beta$ is an order map with the Schwarz property.
It seems possible to consider maps of $*$-semigroups that are not $*$-maps, but that still map the selfadjoint elements
to selfadjoint elements, since there are straightforward examples, satisfying Schwarz or not, of such maps.
However, for our present purposes this is not needed.
The next proposition also holds for the $*$-semigroup~$A$, but we only make use of it for~$D_{1}$.

\begin{Proposition}\label{square}
The ordered $*$-semigroup $(D_{1},\preceq)$ satisfies
$a^{2}\preceq a$
for $a\in D_{1}^{\rm sa}$.
\end{Proposition}

\begin{Proposition}
The $*$-map $\alpha:D_{1}\rightarrow D_{1}$ is an order map satisfying the Schwarz property
$\alpha(a)^{\ast}\alpha(a)\preceq\alpha(a^{\ast}a)
$
for $a\in D_{1}$.
\end{Proposition}

\begin{proof}
For the f\/inal statement $a^{\ast}(1,-1)a\preceq a^{\ast}a$ for $a\in D_{1}$, so applying $\alpha$ yields
$\alpha(a)^{\ast}\alpha(a)=\alpha(a^{\ast}(1,-1)a)\preceq\alpha(a^{\ast}a)$.
\end{proof}

\begin{Proposition}
\label{before}
The $*$-homomorphism $\omega:D_{0}\rightarrow D_{1}$ is an order map of ordered $*$-semi\-groups, so a~Schwarz map.
\end{Proposition}

\begin{proof}
By def\/inition $\omega(a)=(1)a(-1)$ for $a\in D_{0}$.
It is clear that the basic one step relation generating the partial order on $D_{0}^{\rm sa}$ is preserved by this map.
\end{proof}

\begin{Definition}[order representation]An order representation of an ordered $*$-semigroup~$S$ is a~$*$-homomorphism to a~$C^*$-algebra
which is, in addition, an order map to this $C^*$-algebra.
\end{Definition}

Any order representation of an ordered $*$-semigroup in a~$C^*$-algebra clearly satisf\/ies the Schwarz property.

Recall that if~$A$ and~$B$ are selfadjoint elements in a~$C^*$-algebra~$C$ with $-B\leq A\leq B$ then $\Vert
A\Vert \leq\Vert B\Vert$; in particular $0\leq A\leq B$ implies $\Vert A\Vert \leq\Vert
B\Vert$.
Note also if $\beta:S\rightarrow C$ is a~$*$-homomorphism
to a~$C^*$-algebra and $a\in S^{\rm sa}$ then the element
$\beta(a)$ has norm bounded by 1 if and only if $\Vert \beta(a^{2})\Vert \leq\Vert \beta(a)\Vert$.

\begin{Proposition}
\label{gen square}
Let $(S,\preceq)$ be an ordered $*$-semigroup with $a^{2}\preceq a$ for all $a\in S^{\rm sa}$,~$C$ a~$C^*$-algebra and
$\beta:S\rightarrow C$ an order map.
If $\beta$ satisfies the Schwarz inequality $\beta(a)^{\ast}\beta(a)\leq\beta(a^{\ast}a)$ for $a\in S^{\rm sa}$ then $\beta$
maps the selfadjoint elements $S^{\rm sa}$ to the unit ball of the positive cone of~$C$.

If in addition $\beta$
satisfies the Schwarz inequality for all $a\in S$ then $\beta(S)$ is contained in the unit ball of~$C$.
\end{Proposition}

\begin{proof}
For $a\in S^{\rm sa}$ we have $\beta(a)\in C^{\rm sa}$, so
\begin{gather*}
0\leq\beta(a)^{2}\leq\beta\big(a^{2}\big)\leq\beta(a),
\end{gather*}
where the second inequality uses the Schwarz inequality for $a\in S^{\rm sa}$ and the last inequality follows from the
hypothesis.
For $a\in S^{\rm sa}$ the $C^*$-norm condition implies
\begin{gather*}
\Vert \beta(a)\Vert^{2}=\Vert \beta(a)^{2}\Vert \leq\Vert \beta(a)\Vert,
\end{gather*}
and so $\Vert \beta(a)\Vert \leq1$.
If $a\in S$ then using the Schwarz inequality for general $a\in S$ and the $C^*$-norm condition yields
\begin{gather*}
\Vert \beta(a)\Vert^{2}=\Vert \beta(a)^{\ast}\beta (a)\Vert \leq\Vert
\beta(a^{\ast}a)\Vert \leq1,
\end{gather*}
the last inequality following from $a^{\ast}a\in S^{\rm sa}$.
\end{proof}

Since an order representation of an ordered $*$-semigroup has the Schwarz property, the conclusions of
Proposition~\ref{gen square} hold for all order representations of an ordered $*$-semigroup~$S$ satisfying $a^{2}\preceq
a$ for $a\in S^{\rm sa}$.
Thus an order representation of such an ordered $*$-semigroup~$S$ is automatically a~contractive $*$-representation.

By Proposition~\ref{square} the preceding comments apply to the ordered $*$-semigroup $D_{1}:$ an order representation
$\sigma:D_{1}\rightarrow C$ must be contractive, and maps $D_{1}^{\rm sa}$ to the unit ball of the positive cone of the
$C^*$-algebra~$C$.
Similarly for the ordered unital $*$-semigroup $D_{0}$.
Furthermore, since the unit $(-1,1)$ of $D_{0}$ is the unique upper bound for the partial order,
$0\leq\sigma(a)\leq\sigma((-1,1))$ for $a\in D_{0}^{\rm sa}$ and $\sigma:D_{0}\rightarrow C$ an order map.
If $\sigma$ also has the Schwarz property then
\begin{gather*}
0\leq\sigma(a)^{\ast}\sigma(a)\leq\sigma(a^{\ast}a)\leq\sigma((-1,1))
\end{gather*}
for all $a\in D_{0}$.
Thus $\sigma((-1,1))=0$ if and only if $\sigma$ is $0$.

\subsection[$C^*$-algebras of ordered $*$-semigroups]{$\boldsymbol{C^*}$-algebras of ordered $\boldsymbol{*}$-semigroups}

\begin{Definition}
Given $\left(S,\preceq\right) $ an ordered $*$-semigroup, $C^{\ast} (S,\preceq)$ is a~$C^*$-algebra along with an order
representation $\iota:S\rightarrow C^{\ast}(S,\preceq)$ satisfying the following universal property:

for $\gamma:S\rightarrow B$ an order representation to a~$C^*$-algebra~$B$ there is a~unique $*$-homomorphism
$\pi_{\gamma}=\pi:C^{\ast}(S,\preceq)\rightarrow B$ such that $\pi_{\gamma}\circ\iota=\gamma$.
\end{Definition}

That such a~universal $C^*$-algebra of a~given ordered $*$-semigroup exists would follow if there are order
representations of~$S$, and if for each $s\in S$
\begin{gather*}
\sup\{\Vert \beta(s)\Vert \,|\,\beta~\text{is an order representation of}~S\}<\infty.
\end{gather*}
If this supremum is bounded by~$1$ on~$S$ for an ordered $*$-semigroup~$S$ it follows that $C^{\ast}(S,\preceq)$ must be
a~quotient of the universal $C^*$-algebra $C^{\ast}(S)$.
Dif\/fering conditions for various ordered $*$-semigroups~$S$ would guarantee this condition.
For example, Proposition~\ref{gen square} shows if $a^{2}\preceq a$ for $a\in S^{\rm sa}$ then this $\sup$ is bounded by
$1$, so the universal $C^*$-algebra $C^{\ast}(S,\preceq)$ exists and furthermore is a~quotient of $C^{\ast}(S)$.

If $\beta:S\rightarrow T$ is an order homomorphism of ordered $*$-semigroups, where the universal $C^*$-algebra of the
ordered $*$-semigroups exist then there is a~$*$-homomorphism $\pi_{\beta}:C^{\ast}\left(S,\preceq\right)\rightarrow C^{\ast}\left(T,\preceq\right)$.
In particular if $(S,\preceq_{1})\preceq(S,\preceq_{2})$ (Def\/inition~\ref{o on o}) then
$C^{\ast}\left(S,\preceq_{2}\right) $ is a~quotient of $C^{\ast}\left(S,\preceq_{1}\right)$.

\begin{Remark}
\label{orders}
Assume $(S,\preceq)$ is an ordered $*$-semigroup and that $C^{\ast}(S,\preceq)$ exits.
If the cano\-ni\-cal order representation $\iota:S\rightarrow C^{\ast}(S,\preceq)$ is a~monomorphism then one can form the
pull back partial order $\preceq_{\iota}$ on~$S$.
Since $\iota$ is an order representation it is clear that $(S,\preceq)\preceq(S,\preceq_{\iota})$, so
$C^{\ast}\left(S,\preceq_{\iota}\right)$ (with order representation $\iota_{1}:(S,\preceq_{\iota})\rightarrow
C^{\ast}(S,\preceq_{\iota})$) is a~quotient of $C^{\ast}\left(S,\preceq\right)$.
However, $a\preceq_{\iota}b$ in $S^{\rm sa}$ is equivalent to $\iota(a)\leq\iota(b)$, which implies
\begin{gather*}
\beta(a)=\pi_{\beta}\iota(a)\leq\pi_{\beta}\iota(b)=\beta(b)
\end{gather*}
for all order representations $\beta$ of $\left(S,\preceq\right)$.
Therefore every order representation $\beta$ of $(S,\preceq)$ is also an order representation of $(S,\preceq_{\iota})$,
in particular when $\beta$ is the order representation $\iota$.
This yields a~$*$-homomorphism
\begin{gather*}
\pi_{\iota}: \ C^{\ast}(S,\preceq_{\iota})\rightarrow C^{\ast}(S,\preceq)
\end{gather*}
with $\iota_{1}\pi_{\iota}=\iota$, so $C^{\ast}(S,\preceq_{\iota})\cong C^{\ast}(S,\preceq)$.
\end{Remark}

\begin{Example}[$C^{\ast}(G,\preceq)$]
We noted earlier, using our previous description of an ordering on a~discrete group~$G$, that any
$*$-homomorphism of~$G$ to a~$C^*$-algebra must also be an order representation.
Therefore the universal $C^*$-algebra $C^{\ast}(G,\preceq)$ is (isomorphic to) the universal $C^*$-algebra
$C^{\ast}(G)$.
\end{Example}

\begin{Example}
[$C^{\ast}(\mathbb{N},\preceq)$]Consider the $*$-semigroup $\mathbb{N}$ with the partial order $n\preceq m$ if and only
if $n=m$ or there is an $r\in \mathbb{N}$ with $n=m+r$.
Then (using additive notation) $a+a\preceq a$ for $a\in\mathbb{N}$ and it follows that $C^{\ast}(\mathbb{N},\preceq)$ is
the universal $C^*$-algebra generated by a~positive contraction, so is isomorphic to $C_{0}((0,1])$.
With a~similar partial order on $\mathbb{N}_{0}$ it follows that
\begin{gather*}
C^{\ast}(\mathbb{N}_{0},\preceq)\cong C([0,1]),
\end{gather*}
the universal unital $C^*$-algebra generated by a~positive contraction.
\end{Example}

\begin{Example}
[$C^{\ast}(A)$]Consider the ordered $*$-semigroup~$A$.
Noting the form of any element of $A^{\rm sa}$ it is clear that any $*$-representation of~$A$ in a~$C^*$-algebra is also an
order representation of $(A,\preceq)$.
Thus
\begin{gather*}
C^{\ast}(A,\preceq)\cong C^{\ast}(A).
\end{gather*}
\end{Example}

Since $D_{1}$, the ordered $*$-semigroup we are principally interested in, satisf\/ies the order condition $a^{2}\preceq
a$ for $a\in D_{1}^{\rm sa}$ (Proposition~\ref{square}) it follows that $C^{\ast}(D_{1},\preceq)$ exists and is a~quotient
of the universal $C^*$-algebra $C^{\ast}(D_{1})$.
Similarly the universal $C^*$-algebra $C^{\ast}(D_{0},\preceq)$ exists.
Let
\begin{gather*}
i_{0}: \ D_{0}\rightarrow C^{\ast}(D_{0},\preceq)\qquad\text{and}\qquad i_{1}:\ D_{1} \rightarrow C^{\ast}(D_{1},\preceq)
\end{gather*}
denote the canonical order representations.
Proposition~\ref{basic} implies that the $*$-algebra $\mathbb{C}[i_{1}(D_{1})]$ $=\widetilde{i_{1}}\mathbb{C} [D_{1}]$,
where $\widetilde{i_{1}}$ is the linear extension of $i_{1}$ to the $*$-algebra $\mathbb{C}[D_{1}]$, is a~dense
$*$-subalgebra of $C^{\ast} (D_{1},\preceq)$.

\begin{Remark}
\label{dense}
If $a\in\mathbb{C}[(D_{1})]$ then $\widetilde{i_{1}}(a)$ has norm in $C^{\ast}(D_{1},\preceq)$ equal to
\begin{gather*}
\sup\big\{\big\Vert \widetilde{\beta}(a)\big\Vert \;\big|\; \beta~\text{is an order representation of}~D_{1}\big\}.
\end{gather*}
A~similar statement holds for general ordered $*$-semigroups~$S$ whenever this sup is bounded on~$S$.
\end{Remark}

\begin{Remark}
\label{order}
Not every contractive $*$-homomorphism of an ordered $*$-semigroup need be an order representation, so the universal
$C^*$-algebra $C^{\ast}(S,\preceq)$ is, in general, dif\/ferent from $C^{\ast}(S)$.
Consider, for example, the unital ordered \allowbreak$*$-semigroup $D_{0}$, which may be viewed as the free unital
$*$-semigroup on $F=\operatorname*{Irr}\nolimits_{+} (A^{0})\setminus(-1,1)$ (Remark~\ref{free}).
Therefore an arbitrary $*$-map of~$F$ to a~set of elements in the unit ball of a~unital $C^*$-algebra~$C$ def\/ines
a~unital contractive $*$-homomorphism of $D_{0}$ to~$C$ that generally is not an order representation.
\end{Remark}

\begin{Proposition}
\label{extension 0}
Given an order representation $\sigma:D_{0}\rightarrow C$ to a~$C^*$-algebra~$C$, there is an order representation
$\rho:D_{1}\rightarrow C$ that extends $\sigma$.
\end{Proposition}

\begin{proof}
We may assume that~$C$ contains a~nonzero projection~$q$, since if the projection $\sigma(-1,1)$ of~$C$ is $0$ then the
comments following Proposition~\ref{gen square} imply $\sigma=0$, and $\rho=0$ extends $\sigma$.
For $a\in D_{1}$ write
\begin{gather*}
a=(1,-1)a_{1}(1,-1)a_{2}(1,-1)\cdots (1,-1)a_{k}(1,-1)
\end{gather*}
uniquely as a~product, where $a_{i}\in D_{0}$ (where the initial and f\/inal $(1,-1)$ may or may not occur).
Def\/ine
\begin{gather*}
\rho(a)=q\sigma(a_{1})q\sigma(a_{2})q\cdots q\sigma(a_{k})q
\end{gather*}
(again,
where the initial and f\/inal~$q$ may or may not occur).
It follows from uniqueness of the product that $\rho$ is well def\/ined and def\/ines a~$*$-homomorphism extending $\sigma$.
Finally, a~basic step in the order relation def\/ined on $D_{1}$ is preserved under $\rho$ since, for example,
$\sigma(a_{k}^{\ast})q\sigma(a_{k})\leq\sigma(a_{k}^{\ast})\sigma(a_{k})$ in~$C$.
Thus $\rho$ is an order map.
\end{proof}

For $\eta:(D_{0},\preceq)\rightarrow(D_{1},\preceq)$ the natural inclusion, universality yields a~$*$-homomorphism
\begin{gather*}
\pi_{\eta}:\ C^{\ast}\left(D_{0},\preceq\right) \rightarrow C^{\ast} (D_{1},\preceq)
\end{gather*}
with $\pi_{\eta}\circ i_{0}=i_{1}\circ\eta$ (with $i_{k}:D_{k}\rightarrow C^{\ast}(D_{k},\preceq)$ the canonical order
representations).

\begin{Corollary}
\label{injection 0}
The $*$-homomorphism $\pi_{\eta}:C^{\ast}(D_{0},\preceq)\rightarrow C^{\ast}(D_{1},\preceq)$ is an injection of $C^*$-al\-gebras.
\end{Corollary}

\begin{proof}
If $a\in\mathbb{C}[D_{0}]^{\rm sa}$ then (Remark~\ref{dense}) the norm
\begin{gather*}
\Vert \widetilde{i_{0}}(a)\Vert
=\sup\{\Vert \widetilde {\sigma}(a)\Vert \,|\,\sigma~\text{is an order representation of}~D_{0}\},
\end{gather*}
where $\widetilde{\;}$ denotes the linear extensions to $\mathbb{C}[D_{0}]$.
For such a~$\sigma$ there is, by Proposition~\ref{extension 0}, an order representation $\rho_{\sigma}$ of $D_{1}$
with $\rho_{\sigma}\circ \eta=\sigma$.
Therefore
\begin{gather*}
\Vert \widetilde{i_{0}}(a)\Vert =\sup\{\Vert \widetilde {\rho_{\sigma}}
\widetilde{\eta}(a)\Vert \,|\,\sigma=\rho_{\sigma} \circ\eta~\text{is an order representation of}~D_{0}\}
\\
\phantom{\Vert \widetilde{i_{0}}(a)\Vert}
\leq\sup\{\Vert \widetilde{\rho}
\widetilde{\eta}(a)\Vert \,|\, \rho~\text{is an order representation of}~D_{1}\}.
\end{gather*}
The latter is equal to $\Vert \widetilde{\iota_{1}}(\widetilde{\eta}(a))\Vert =\Vert
\pi_{\eta}(\widetilde{i_{0}}(a))\Vert$.
Since $\pi_{\eta}$ is also norm reducing, it must be an isometry on the dense $*$-subalgebra $\mathbb{C}[i_{0}(D_{0})]$
of $C^{\ast}(D_{0},\preceq)$.
\end{proof}

\section[A~$*$-semigroup matricial order]{A~$\boldsymbol{*}$-semigroup matricial order}
\label{section4}

Our main interest is to extend the order on the $*$-semigroups considered above to a~partial order on certain
selfadjoint elements in matrices over the semigroup.
The present approach suf\/f\/ices for our context and provides indications of what matricial partial orders for general
$*$-semigroups should entail; a~more general approach could see future ref\/inements as further examples are investigated.

For a~$*$-semigroup~$S$ we introduce matricial orders as well as their associated maps on these semigroups, namely
complete order maps and order representations,~$k$-amplif\/ications that have the Schwarz property, and completely
positive maps to a~$C^*$-algebra.
These lead to a~universal $C^*$-algebra $C^{\ast}((S,\preceq,\mathcal{M}))$ for a~matricially ordered $*$-semigroup
$(S,\preceq,\mathcal{M})$.
With these def\/initions a~$*$-homomorphism of a~$*$-semigroup to a~$C^*$-algebra will always satisfy the Schwarz property
for all~$k$, and a~$*$-map to a~$C^*$-algebra with the Schwarz property for all~$k$ must be completely positive.
A~completely positive map on a~$*$-semigroup~$S$ to a~$C^*$-algebra~$C$ determines a~right Hilbert module over~$C$.

We def\/ine a~matricial order on the $*$-semigroup $D_{1}$, and after exploring this partial order, we show that there is
a~$C^*$-correspondence $\mathcal{E}$ (Theorem~\ref{left action}) over the universal $C^*$-algebra $C^{\ast}
((D_{1},\preceq,\mathcal{M}))$ for complete order maps on $D_{1}$.
The f\/inal part of this section involves some detailed matricial order properties used to obtain
Proposition~\ref{extension}, a~complete order analogue of the previous extension result Proposition~\ref{extension 0}.
These properties will prove essential in identifying an isomorphic representation of the universal $C^*$-algebra
$C^{\ast}((D_{1},\preceq,\mathcal{M}))$ (Corollaries~\ref{need} and~\ref{iso image}).

\subsection{Matricial order and completely positive maps}

Let $(S,\preceq)$ be an ordered $*$-semigroup.
To extend partial orders on~$S$ to matrices over~$S$ initially appears problematic since the collection of matrices
over~$S$ do not inherit suf\/f\/icient algebraic structure from~$S$.
However, for specif\/ic types of matrices over~$S$ this is possible.
Additionally this requires considering how these specif\/ic types of matrices in $\mathrm{M}_{d}(S)$ may be viewed as
matrices in $\mathrm{M}_{k}(S)$.

For $k\in\mathbb{N}$, $[n_{i}]$ denotes an element $[n_{1},\dots,n_{k}]\in\mathrm{M}_{1,k}(S)$.
Then $[n_{i}]^{\ast}\in\mathrm{M}_{k,1}(S)$.
One of the basic types of matrices over~$S$ considered has the form
$[n_{i}]^{\ast}[n_{j}]=[n_{i}^{\ast}n_{j}]\in\mathrm{M}_{k}(S)$.

\begin{Notation}
For $d,k\in\mathbb{N}$ and $d\leq k$, set
\begin{gather*}
\mathcal{P}(d,k)=\left\{(t_{1},\dots,t_{d})\in(\mathbb{N}_{0})^{d}\;\Bigg|\; \sum\limits_{r=1}^{d}t_{r}=k\right\}.
\end{gather*}
Set
$\delta=(1,\dots,1)$
the unique element of $\mathcal{P}(k,k)$ with nonzero entries.
\end{Notation}

One can view the~$d$-tuples in $\mathcal{P}(d,k)$ as ordered partitions of~$k$, where zero summands are allowed.
The set $\mathcal{P}(1,k)$ consists of the single element $(k)$.

Elements of $\mathcal{P}(d,k)$ are applied in two ways, one to view elements $[n_{i}]\in\mathrm{M}_{1,k}(S)$ as certain
formal $d\times k$ matrices, and secondly as maps $\mathrm{M}_{d}(S)\rightarrow\mathrm{M}_{k}(S)$.
The context should make these dif\/ferent uses clear.

For $[n_{i}]\in\mathrm{M}_{1,k}(S)$ and $\tau=(t_{1},\dots,t_{d})\in \mathcal{P}(d,k)$, $[n_{i}]_{\tau}$ denotes the
formal $d\times k$ matrix whose entries from~$S$ in the~$i$-th row is the length $t_{i}$ string $(n_{t_{1}+\dots
+t_{i-1}+1},\dots,n_{t_{1}+\dots +t_{i}})$, starting with $n_{t_{1}+\dots +t_{i-1}+1}$ at the $t_{1}+\dots +t_{i-1}+1$ spot.
All other elements are denoted $0$.
If necessary set $t_{0}=0$.
The word formal is used since these matrices are not strictly over~$S$ unless~$S$ has an absorbing element $0$.
Basically the~$k$ elements of $[n_{i}]$ are displayed in~$d$ rows using the partition~$\tau$.
The matrix
\begin{gather*}
[n_{i}]_{\tau}=\left[
\begin{array}[c]{@{}lllll@{}}%
n_{1},\dots ,n_{t_{1}}, & 0,\dots \text{ \ \ \ \ \ \ \ \ \ \ \ \ \ \ } & \dots & \dots &
\dots \text{ \ \ \ \ \ \ \ \ \ \ \ \ \ \ \ \ \ },0\\
0,\dots \text{ \ },0 & n_{t_{1}+1},\dots ,n_{t_{1}+t_{2}} & 0,\dots & \dots & \dots \text{
\ \ \ \ \ \ \ \ \ \ \ \ \ \ \ \ \ },0\\
0,\dots  &  &  &  & \dots \text{ \ \ \ \ \ \ \ \ \ \ \ \ \ \ \ \ \ },0\\
\dots  & \dots  &  &  & \\
0,\dots \text{ \ } & \dots \text{ \ \ \ \ \ \ \ \ \ \ \ } & ,0, & .,n_{t_{1}%
+\dots +t_{d-1}}, & 0,\dots\text{ \ \ \ \ \ \ \ \ \ \ \ \ \ \ },0\\
0,\dots \text{ \ } & \dots \text{ \ \ \ \ \ \ \ \ \ \ \ } &  & \text{
\ \ \ \ \ \ \ }\dots ,0, & n_{t_{1}+\dots +t_{d-1}+1},\dots ,n_{k}%
\end{array}
\right].
\end{gather*}
For $[n_{i}]\in\mathrm{M}_{1,k}(S)$, and $(k)$ the element of $P(1,k)$, $[n_{i}]_{(k)}=[n_{i}]$.
For the element $\delta\in\mathcal{P} (k,k)$, $[n_{i}]_{\delta}$ is the $k\times k$ diagonal matrix with $n_{i}$ in the
$(i,i)$ entry.

Each $\tau=(t_{1},\dots,t_{d})\in\mathcal{P}(d,k)$ also determines a~$*$-map
$\iota_{\tau}:\mathrm{M}_{d}(S)\rightarrow\mathrm{M}_{k}(S)$ of square matrices with entries in~$S$.
An element $\left[a_{i,j}\right] \in\mathrm{M}_{d}(S)$ is mapped to
$\iota_{\tau}(\left[a_{i,j}\right]):=\left[a_{i,j}\right]_{\tau}\in\mathrm{M}_{k}(S)$, where $\left[a_{i,j}\right]
_{\tau}$ is a~matrix described via matrix blocks.
The $i$, $j$ block of $\left[a_{i,j}\right]_{\tau}$ is the $t_{i}\times t_{j}$ matrix with the constant entry $a_{i,j}:$
\begin{gather*}
\left[
\begin{array}[c]{@{}rrrrr@{}}
\left[
\begin{array}[c]{@{}l@{}}
a_{11},\dots ,a_{11}\\
\dots \\
a_{11},\dots ,a_{11}
\end{array}
\right]  & \dots  &
\begin{array}
[c]{r}%
j\text{-th column}
\\
\multicolumn{1}{c}{\quad\text{width}}
\\
\multicolumn{1}{c}{\longleftarrow t_{j}\longrightarrow}%
\end{array}
& \dots  & \left[
\begin{array}[c]{@{}l@{}}
a_{1d},\dots ,a_{1d}\\
\dots \\
a_{1d},\dots ,a_{1d}
\end{array}
\right] \\
\dots  &  &  & \dots  & \dots \\
i\text{-th row},%
\begin{array}[c]{@{}c@{}}
\uparrow\\
t_{i}\\
\downarrow
\end{array}
&  & \left[
\begin{array}[c]{@{}l@{}}
a_{ij},\dots ,a_{ij}\\
\dots \\
a_{ij},\dots ,a_{ij}%
\end{array}
\right]  &  & \\
& \dots  & \dots  &  &
\end{array}
\right].
\end{gather*}
The map $\iota_{\delta}$ is the identity map on $\mathrm{M}_{k}(S)$.
If $\tau=(t_{1},\dots,t_{d})\in\mathcal{P}(d,k)$ has zero entries then the map $\iota_{\tau}$ basically restricts to an
element $\widetilde{\tau} \in\mathcal{P}(\widetilde{d},k)$ for an appropriate $\widetilde{d}<d$, so $\iota_{\tau}$ may
actually be viewed as a~map from $\mathrm{M}_{\widetilde {d}}(S)$ to $\mathrm{M}_{k}(S)$.

There is a~(partial) composition that can be def\/ined for elements $\tau =(t_{1},\dots,t_{d})\in\mathcal{P}(d,h)$ and
$\sigma=(s_{1},\dots,s_{h})\in\mathcal{P}(h,k)$ to obtain an element $\gamma=\sigma\circ\tau$, where
\begin{gather*}
\gamma=(g_{1},\dots,g_{d})\in\mathcal{P}(d,k)
\qquad
\text{with}
\quad
g_{j}=\sum_{l=1}^{t_{j}}s_{t_{j-1}+l}.
\end{gather*}
We have $\iota_{\sigma}\circ\iota_{\tau}=\iota_{\sigma\circ\tau}$.

For $\left[a_{i,j}\right] \in\mathrm{M}_{d}(S)$, $[m_{i}]$, $[n_{i}]\in\mathrm{M}_{1,k}(S)$ and $\tau\in\mathcal{P}(d,k)$
the $k\times k$ matrix
\begin{gather*}
[m_{i}]_{\tau}^{\ast}\left[a_{i,j}\right] [n_{i}]_{\tau} =[m_{i}]_{\delta}^{\ast}\left[a_{i,j}\right]_{\tau}[n_{i}]_{\delta}.
\end{gather*}
In particular
\begin{gather*}
[n_{i}]^{\ast}a[n_{j}]=[n_{i}]_{\delta}^{\ast}\left[a\right]_{(k)}[n_{j}]_{\delta}=[n_{i}^{\ast}an_{j}]
\end{gather*}
for $a\in S=\mathrm{M}_{1}(S)$, and
\begin{gather*}
[n_{i}]_{\delta}^{\ast}\left[a_{i,j}\right] [n_{j}]_{\delta} =[n_{i}^{\ast}a_{i,j}n_{j}].
\end{gather*}
Although $[n_{i}]_{\delta}$ is generally not in $\mathrm{M}_{k}(S)$ (unless~$S$ has an absorbing zero element), the
matrices $\left[a_{i,j}\right]_{\tau} [n_{j}]_{\delta}$ and $[m_{i}]_{\delta}^{\ast}\left[a\right]_{(k)}
[n_{j}]_{\delta}$ are always in $\mathrm{M}_{k}(S)$.

\begin{Definition}
A~$*$-semigroup~$S$ is matricially ordered, write $(S,\preceq,\mathcal{M})$, if there is a~sequence of partially ordered
sets $(\mathcal{M}_{k}(S),\preceq)$, $\mathcal{M}_{k}(S)\subseteq\mathrm{M}_{k}(S)^{\rm sa}$ ($k\in\mathbb{N}$), and
$\mathcal{M}_{1}(S)=S^{\rm sa}$, satisfying (for $[n_{i}]\in\mathrm{M}_{1,k}(S)$)
\begin{enumerate}\itemsep=0pt
\item[a)] $[n_{i}]^{\ast}[n_{j}]=[n_{i}^{\ast}n_{j}]\in\mathcal{M}_{k}(S)$,

\item[b)] $\mathcal{M}_{k}(S)$ is closed under conjugation by the elements $[n_{i}]_{\delta}$,

\item[c)] if $\left[a_{i,j}\right] \preceq\left[b_{i,j}\right] $ in $\mathcal{M}_{k}(S)$ then
$[n_{i}]_{\delta}^{\ast}\left[a_{i,j}\right] [n_{j}]_{\delta}\preceq[n_{i}]_{\delta}^{\ast}\left[b_{i,j}\right][n_{j}]_{\delta}$,

\item[d)]
for $\tau\in\mathcal{P}(d,k)$, $\iota_{\tau}:\mathcal{M}_{d}(S)\rightarrow \mathcal{M}_{k}(S)$ are order maps.
\end{enumerate}
\end{Definition}

\begin{Example}[$C^*$-algebras]
For~$C$ a~$C^*$-algebra set $\mathcal{M}_{k}(C)$ to be the usual partially ordered set~$\mathrm{M}_{k}(C)^{\rm sa}$.
To see that $\mathcal{M}_{k}(C)$ satisf\/ies the requirements f\/irst note that the maps
$\iota_{\tau}:\mathcal{M}_{d}(C)\rightarrow\mathcal{M}_{k}(C)$ are linear.
To consequently see these maps are order maps, it suf\/f\/ices that $\iota_{\tau}$ maps the positive cone
$\mathrm{M}_{d}(C)^{+}$ to $\mathrm{M}_{k}(C)^{+}$.
\end{Example}

This is clear from the following lemma.

\begin{Lemma}
\label{products}
Let~$S$ be a~$*$-semigroup, $\tau=(t_{1},\dots,t_{d})\in\mathcal{P}(d,k)$ and $\left[b_{i,j}\right]
\in\mathrm{M}_{r,d}(S)$.
There is $\left[c_{i,j}\right] \in\mathrm{M}_{r,k}(S)$, whose entries appear in $\left[b_{i,j}\right] $, such that
\begin{gather*}
\iota_{\tau}(\left[b_{i,j}\right]^{\ast}\left[b_{i,j}\right])=\left[c_{i,j}\right]^{\ast}\left[c_{i,j}\right].
\end{gather*}
In particular if $[n_{i}]\in\mathrm{M}_{1,d}(S)$ then $\iota_{\tau} ([n_{i}]^{\ast}[n_{j}])=[m_{i}]^{\ast}[m_{j}]$ for
some $[m_{i}]\in \mathrm{M}_{1,k}(S)$.
The entries of $[m_{i}]$ appear in $[n_{i}]$.
\end{Lemma}

\begin{proof}
For $1\leq i\leq r$ let the $r\times k$ matrix $\left[c_{i,j}\right] $ have~$i$-th row
\begin{gather*}
[b_{i1},\dots,b_{i1},b_{i2},\dots,b_{i2},\dots,b_{id},\dots,b_{id}],
\end{gather*}
where each element $b_{ij}$ appears repeated $t_{j}$ consecutive times.
In particular if $\left[b_{i,j}\right] $ is $[n_{i}]\in\mathrm{M}_{1,d}(S)$ then
\begin{gather*}
\iota_{\tau}([n_{i}]^{\ast}[n_{j}])=\iota_{\tau}([n_{i}^{\ast}n_{j}])=[m_{i}]^{\ast}[m_{j}]
\end{gather*}
with $[m_{i}]=[n_{1},\dots,n_{1},\dots,n_{d},\dots,n_{d}]\in\mathrm{M}_{1,k}(S)$, each $n_{j}$ repeated $t_{j}$ times.
\end{proof}

If $\beta:S\rightarrow T$ is a~$*$-map of $*$-semigroups (not necessarily a~homomorphism) then def\/ine
the~$k$-amplif\/ication
\begin{gather*}
\beta_{k}:\ \mathrm{M}_{k}(S)\rightarrow\mathrm{M}_{k}(T)\qquad\text{by}\quad \beta_{k}[n_{i,j}]=[\beta(n_{i,j})].
\end{gather*}
When convenient interpret $\beta_{k}([n_{i}])=[\beta(n_{i})]$ in this sense also.
The amplif\/ication maps $\beta_{k}$ behave well with the maps $\iota_{\tau}$,
\begin{gather*}
\beta_{k}\circ\iota_{\tau}=\iota_{\tau}\circ\beta_{d}
\qquad
\text{for}
\quad
\tau \in\mathcal{P}(d,k).
\end{gather*}

Depending on how the sets $\mathcal{M}_{k}(S)$ and $\mathcal{M}_{k}(T)$ are def\/ined $\beta_{k}$ may not def\/ine a~map
$\beta_{k}:\mathcal{M}_{k} (S)\rightarrow\mathcal{M}_{k}(T)$.
However if $\beta:S\rightarrow C$ is a~$*$-map to a~$C^*$-algebra~$C$ then $\beta_{k}:\mathcal{M}_{k}(S)\rightarrow
\mathcal{M}_{k}(C)$ is always def\/ined.

\begin{Definition}
A~$*$-map $\beta:S\rightarrow T$ of matricially ordered $*$-semi\-groups~$S$ and~$T$ is called a~$k$-order map if
$\beta_{k}:\mathcal{M}_{k}(S)\rightarrow \mathcal{M}_{k}(T)$ is def\/ined, and is an order map of partially ordered sets.
The map $\beta$ is a~complete order map if $\beta$ is a~$k$-order map for all $k\in\mathbb{N}$.
If~$T$ is a~$C^*$-algebra, a~complete order map $\beta$ which is a~representation is called a~complete order
representation.
\end{Definition}

If $\beta:B\rightarrow C$ is a~$*$-homomorphism of $C^*$-algebras then $\beta$ is a~complete order representation, and
any completely positive map of $C^*$-algebras is a~complete order map.

\begin{Definition}
Let $\beta:S\rightarrow T$ be a~$*$-map of a~$*$-semi\-group~$S$ to a~matricially ordered $*$-semi\-group $(T,\preceq,\mathcal{M})$.
The map $\beta_{k}$ has the Schwarz property for $k$, if
\begin{gather*}
\beta_{k}([n_{i}])^{\ast}\beta_{k}([n_{j}])\preceq\beta_{k}([n_{i}]^{\ast}[n_{j}])
\end{gather*}
in $\mathcal{M}_{k}(T)$ for $[n_{i}]\in\mathrm{M}_{1,k}(S)$.
Here $\beta_{k}([n_{i}])^{\ast}\beta_{k}([n_{j}])$ is the selfadjoint element $[\beta(n_{i})^{\ast}\beta(n_{j})]$ in~$\mathcal{M}_{k}(T)$.
\end{Definition}

A~composition of~$k$-order maps of matricially ordered $*$-semi\-groups is again a~$k$-order map as long as the
composition is def\/ined, so a~composition of complete order maps is a~complete order map.
If $\beta:R\rightarrow S$ and $\alpha:S\rightarrow T$ are $*$-maps of matricially ordered $*$-semi\-groups having the
Schwarz property for~$k$, and $\alpha$ is additionally a~$k$-order map, then their composition $\alpha\circ\beta$ also
has the Schwarz property for~$k$.
If $\sigma:S\rightarrow T$ is a~$*$-homomorphism of $*$-semi\-groups then
$\sigma_{k}([n_{i}])^{\ast}\sigma_{k}([n_{j}])=\sigma_{k}([n_{i}]^{\ast} [n_{j}])$ for $[n_{i}]\in\mathrm{M}_{1,k}(S)$
and the map $\sigma_{k}:\mathcal{M}_{k}(S)\rightarrow\mathcal{M}_{k}(T)$ (when def\/ined) automatically has the Schwarz
property.
Therefore, if $\beta:R\rightarrow S$ is a~$k$-order map with the Schwarz property and $\sigma:S\rightarrow T$ is
a~$*$-semi\-group homomorphism as well as a~$k$-order map, then $\sigma\circ\beta$ is a~$k$-order map with the Schwarz
property for~$k$.

\begin{Definition}
\label{cpd}
A~$*$-map $\beta$ from a~$*$-semi\-group~$S$ to a~$C^*$-algebra~$C$ is completely positive if $[\beta(n_{i}^{\ast}n_{j})]$
is positive in $\mathrm{M}_{k}(C)$ for any f\/inite set of elements $n_{1},\dots,n_{k}$ of~$S$.
\end{Definition}

Let $\beta$ be a~$*$-map of a~$*$-semi\-group~$S$ to a~$C^*$-algebra~$C$ with $\beta_{k}$ having the Schwarz pro\-perty for
all~$k$.
Then, since $\beta_{k}([n_{i}])^{\ast}\beta_{k}([n_{j}])$ is positive in the $C^*$-algebra $\mathrm{M}_{k}(C)$, the
matrix $\beta_{k}([n_{i}]^{\ast}[n_{j}])=[\beta (n_{i}^{\ast}n_{j})]$ must also be positive in $\mathrm{M}_{k}(C)$, and
$\beta$ must be a~completely positive map of the $*$-semi\-group~$S$.
In particular, if $\beta$ is a~$*$-homomorphism of a~$*$-semi\-group~$S$ to a~$C^*$-algebra~$C$, so also has the Schwarz
property for all~$k$, it must be completely positive.
Completely positive maps enable the construction (Lemma~\ref{H module}) of Hilbert modules $\mathcal{E}_{C}$ over the
$C^*$-algebra~$C$.

\begin{Remark}
If a~map $\beta:S\rightarrow C$ from a~$*$-semi\-group~$S$ to a~$C^*$-algebra~$C$ has the form $\beta=\rho\circ\gamma$,
where $\gamma$ is a~contractive $*$-representation of~$S$ in the $C^*$-algebra~$C$ and $\rho:C\rightarrow C$ is
a~completely positive map of $C^*$-algebras, then $\beta$ must be a~completely positive map.
This follows from
\begin{gather*}
[\beta(n_{i}^{\ast}n_{j})]=[\rho(\gamma(n_{i}^{\ast}n_{j}))]=\rho_{k}([\gamma(n_{i})]^{\ast}[\gamma(n_{j})]).
\end{gather*}
\end{Remark}

\subsection[$C^*$-correspondence]{$\boldsymbol{C^*}$-correspondence}

For $\beta:S\rightarrow C$ a~completely positive map from a~$*$-semi\-group~$S$ to a~$C^*$-algebra~$C$ form the algebraic
tensor product of complex vector spaces $X=\mathbb{C}[S]\otimes_{\rm alg}C$.
For two simple tensors $x=s\otimes c$, $y=t\otimes d$, with $s,t\in S,c,d$ in~$C$ def\/ine $\langle x,y\rangle
\in C$ by $\langle c,\beta(s^{\ast}t)d\rangle =c^{\ast}\beta(s^{\ast}t)d$ and extend to
a~sesqui-linear~$C$-valued map on $X\times X$.
The following proof uses a~standard approach for establishing Stinespring's dilation theorem (cf.~\cite{p}).

\begin{Lemma}
\label{H module}
For a~$C^*$-algebra~$C$ and $\beta:S\rightarrow C$ a~completely positive-map on a~$*$-semi\-group~$S$ let
$X=\mathbb{C}[S]\otimes_{\rm alg}C$ equipped with the sesqui-linear~$C$-valued map $\langle\,
,\, \rangle $ on $X\times X$ defined, as above, through $\beta$.
Then $\langle x,x\rangle \geq0$ in~$C$, $(x\in X)$.
\end{Lemma}

\begin{proof}
Let $x=\sum\limits_{i=1}^{k}s_{i}\otimes_{\rm alg}c_{i}\in X$, where $s_{i}\in S$ and $c_{i}\in C$.
Compute
\begin{gather*}
\langle x,x\rangle =\sum\langle c_{i},\beta(s_{i}^{\ast} s_{j})c_{j}\rangle =\sum
c_{i}^{\ast}\beta(s_{i}^{\ast}s_{j})c_{j}=\langle \overrightarrow{c},T\overrightarrow{c}\rangle,
\end{gather*}
where~$T$ is the positive matrix $[\beta(s_{i}^{\ast}s_{j})]$ in $\mathrm{M}_{k}(C)$ and $\overrightarrow{c}$ is the
vector $(c_{1},\dots,c_{k})$ in the Hil\-bert~$C$-module $\oplus_{1}^{k}C$.
The latter inner product in this Hilbert~$C$-module is equal to $\langle \sqrt{T}c,\sqrt{T}c\rangle $, which
is positive in~$C$ \cite[Lemma~4.1]{l}.
\end{proof}

For $x\in X$ and $d\in C$ we have $\langle xd,xd\rangle =d^{\ast}\langle x,x\rangle
d\leq\Vert \langle x,x\rangle \Vert d^{\ast}d$ (cf.~\cite{l}), so if $x\in N=\left\{x\in
X\,|\,\langle x,x\rangle =0\right\} $, or if $d=0$, then $xd\in N$.
In particular this allows the right action of~$C$ on~$X$ to become a~well def\/ined right action of~$C$ on the quotient
space $X/N$.
Set $\overline {x}d=\overline{xd}$ for $\overline{x}\in X/N$ and $d\in C$, where $\overline{x}$ denotes the class of~$x$
in $X/N$.
The quotient space $X/N$ is now an inner product~$C$-module, and its completion $\mathcal{E}_{C}$ becomes a~right
Hilbert~$C$-module (cf.~\cite{l}).

Since $S\mathbb{C}[S]\subseteq\mathbb{C}[S]$ there is a~left action of~$S$ on~$X$.
The Cauchy Schwarz inequality  \cite[Proposition~1.1]{l} for the semi-inner product on~$X$ shows
\begin{gather*}
\langle sx,sx\rangle =\langle s^{\ast}sx,x\rangle \leq\langle s^{\ast}sx,s^{\ast}sx\rangle^{1/2}\langle x,x\rangle^{1/2}
\qquad
\text{for}
\quad
x\in X,
\quad
s\in S
\end{gather*}
in~$C$.
Thus if $x\in N$ then $sx\in N$.
This yields a~well def\/ined left action of~$S$, and therefore also of $\mathbb{C}[S]$, on the inner product~$C$-module
$X/N$.
This is a~pre-$*$-representation (borrowing terminology from~\cite{f,f2}) on the pre-Hilbert module $X/N$.

If $X/N$ is already complete, for example if it is f\/inite dimensional, or if this left action is by bounded operators on
$X/N$, so extendable to an action on $\mathcal{E}_{C}$, there is a~$*$-representation of~$S$ by adjointable operators in
$\mathcal{L}(\mathcal{E}_{C})$.
However, it is not, a~priori, contractive.

In conclusion, a~completely positive-map $\beta$ on a~$*$-semi\-group~$S$ implies the existence of an associated Hilbert
module, along with a~$*$-action of the semigroup on a~dense submodule.
As there are examples of adjointable operators on dense submodules of a~Hilbert module which are not bounded, $\beta$
does not a~priory seem to yield a~bounded, let alone a~contractive, $*$-representation of the semigroup in this module
without some additional structural conditions.
Making use of a~matricial partial order on a~$*$-semi\-group~$S$ is one approach to naturally include some necessary
structure.

\subsection[Matricial orders on~$A$, $D_{1}$, and $D_{0}$]
{Matricial orders on~$\boldsymbol{A}$, $\boldsymbol{D_{1}}$, and $\boldsymbol{D_{0}}$}

The following proposition describes conditions for uniqueness of factorization, allowing for an approach to def\/ining
matricial partial orders on the $*$-semi\-groups under consideration.

\begin{Proposition}
\label{pain}
Let $\mathbf{n}=[n_{ij}]\in\mathcal{M}_{k}(A)$ with $\mathbf{n}=[w_{i}]^{\ast}[w_{j}]=[v_{i}]^{\ast}[v_{j}]$, where
$[w_{i}], [v_{i}]\in\mathrm{M}_{1,k}(A)$ has reduced entries in~$A$.
Assume $[w_{i}]$ satisfies one of the two conditions
\begin{gather*}
1)\quad(-1,1)[w_{i}]=[w_{i}]
\qquad \text{or}
\qquad
2) \quad(1,-1)[w_{i}]=[w_{i}].
\end{gather*}
Then $[v_{i}]=[w_{i}]$, or if $[v_{i}]\neq[w_{i}]$ then
$[v_{i}]=(1)[w_{i}]$ if $(-1,1)[w_{i}]=[w_{i}]$ or $[v_{i}]=(-1)[w_{i}]$ if $(1,-1)[w_{i}]=[w_{i}]$.
In particular $[v_{i}]$ must then also satisfy one of the two conditions $(-1,1)[v_{i}]=[v_{i}]$ or
$(1,-1)[v_{i}]=[v_{i}]$.

If $[w_{i}]$ does not satisfy one of these two conditions then
$[v_{i}]=[w_{i}]$.
\end{Proposition}

\begin{proof}
Write $w_{i}=(n_{i,0},\dots,n_{i,h_{i}})$ and $v_{i}=(m_{i,0},\dots,m_{i,l_{i}})$ reduced in~$A$.
With~$l$ the length function on~$A$, $l(w_{i}^{\ast}w_{i})$ is either $2l(w_{i})$ or $2l(w_{i})-2$.
The latter occurs exactly when $l(w_{i})\geq2$ and $\left\vert n_{i,0}\right\vert =1$.
Similarly for $v_{i}$.
Since $v_{i}^{\ast}v_{i}=w_{i}^{\ast}w_{i}$ for all~$i$, it follows that for each~$i$ one of following three
possibilities involving $l(v_{i})$ and $l(w_{i})$ must hold; $l(v_{i})=l(w_{i})-1$, occurring exactly if $l(w_{i})\geq2$
and $\left\vert n_{i,0}\right\vert =1$, $l(w_{i})=l(v_{i})-1$, occurring exactly if $l(v_{i})\geq2$ and $\left\vert
m_{i,0}\right\vert =1$, or $l(v_{i})=l(w_{i})$.
From $v_{i}^{\ast}v_{i}=w_{i}^{\ast}w_{i}$, it then follows that for each~$i$ one of the following three must hold:
$w_{i} =(n_{i,0})(v_{i})$, which occurs exactly if $\left\vert n_{i,0}\right\vert =1$; $v_{i}=(m_{i,0})w_{i}$, which
occurs exactly if $\left\vert m_{i,0}\right\vert =1$; or $v_{i}=w_{i}$.

For now assume that $(-1,1)[w_{i}]=[w_{i}]$, so $n_{i,0}\leq-1$ for all $i$; the case with $(1,-1)[w_{i}]=[w_{i}]$ is similar.
By the f\/inal sentence of the last paragraph, whenever $v_{j}\neq w_{j}$ then either
$w_{j}=(n_{j,0})(v_{j})=(-1)(v_{j})$, or $v_{j}=(m_{j,0})w_{j}=(1)w_{j}$.
Note that if $w_{j}=(-1)(v_{j})$ then it must be the case that $m_{j,0}\geq1$, since if $m_{j,0}\leq-1$ then
$n_{j,0}<m_{j,0}\leq-1$ contradicting $v_{j}^{\ast} v_{j}=w_{j}^{\ast}w_{j}$.
Suppose $v_{i}=w_{i}$ for some~$i$, and $v_{j}\neq w_{j}$ for some $j\neq i$.
If, for example, $v_{j}=(1)w_{j}$ then $m_{i,0}=1$ and
$w_{j}^{\ast}w_{i}=n_{ji}=v_{j}^{\ast}v_{i}=w_{j}^{\ast}(-1)w_{i}$, which cannot be the case as both $n_{i,0}$ and
$n_{j,0}\leq-1$.
On the other hand, if $w_{j}=(-1)(v_{j})$ then we have already noted $m_{j,0}\geq1$.
Therefore $v_{j}^{\ast}v_{i}=n_{ji}=w_{j}^{\ast}w_{i}=v_{j}^{\ast}(1)v_{i}$, which cannot be the case as $m_{i,0}$ (with $=n_{i,0}$)
and $-m_{j,0}$ are both $\leq-1$.
Thus either $[v_{i}]=[w_{i}]$, or for each~$i$ then $w_{i} =(-1)(v_{i})$ (so $m_{i,0}\geq1$) or $v_{i}=(1)w_{i}$
(which also implies $m_{i,0}\geq1$).
Therefore, if $[v_{i}]\neq[w_{i}]$ we must have $m_{i,0}\geq1$ for all~$i$, so $(1,-1)[v_{i}]=[v_{i}]$.

Next we show that if $[v_{i}]\neq[w_{i}]$ then $[v_{i}]=(1)[w_{i}]$.
Consider those possible~$j$ with $w_{j}=(-1)(v_{j})$.
We already know that $m_{j,0}\geq1$.
If $m_{j,0}=1$ then, since $w_{j}$ is reduced, $v_{j}$ must be $(1)$ and so $w_{j}=(-1,1)$.
Thus $v_{j}=(1)=(1)w_{j}$.
It follows that if $[v_{i}]\neq[w_{i}]$ then $v_{j}=(1)w_{j}$ for every~$j$, and $[v_{i}]=(1)[w_{i}]$.

The case remains when $[w_{i}]$ does not satisfy either of the conditions, so necessarily $k\geq2$.
Assume $[v_{i}]\neq[w_{i}]$, so there is a~$j$ with $v_{j}\neq w_{j}$.
Since $[w_{i}]$ does not satisfy either condition the sign of $n_{i,0}$ is not constant over~$i$, so if $n_{j,0}\leq-1$
say, then there is a~$g$ with $n_{g,0}\geq1$.
Since $v_{j}\neq w_{j}$ then $l(w_{j})\neq l(v_{j})$ and the f\/irst paragraph implies that either $w_{j}=(-1)(v_{j})$
with $m_{j,0}\geq1$, or $v_{j}=(1)w_{j}$ (Note that $v_{j}=(-1)w_{j}$ is not a~possibility since $l(w_{j})\neq l(v_{j})$.).
For each one of these two possibilities for the~$j$ coordinate, there are the three possibilities for the~$g$
coordinate: $w_{g}=v_{g}$; $w_{g}=(1)(v_{g})$ with $m_{g,0}\leq-1$ (since $n_{g,0}\geq1$); $v_{g}=(-1)w_{g}$.
Considering the $(j,g)$ coordinate of the matrix $\mathbf{n}$, so $w_{j}^{\ast}w_{g}=n_{jg}=v_{j}^{\ast}v_{g}$, we see
that none of these cases are possible.
For example, if $w_{g}=v_{g}$ then $v_{j}^{\ast}v_{g}=v_{j}^{\ast}(1)v_{g}$ when $w_{j}=(-1)(v_{j})$, or
$w_{j}^{\ast}w_{g}=w_{j}^{\ast}(-1)w_{g}$ when $v_{g}=(-1)w_{g}$, both not possible.
The other cases are similarly not possible.
Thus $[v_{i}]=[w_{i}]$.
\end{proof}

We are now in a~position to specify a~matricial order $(A,\preceq,\mathcal{M}_{k})$ for the ordered $*$-semi\-group~$A$.

\begin{Definition}
For the $*$-semi\-group~$A$ def\/ine the set
\begin{gather*}
\mathcal{M}_{k}(A)=\{[w_{i}]^{\ast}[w_{j}]\,|\,[w_{i}]\in\mathrm{M}_{1,k}(A),\;
[w_{i}]
~\text{has reduced entries in}~A\}.
\end{gather*}
\end{Definition}

For an element $\mathbf{n}\in\mathcal{M}_{k}(A)$ Proposition~\ref{pain} provides conditions that yield two such
factorizations $[w_{i}]^{\ast} [w_{j}]$, otherwise the factorization is unique.
Suppose $\mathbf{n} =[n_{ij}]=[w_{i}]^{\ast}[w_{j}]$, where each entry of $[w_{i}]\in \mathrm{M}_{1,k}(A)$ is a~reduced
word in~$A$.
We f\/irst consider the two cases occurring in Proposition~\ref{pain}; the case where $(-1,1)[w_{i}]=[w_{i}]$, and the
case $(1,-1)[w_{i}]=[w_{i}]$.

Consider the f\/irst case $(-1,1)[w_{i}]=[w_{i}]$.
If $(-1,1)[w_{i}]=[w_{i}]$ choose any $[w_{i}]_{+}\in \mathrm{M}_{1,k}(A)$ so that $[w_{i}]=(-1)[w_{i}]_{+}$.
Note that if $w_{i}=(-1)$ for some~$i$ then, as before, the only possibility is to set $(w_{i})_{+}=(1,-1)=(1)w_{i}$.
The element $(w_{i})_{+}=(1)w_{i}$ (reduced) is in fact always a~possible choice for $(w_{i})_{+}$ since
$w_{i}=(-1)(w_{i})_{+}$ is satisf\/ied under our assumption that $(-1,1)w_{i} =w_{i}$.
For this latter choice of $(w_{i})_{+}$ we have $(w_{i})_{+}^{\ast}(w_{i})_{+}=w_{i}^{\ast}w_{i}$.
However, if the choice $(w_{i})_{+}=(1)w_{i}$ is made for all~$i$ then
$[w_{i}]^{\ast}[w_{j}]=[w_{i}]_{+}^{\ast}[w_{j}]_{+}$ and one obtains the only other possible factorization of $\mathbf{n}=
[v_{i}]^{\ast}[v_{j}]$ with $[v_{i}]=(1)[w_{i}]$ and $(1,-1)[v_{i}]=[v_{i}]$ (Proposition~\ref{pain}).

Generally there are many possible choices for an element $[w_{i}]_{+} \in\mathrm{M}_{1,k}(A)$, since for each~$i$ there
are two choices.
However there is a~canonical `extreme' choice for $[w_{i}]_{+}$, where for each~$i$ the choice $(w_{i})_{+}$ satisf\/ies
$w_{i}^{\ast}w_{i}\prec(w_{i})_{+}^{\ast}(w_{i})_{+}$ whenever $w_{i}\neq(-1,\dots )$ as a~reduced word in~$A$.
If this is not possible for a~particular $w_{i}$, necessarily of the form $(-1,\dots)$, then $(w_{i})_{+}$ must satisfy
$w_{i}^{\ast}w_{i}=(w_{i})_{+}^{\ast} (w_{i})_{+}$.
Refer to the resulting $[w_{i}]_{+}$, when this extreme choice for each of the entries $(w_{i})_{+}$ is made whenever
possible, as the `extreme' $[w_{i}]_{+}$.

Now consider the second case $(1,-1)[w_{i}]=[w_{i}]$.
If $(1,-1)[w_{i}]=[w_{i}]$ choose any $[w_{i}]_{-}\in \mathrm{M}_{1,k}(A)$ so that $[w_{i}]
=(1)[w_{i}]_{-}$.
Again, if $w_{i}=(1)$ for some~$i$ then set $(w_{i})_{-}=(-1,1)$.
As before, the element $(w_{i})_{-}=(-1)w_{i}$ (reduced) is in fact always a~possible choice for $(w_{i})_{-}$ since
$w_{i}=(1)(w_{i})_{-}$ is satisf\/ied under our assumption that $(1,-1)w_{i}=w_{i}$.
For this latter choice of~$(w_{i})_{-}$ we have $(w_{i})_{-}^{\ast}(w_{i})_{-}=w_{i}^{\ast}w_{i}$.

Depending on which situation is considered, so $(-1,1)[w_{i}]=[w_{i}]$ or $(1,-1)[w_{i}]=[w_{i}]$, the matrix
\begin{gather*}
[w_{i}]^{\ast}[w_{j}]=[w_{i}]_{\pm}^{\ast}(\pm1,\mp1)[w_{j}]_{\pm}.
\end{gather*}

If neither of the two conditions on $[w_{i}]$, $(-1,1)[w_{i}]=[w_{i}]$ nor $(1,-1)[w_{i}]=[w_{i}]$, is satisf\/ied then
Proposition~\ref{pain} shows $\mathbf{n}=[w_{i}]^{\ast}[w_{j}]$ has a~unique factorization.

\begin{Definition}
Let $\mathbf{n}=[w_{i}]^{\ast}[w_{j}]\in\mathcal{M}_{k}(A)$.

If $(-1,1)[w_{i}]=[w_{i}]$ then choose $[w_{i}]_{+}$ with $[w_{i}]=(-1)[w_{i}]_{+}$.
Def\/ine $[w_{i}]^{\ast}[w_{j}]\preceq[w_{i}]_{+}^{\ast} [w_{j}]_{+}$.

If $(1,-1)[w_{i}]=[w_{i}]$
then choose $[w_{i}]_{-}$ with $[w_{i}]=(-1)[w_{i}]_{-}$.
Def\/ine $[w_{i}]^{\ast}[w_{j}]\preceq[w_{i}]_{-}^{\ast}[w_{j}]_{-}$.

In case neither of the two
conditions on $[w_{i}]$, $(-1,1)[w_{i}]=[w_{i}]$ nor $(1,-1)[w_{i}]=[w_{i}]$, is satisf\/ied $\mathbf{n}$ is def\/ined to be
a~maximal element in the partial order.

The partial order on $\mathcal{M}_{k}(A)$ is that generated by these basic relations.
Write $\mathbf{n}\prec\mathbf{m}$ if $\mathbf{n} \preceq\mathbf{m}$ but $\mathbf{n\neq m.}$
\end{Definition}

Proposition~\ref{pain} shows that if $\mathbf{n}=[w_{i}]^{\ast}[w_{j}]$ and the condition, $(-1,1)[w_{i}]=[w_{i}]$
or
$(1,-1)[w_{i}]$ $=[w_{i}]$, is satisf\/ied then it is always possible to write down a~second alternative
factorization of $\mathbf{n}$, where the other condition holds.
For example if $(-1,1)[w_{i}]=[w_{i}]$ then $\mathbf{n=}[v_{i}]^{\ast}[v_{j}]$ with $(1,-1)[v_{i}]=[v_{i}]$, where
$[v_{i}]=(1)[w_{i}]$.

If, for example, $[w_{i}]\in\mathrm{M}_{1,k}(A)$ with $(-1,1)[w_{i}]=[w_{i}]$ then choosing $(w_{i})_{+}=(1)w_{i}$ for
some~$i$ but not all~$i$ will generally ensure that neither $(-1,1)[w_{i}]_{+}=[w_{i}]_{+}$ nor
$(1,-1)[w_{i}]_{+}=[w_{i}]_{+}$ is satisf\/ied, so the matrix $[w_{i}]_{+}^{\ast}[w_{j}]_{+}$ is therefore maximal in
the partial order on $\mathcal{M}_{k}(A)$.
For $k\geq2$, it is straightforward that there are therefore an inf\/inite number of maximal elements in the matricial
partial ordering on $\mathcal{M}_{k}(A)$.
This points to an immediate structural dif\/ference with the partial ordering on $(A^{\rm sa},\preceq)$, which has exactly two maximal elements.

\begin{Proposition}
The $*$-semi\-group~$A$ with the partially ordered sets
\begin{gather*}
\mathcal{M}_{k}(A)=\{[w_{i}]^{\ast}[w_{j}]\,|\,[w_{i}]\in\mathrm{M}_{1,k}(A)\}
\end{gather*}
is a~matricially ordered $*$-semi\-group $(A,\preceq,\mathcal{M})$.
\end{Proposition}

\begin{proof}
We show $[n_{i}]_{\delta}^{\ast}\left[a_{i,j}\right] [n_{j}]_{\delta}\preceq[n_{i}]_{\delta}^{\ast}\left[b_{i,j}\right] [n_{j}]_{\delta}$
whenever $\left[a_{i,j}\right] \preceq\left[b_{i,j}\right]
$ in $\mathcal{M}_{k}(A)$ and $[n_{i}]\in$ $\mathrm{M}_{1,k}(A)$.
We may assume $\left[a_{i,j}\right] =[w_{i}]^{\ast}[w_{j}]$, where $[w_{i}]\in \mathrm{M}_{1,k}(A)$ with each $w_{i}$
a~reduced word and, for example, $(-1,1)[w_{i}]=[w_{i}]$ and $\left[b_{i,j}\right] =[w_{i}]_{+}^{\ast} [w_{j}]_{+}$.

Def\/ine
\begin{gather*}
[m_{i}]=[w_{1}n_{1},\dots,w_{k}n_{k}]\in\mathrm{M}_{1,k}(A).
\end{gather*}
Then
\begin{gather*}
[m_{i}]^{\ast}[m_{j}]=[n_{i}]_{\delta}^{\ast}([w_{i}]^{\ast}
[w_{j}])[n_{j}]_{\delta}=[n_{i}]_{\delta}^{\ast}\left[a_{i,j}\right] [n_{j}]_{\delta}.
\end{gather*}
We have
\begin{gather*}
(m_{i})_{+}=(w_{i}n_{i})_{+}=(w_{i})_{+}n_{i},
\end{gather*}
where if $w_{i}{}_{+}$ was chosen to be $(1,w_{i})$ then, to be consistent, set $(w_{i}n_{i})_{+}=(1,w_{i}n_{i})$.
In particular, if $w_{i}=(-1)$ then $(w_{i})_{+}$ must be $(1,-1)=(1,w_{i})$, so set $(m_{i})_{+}=(1,m_{i})$.
Then
\begin{gather*}
[n_{i}]_{\delta}^{\ast}\left[a_{i,j}\right] [n_{j}]_{\delta}
=[m_{i}]^{\ast}[m_{j}]\preceq[m_{i}]_{+}^{\ast}[m_{j}]_{+}
=([w_{i}]_{+}[n_{i}])^{\ast}[w_{j}]_{+}[n_{j}]=[n_{i}]_{\delta}^{\ast}\left[b_{i,j}\right]
[n_{j}]_{\delta}.
\end{gather*}

That the maps $i_{\tau}:$ $\mathcal{M}_{d}(A)\rightarrow\mathcal{M}_{k}(A)$ are order maps for $\tau\in\mathcal{P}(d,k)$
follows from Lemma~\ref{products} with $r=1$.
\end{proof}

\begin{Definition}
For the $*$-subsemigroup $D_{1}$ def\/ine $\mathcal{M}_{k}(D_{1})=\mathcal{M}_{k}(A)\cap\mathrm{M}_{k}(D_{1})$.
Therefore
\begin{gather*}
\mathcal{M}_{k}(D_{1})=\big\{[w_{i}]^{\ast}[w_{j}]\,|\,[w_{i}]\in \mathrm{M}_{1,k}(A),w_{i}^{\ast}w_{j}\in
D_{1}\ \text{for all}\  i,j\big\}.
\end{gather*}
\end{Definition}

{\sloppy The characterization of $D_{1}$ in Remark~\ref{b} provides a~check that each of the entries
$(w_{i})_{\pm}^{\ast}(w_{j})_{\pm}$ $\in D_{1}$ if the entries $(w_{i})^{\ast}(w_{j})\in D_{1}$.
Therefore
\begin{gather*}
[w_{i}]_{\pm}^{\ast}[w_{j}]_{\pm}\in\mathcal{M}_{k}(D_{1})\qquad\text{if}\quad \mathbf{n}=[n_{ij}]=[w_{i}]^{\ast}[w_{j}]\in\mathcal{M}_{k}(D_{1}).
\end{gather*}

}

\begin{Corollary}\label{mo}
The $*$-semi\-group $D_{1}$ with the sequence of partially ordered sets
\begin{gather*}
\mathcal{M}_{k}(D_{1})=\big\{[w_{i}]^{\ast}[w_{j}]\,|\,[w_{i}]\in \mathrm{M}_{1,k}(A),w_{i}^{\ast}w_{j}\in D_{1}\ \text{for all}\  i,j\big\}
\end{gather*}
is a~matricially ordered $*$-semi\-group $(D_{1},\preceq,\mathcal{M})$.
\end{Corollary}

In contrast to the situation for $k=1$, two dif\/ferent maximal elements, each strictly larger than a~given matrix,
actually may occur, as the following Example~\ref{hmm} shows.

\begin{Example}
\label{hmm}
Consider $\mathbf{n}=[\upsilon_{i}]^{\ast}[\upsilon_{j}]\in\mathcal{M}_{2}(D_{1})$, where
\begin{gather*}
[\upsilon_{i}]=[(-2,3),(-3,4)]\in\mathrm{M}_{1,2}(A).
\end{gather*}
The various choices in forming an element $[\upsilon_{i}]_{+}$ yield several elements in $\mathcal{M}_{2}(D_{1})$
strictly larger than $\mathbf{n}$, two of which are maximal in the partial order.
With
\begin{gather*}
[\upsilon_{i}]_{+}=[(-1,3),(-2,4)]
\end{gather*}
the resultant matrix $[\upsilon_{i}]_{+}^{\ast}[\upsilon_{j}]_{+}$ is not maximal in the partial order, while choosing
\begin{gather*}
[\upsilon_{i}]_{+}=[(-1,3),(1,-3,4)]\qquad\text{or}\qquad [\upsilon_{i}]_{+}=[(1,-2,3),(-2,4)]
\end{gather*}
one obtains two dif\/ferent maximal elements larger than $\mathbf{n}$.
These matrices are actually in $\mathrm{M}_{2}(\operatorname*{Irr}(D_{0}))^{\rm sa}$, so the set $\mathcal{M}_{k}(D_{1})$
properly contains the matrices $[n_{i}]^{\ast}[n_{i}]=[n_{i}^{\ast}n_{j}]$ and
$[n_{i}]^{\ast}a[n_{i}]$ $=[n_{i}^{\ast}an_{j}]$ ($a\in D_{1}^{\rm sa}$ and $[n_{i}]\in\mathrm{M}_{1,k}(D_{1})$).

The other factorization of $\mathbf{n}$ using $(1)[\upsilon_{i}]$ (cf.\ Proposition~\ref{pain}) does not yield any other
matrices immediately larger than $[\upsilon_{i}]^{\ast}[\upsilon_{j}]$.
However, it is possible that each of the two facto\-ri\-za\-tions of an $\mathbf{n}\in\mathcal{M}_{k}(D_{1})$ may lead to
dif\/ferent maximal elements strictly larger than $\mathbf{n}$.
For example, the elements
\begin{gather*}
[\upsilon_{i}]=[(-1,2,-5,6),(-3,5)]\in\mathrm{M}_{1,2}(A)
\end{gather*}
and
\begin{gather*}
(1)[\upsilon_{i}]=[(2,-5,6),(1,-3,5)]\in\mathrm{M}_{1,2}(A)
\end{gather*}
both yield $\mathbf{n=}[\upsilon_{i}]^{\ast}[\upsilon_{j}]\in\mathrm{M}_{2}(\operatorname*{Irr}(D_{0}))^{\rm sa}$.
Using the factorization $[\upsilon_{i}]^{\ast}[\upsilon_{j}]$ we have $\mathbf{n}\preceq[w_{i}]^{\ast}[w_{j}]$,
where
\begin{gather*}
[w_{i}]=[(2,-5,6),(-2,5)]\in\mathrm{M}_{1,2}(A).
\end{gather*}
Also $\mathbf{n}\preceq[u_{i}]^{\ast}[u_{j}]$ with
\begin{gather*}
[u_{i}]=[(1,-5,6),(-3,5)]\in\mathrm{M}_{1,2}(A).
\end{gather*}
Note that these matrices are also in $\mathrm{M}_{2}(\operatorname*{Irr} (D_{0}))^{\rm sa}$.
\end{Example}

The example above is an illustration of the following structure for the partial order def\/ined on
$\mathcal{M}_{k}(D_{1})$.

\begin{Proposition}
Let $\mathbf{n}=[n_{ij}]\in\mathcal{M}_{k}(D_{1})$.
Consider two different factorizations $\mathbf{n}=[w_{i}]^{\ast}[w_{j}]$ and $\mathbf{n} =[v_{i}]^{\ast}[v_{j}]$, where
$[w_{i}]$, $[v_{i}]\in\mathrm{M}_{1,k}(A)$, $(-1,1)[w_{i}]=[w_{i}]$ and $(1,-1)[v_{i}]=[v_{i}]$.
With $[w_{i}]_{+}$ and $[v_{i}]_{-}$ the extreme choices one of the following holds:
\begin{enumerate}\itemsep=0pt
\item[$a)$]
$\mathbf{n=}[w_{i}]_{+}^{\ast}[w_{j}]_{+}$ or $\mathbf{n=}[v_{i}]_{-}^{\ast}[v_{j}]_{-}$,

\item[$b)$]
$[w_{i}]_{+}^{\ast}[w_{j}]_{+}$ and $[v_{i}]_{-}^{\ast}[v_{j}]_{-}$ are distinct maximal elements strictly larger than
$\mathbf{n}$.
This can only occur if $k\geq2$.
\end{enumerate}
\end{Proposition}

\begin{proof}
Let $[w_{i}]$, $[v_{i}]\in\mathrm{M}_{1,k}(A)$ with $w_{i}=(n_{i,0},\dots,n_{i,h_{i}})$ and
$v_{i}=(m_{i,0},\dots,m_{i,l_{i}})$ reduced in~$A$.
The hypotheses imply that for all~$i$ each f\/irst entry of $w_{i}$, $n_{i,0} \leq-1$, while for $v_{i}$, $m_{i,0}\geq1$.
If $n_{i,0}$ $=-1$ for all~$i$, then $w_{i}=(-1)(n_{i,1},\dots,n_{i,h_{i}})$, which is also $(-1)w_{i}{}_{+}$ by
def\/inition.
Thus
\begin{gather*}
\mathbf{n}=[w_{i}^{\ast}w_{j}]=[w_{i}{}_{+}^{\ast}(1,-1)w_{j+}]=[w_{i+}^{\ast}w_{j+}]=[w_{i}]_{+}^{\ast}[w_{j}]_{+}.
\end{gather*}
We used that $[w_{i}]_{+}$ satisf\/ies $(1,-1)[w_{i}]_{+}=[w_{i}]_{+}$ since $n_{i,1}\geq1$ for all~$i$.
The situation with each f\/irst entry of $v_{i}=1$ follows similarly.
Note if each $w_{i}=(-1)$ then $v_{i}=(1,-1)$ (similarly if each $v_{i}=(1)$ then $w_{i}=(-1,1)$) for all~$i$.

Proposition~\ref{pain} applied to the hypothesis implies $[v_{i}]=(1)[w_{i}]$ and $[w_{i}]=(-1)[v_{i}]$.
If the f\/irst entries of $w_{i}$, $n_{i,0}\leq-2$ for all~$i$, then $[v_{i}]=(1)[w_{i}]$ implies that the f\/irst entries
$m_{i,0}$ of $v_{i}$ all equal~$1$, a~case already dealt with.

We have reduced to the case, for $k\geq2$, that there are~$i$ and~$j$ with $n_{i,0}=-1$, and $n_{j,0}\leq-2$.
Thus the~$i$-th coordinate of $[w_{i}]_{+}$ is positive while the~$j$-th coordinate of $[w_{i}]_{+}$ is negative, and
Proposition~\ref{pain} implies there is no factorization of $[w_{i}]_{+}^{\ast} [w_{j}]_{+}$ satisfying either of the
conditions needed to implement a~basic step in the partial ordering.
Therefore $[w_{i}]_{+}^{\ast}[w_{j}]_{+}$ is a~maximal element.

Similarly the alternative
factorization $\mathbf{n}=[v_{i}]^{\ast}[v_{j}]$ yields a~second maximal element
$[v_{i}]_{-}^{\ast}[v_{j}]_{-}$ strictly larger than $\mathbf{n.}$ By considering the~$i$ and~$j$ elements of
$[w_{i}]_{+}$ and $[v_{i}]_{-}$ one can check that $[w_{i}]_{+}^{\ast} [w_{j}]_{+}$ and $[v_{i}]_{-}^{\ast}[v_{j}]_{-}$
have dif\/ferent $(i,j)$ entries, so are distinct.
\end{proof}

Some further structure for the partial order on $\mathcal{M}_{k}(D_{1})$ is readily available, namely there are at least
two elements in $\mathcal{M}_{k}(D_{1})$ immediately less than any given element.

\begin{Proposition}
Let $\mathbf{n}=[n_{ij}]=[v_{i}]^{\ast}[v_{j}]\in\mathcal{M}_{k}(D_{1})$, where $[v_{i}]\in\mathrm{M}_{1,k}(A)$.
There are exactly two elements $[a_{i}]^{\ast}[a_{j}]$ and $[b_{i}]^{\ast}[b_{j}]\in\mathcal{M}_{k}(D_{1})$ with
\begin{gather*}
[a_{i}]^{\ast}[a_{j}]\prec\mathbf{n},
\qquad
[b_{i}]^{\ast}[b_{j}]\prec\mathbf{n},
\end{gather*}
and where
$[a_{i}]_{+}$ and $[b_{i}]_{-}$ are choices with
$[a_{i}]_{+}^{\ast}[a_{j}]_{+}=[b_{i}]_{-}^{\ast}[b_{j}]_{-}=\mathbf{n.}$
\end{Proposition}

\begin{proof}
To form $[a_{i}]_{+}$ and $[b_{i}]_{-}$ we implicitly assume $(-1,1)[a_{i}]=[a_{i}]$, $(1,-1)[b_{i}]=[b_{i}]$.
Let $[w_{i}]=(-1)[v_{i}]$, so $(-1,1)[w_{i}]=[w_{i}]$.
If $[w_{i}]^{\ast}[w_{j}]\neq[v_{i}]^{\ast}[v_{j}]$ set $[a_{i}]=[w_{i}]$.
With $[a_{i}]_{+}$ def\/ined by $[a_{i}]=(-1)[a_{i}]_{+}$ we have $[a_{i}]_{+}=[v_{i}]$ is a~choice for $[a_{i}]_{+}$ and
$[a_{i}]_{+}^{\ast}[a_{j}]_{+}=[v_{i}]^{\ast}[v_{j}]=\mathbf{n}$.
If $[w_{i}]^{\ast}[w_{j}]=[v_{i}]^{\ast}[v_{j}]$ set $[a_{i}]=(-2)[v_{i}]$.
Then $[a_{i}]_{+}=(-1)[v_{i}]=[w_{i}]$ is a~valid choice for $[a_{i}]_{+}$ and
$[a_{i}]_{+}^{\ast}[a_{j}]_{+}=[w_{i}]^{\ast}[w_{j}]=[v_{i}]^{\ast}[v_{j}]$.
We have
\begin{gather*}
(-1,1)[a_{i}]=[a_{i}]\qquad\text{and}\qquad [a_{i}]^{\ast}[a_{j}]\prec\mathbf{n=} [a_{i}]_{+}^{\ast}[a_{j}]_{+}.
\end{gather*}

To obtain a~second element $[b_{i}]^{\ast}[b_{j}]\in\mathcal{M}_{k}(D_{1})$ set $[u_{i}]=(1)[v_{i}]$, so
$(1,-1)[u_{i}]=[u_{i}]$.
With a~similar argument as above set $[b_{i}]=[u_{i}]$ if $[u_{i}]^{\ast}[u_{j}]\neq[v_{i}]^{\ast}[v_{j}]$, or
$[b_{i}]=(2)[v_{i}]$ if $[u_{i}]^{\ast} [u_{j}]=[v_{i}]^{\ast}[v_{j}]$.
Then
\begin{gather*}
(1,-1)[b_{i}]=[b_{i}]\qquad\text{and}\qquad [b_{i}]^{\ast}[b_{j}]\prec\mathbf{n=} [b_{i}]_{-}^{\ast}[b_{j}]_{-}.
\end{gather*}

Note that if $\mathbf{n}$ had two dif\/fering factorizations (the precise conditions of this are detailed in
Proposition~\ref{pain}) then for each of these two factorizations the two elements obtained in the manner described
above eventually result in the same set of two elements.
Also, using Proposition~\ref{pain} for the various scenarios implies $[a_{i}]^{\ast}[a_{j}]\neq[b_{i}]^{\ast}[b_{j}]$.

To show that there are exactly two such elements, assume that $[c_{i}]\in\mathrm{M}_{1,k}(A)$ with
\begin{gather*}
[c_{i}]^{\ast}[c_{j}]\prec\mathbf{n}=[c_{i}]_{+}^{\ast}[c_{j}]_{+}.
\end{gather*}
By def\/inition $[c_{i}]=(-1)[c_{i}]_{+}$ for any choice of $[c_{i}]_{+}$.
Since $[a_{i}]_{+}^{\ast}[a_{j}]_{+}=\mathbf{n=}[c_{i}]_{+}^{\ast}[c_{j}]_{+}$, by considering cases, including the
possibility that $\mathbf{n}$ has two factorizations by reduced words $\mathrm{M}_{1,k}(A)$, Proposition~\ref{pain}
implies that $[a_{i}]_{+}=[c_{i}]_{+}$, so $[a_{i}]=(-1)[a_{i}]_{+}=(-1)[c_{i}]_{+}=[c_{i}]$ and therefore
$[a_{i}]^{\ast}[a_{j}]=[c_{i}]^{\ast}[c_{j}]$.
We illustrate this argument in one case.
If $(-1,1)[a_{i}]_{+}=[a_{i}]_{+}$ but $[a_{i}]_{+}\neq[c_{i}]_{+}$ then Proposition~\ref{pain} implies
$[c_{i}]_{+}=(1)[a_{i}]_{+}$.
Thus $[c_{i}]=(-1)[c_{i}]_{+}=(-1,1)[a_{i}]_{+}=[a_{i}]_{+}$ and
$[c_{i}]^{\ast}[c_{j}]=[a_{i}]_{+}^{\ast}[a_{j}]_{+}=\mathbf{n}$, contradicting the assumption on
$[c_{i}]^{\ast}[c_{j}]$.
\end{proof}

\begin{Example}
Consider $\mathbf{n=}[w_{i}]^{\ast}[w_{j}]$, where
\begin{gather*}
[w_{i}]=[(2,-5,6),(-2,5)]\in\mathrm{M}_{1,2}(A)
\end{gather*}
arose in the previous example.
Proceeding along the lines of the last proposition form
\begin{gather*}
[a_{i}]=(-1)[w_{i}]=[(-1,2,-5,6),(-3,5)]
\end{gather*}
since $[a_{i}]^{\ast}[a_{j}]\neq[w_{i}]^{\ast}[w_{j}]$.
We see $[a_{i}]^{\ast}[a_{j}]\prec\mathbf{n}$, in fact $[a_{i}]$ is the element $[v_{i}]$ from the preceding example.

To obtain the second element just below $\mathbf{n}$ set
\begin{gather*}
[b_{i}]=(1)[w_{i}]=[(3,-5,6),(1,-2,5)]
\end{gather*}
(since $[b_{i}]^{\ast}[b_{j}]\neq[w_{i}]^{\ast}[w_{j}]$) and then $[b_{i}]^{\ast}[b_{j}]\prec\mathbf{n}$.

Generally there may be elements other than $\mathbf{n}$ immediately above $[a_{i}]^{\ast}[a_{j}]$ or
$[b_{i}]^{\ast}[b_{j}]$, as is shown in the preceding example for $[a_{i}]^{\ast}[a_{j}]$ ($=[v_{i}]^{\ast}[v_{j}]$ in
that case).
\end{Example}

The matricial partial order on $D_{1}$ restricts to a~matricial partial ordering on the unital $*$-semi\-group $D_{0}$,
where
\begin{gather*}
\mathcal{M}_{k}(D_{0})=\mathcal{M}_{k}(D_{1})\cap\mathrm{M}_{k}(D_{0}).
\end{gather*}
This follows by noting that if $[w_{i}]\in\mathrm{M}_{1,k}(A)$ with $[w_{i}]^{\ast}[w_{j}]\in\mathcal{M}_{k}(D_{0})$, so
that $w_{i}^{\ast} w_{j}\in D_{0}$ for all pairs $(i,j)$, then whenever $(-1,1)[w_{i}]=[w_{i}]$ or
$(1,-1)[w_{i}]=[w_{i}]$ so that elements $[w_{i}]_{\pm}^{\ast}[w_{j}]_{\pm}$ can actually be formed, then
$[w_{i}]_{\pm}^{\ast}[w_{j}]_{\pm}$ is again in $\mathcal{M}_{k}(D_{0})$.
That this is the case follows again from the description of $D_{0}$ as $\left\{n\in
A_{+}^{0}\,|\,\sigma_{r}(n)\leq0~\text{for all}~r\geq0\right\}$ (Proposition~\ref{D 1}) along with a~more detailed
checking of the various possibilities.

As in the order situation, there are examples of complete order maps.

\begin{Proposition}
\label{0 to 1}
The $*$-homomorphism $\omega:D_{0}\rightarrow D_{1}$ is a~complete order map of matricially ordered $*$-semi\-groups.
\end{Proposition}

\begin{proof}
Since $\omega(a)=(1)a(-1)$ for $a\in D_{0}$ it follows as before (cf.\ Proposition~\ref{before}), that the basic one
step relation that generates the partial orders on the $\mathcal{M}_{k}(D_{0})$ are preserved by this map.
\end{proof}

\begin{Remark}
\label{comp order}
Not every order representation of a~completely ordered $*$-semi\-group $(S,\preceq,\mathcal{M})$ is a~complete order
representation.
The unital completely ordered $*$-semi\-group $D_{0}$ provides an example.

With~$F$ the $*$-closed set
$\operatorname*{Irr}\nolimits_{+}(A^{0})\setminus(-1,1)$ of $D_{0}$, def\/ine a~$*$-map of~$F$ to a~set of elements in the unit
ball of a~unital $C^*$-algebra by mapping the ordered chain of three elements $\{(-k,k)\,|\,4\leq k\leq2\}$ in $F^{\rm sa}$
to an ordered chain of 3 nonzero positive elements in a~unital $C^*$-algebra~$C$.
Note that this requires that all powers of these selfadjoint elements have nonzero images in~$C$, and therefore all
elements in $C^{\rm sa}$ larger than these powers are also nonzero.
Thus
\begin{gather*}
(-4,4)^{2}\preceq(-4,3,-3,4)\preceq(-4,2,-2,4)\preceq(-4,4)
\end{gather*}
and
\begin{gather*}
(-3,3)^{2}\preceq(-3,2,-2,3)\preceq(-3,3)
\end{gather*}
have nonzero images in $C^{\rm sa}$.
However, we can stipulate that in~$C$ the images of $(-4,3,-3,4)$ and $(-4,2,-2,4)$ agree while the image of $(-3,3)$ is
strictly greater than the image of $(-3,2,-2,3)$.
Those elements of~$F$ that are not forced to be nonzero can all be mapped to $0$ in~$C$.
In particular $(-k,k)$ for $k\leq5$, as well as $(-4,3,-2,3)$ can be mapped to zero.
This def\/ines a~unital $*$-semi\-group homomorphism of $D_{0}$ to~$C$ which is an order map, so an order representation~$\pi$.
The map~$\pi_{2}$ is not, however, a~$2$-order map of $D_{0}$ to~$C$.
With $[w_{i}]=((-2,3),(-3,4))\in \mathrm{M}_{1,2}(A)$, then $[w_{i}]_{+}=((-1,3),(-2,4))$
and so
$[w_{i}]^{\ast}[w_{j}]\preceq[w_{i}]_{+}^{\ast}[w_{j}]_{+}$ in $\mathcal{M}_{2}(D_{0})$.
The image $\pi_{2}([w_{i}]^{\ast}[w_{j}])$ is
\begin{gather*}
\left[
\begin{matrix}
\pi(-3,2,-2,3) & \pi(-3,2,-3,4)
\\
\pi(-4,3,-2,3) & \pi(-4,3,-3,4)
\end{matrix}
\right]
\end{gather*}
while $\pi_{2}([w_{i}]_{+}^{\ast}[w_{j}]_{+})$ is the element
\begin{gather*}
\left[
\begin{matrix}
\pi(-3,3) & \pi(-4,4)
\\
\pi(-4,4) & \pi(-4,2,-2,4)
\end{matrix}
\right]
\end{gather*}
of $\mathrm{M}_{2}(C)$.
Since $\pi(-4,4)\neq0$ their dif\/ference
\begin{gather*}
\pi_{2}([w_{i}]_{+}^{\ast}[w_{j}]_{+})-\pi_{2}([w_{i}]^{\ast}[w_{j}])=\left[
\begin{matrix}
\pi(-3,3)-\pi(-3,2,-2,3) & \pi(-4,4)
\\
\pi(-4,4) & 0
\end{matrix}
\right]
\end{gather*}
cannot be positive in $\mathrm{M}_{2}(C)$ (cf.~\cite[Exercise~3.2]{p}).
\end{Remark}

The next result is crucial in obtaining Hilbert modules via the matricially ordered $*$-semi\-group $(D_{1},\preceq, \mathcal{M})$.

\begin{Proposition}
The $*$-map $\alpha:D_{1}\rightarrow D_{1}$ is a~complete order map of $D_{1}$ satisfying the Schwarz inequality for
each $k\in\mathbb{N}$.
\end{Proposition}

\begin{proof}
With $[n_{i,j}]=[w_{i}]^{\ast}[w_{j}]\in\mathcal{M}_{k}(D_{1})$ (where $[w_{i}]\in\mathrm{M}_{1,k}(A)$) it is clear that
$\alpha_{k}([w_{i}]^{\ast}[w_{j}])=[\alpha(n_{i,j})]=([w_{i}](1))^{\ast}[w_{j}](1)\in\mathcal{M}_{k}(D_{1})$.
If either of the conditions initiating a~basic order relation on $[w_{i}]$ is satisf\/ied then this remains the case when
considering $[\alpha(n_{i,j})]$.
The basic process to obtain $[w_{i}]_{\pm}^{\ast} [w_{j}]_{\pm}$ given $[w_{i}]^{\ast}[w_{j}]$ will also yield $\alpha
_{k}([w_{i}]_{\pm}^{\ast}[w_{j}]_{\pm})$ given $[\alpha(n_{i,j})]=([w_{i}](1))^{\ast}[w_{j}](1)$, and $\alpha$ is a~complete order map.

For $[n_{i}]\in\mathrm{M}_{1,k}(D_{1})$, $[n_{i}]^{\ast}(1,-1)[n_{i}]\preceq[n_{i}]^{\ast}[n_{i}]$ in
$\mathcal{M}_{k}(D_{1})$.
Applying the order map $\alpha_{k}$ and noting Remark~\ref{b} the Schwarz inequality follows.
\end{proof}

We now consider Hilbert modules arising from this context for the ordered $*$-semi\-group $(D_{1},\preceq)$.
If $\sigma:D_{1}\rightarrow C$ is an order representation of the ordered $*$-semi\-group $(D_{1},\preceq)$ in
a~$C^*$-algebra~$C$ then the composition
\begin{gather*}
\beta=\sigma\circ\alpha: \ D_{1}\rightarrow C
\end{gather*}
is an order map which also satisf\/ies the Schwarz inequality for all $k\in\mathbb{N}$.
In particular the comments preceding Def\/inition~\ref{cpd} imply $\beta:D_{1}\rightarrow C$ is completely positive and
therefore (Lemma~\ref{H module}) the space $X=\mathbb{C}[D_{1}]\otimes_{\rm alg}C$ along with the sesqui-linear map
$\langle\,,\, \rangle$ yields a~Hilbert module $\mathcal{E}_{C}$.

\begin{Lemma}
\label{bnd}
Let $\sigma:D_{1}\rightarrow C$ be an order representation of the ordered $*$-semi\-group $(D_{1},\preceq)$ in
a~$C^*$-algebra~$C$.
If $\beta:D_{1}\rightarrow C$ is the order map $\sigma\circ\alpha$ defined on the ordered $*$-semi\-group
$(D_{1},\preceq)$ then
\begin{gather*}
\langle (1,-1)x,(1,-1)x\rangle \leq\langle x,x\rangle
\qquad
\text{and}
\qquad
\langle(-1,1)x,(-1,1)x\rangle \leq\langle x,x\rangle
\end{gather*}
for $x\in X=\mathbb{C}[D_{1}]\otimes_{\rm alg}C$.
\end{Lemma}

\begin{proof}
Let $x=\sum\limits_{i=1}^{k}s_{i}\otimes_{\rm alg}c_{i}\in X$, where $s_{i}\in D_{1}$ and $c_{i}\in C$.
Then
\begin{gather*}
\langle (1,-1)x,(1,-1)x\rangle =\sum\langle c_{i},\beta (s_{i}^{\ast}(1,-1)s_{j})c_{j}\rangle
=\sum c_{i}^{\ast}\sigma (\alpha(s_{i}^{\ast})\alpha(s_{j}))c_{j}
\end{gather*}
by Remark~\ref{b}.
This equals $\langle \overrightarrow{c},T\overrightarrow{c}\rangle $ where~$T$ is the positive matrix
$\beta_{k}([s_{i}])^{\ast}\beta_{k}([s_{j}])$ in $\mathrm{M}_{k}(C)$ with $\overrightarrow{c}\in\oplus_{1}^{k}C$,
$[s_{i}]\in\mathrm{M}_{1,k}(D_{1})$ appropriately def\/ined.
Since $\beta$ satisf\/ies the Schwarz inequality for all~$k$,
\begin{gather*}
\beta_{k}([s_{i}])^{\ast}\beta_{k}([s_{j}])\preceq\beta_{k}([s_{i}]^{\ast}[s_{j}])
\end{gather*}
in $\mathrm{M}_{k}(C)$, so
\begin{gather*}
\langle \overrightarrow{c},T\overrightarrow{c}\rangle \leq\langle
\overrightarrow{c},\beta_{k}([s_{i}]^{\ast}[s_{j}])\overrightarrow{c}\rangle =\sum c_{i}^{\ast}\beta(s_{i}^{\ast}
s_{j})c_{j}=\sum\langle c_{i},\beta(s_{i}^{\ast}s_{j})c_{j}\rangle =\langle x,x\rangle.
\end{gather*}

For the second inequality note
\begin{gather*}
\langle (-1,1)x,(-1,1)x\rangle
=\sum\limits_{k}\sum\limits_{J_{k}}c_{i}^{\ast}\beta(s_{i}^{\ast}(-1,1)s_{j})c_{j},
\end{gather*}
where $J_{k}$, $1\leq k\leq$ $4$, are the index sets
\begin{gather*}
J_{1}=\{(i,j)\,|\,(1,-1)s_{i}=s_{i},~\text{and}~(1,-1)s_{j}=s_{j}\},
\\
J_{2}=\{(i,j)\,|\,(-1,1)s_{i}=s_{i},~\text{and}~(-1,1)s_{j}=s_{j}\},
\\
J_{3}=\{(i,j)\,|\,(-1,1)s_{i}=s_{i},~\text{and}~(1,-1)s_{j}=s_{j}\},
\\
J_{4}=\{(i,j)\,|\,(1,-1)s_{i}=s_{i},~\text{and}~(-1,1)s_{j}=s_{j}\}.
\end{gather*}
An inspection shows that the sums over $J_{k}$, where $k=2,3,4$, satisfy
\begin{gather*}
\sum\limits_{J_{k}}c_{i}^{\ast}\beta(s_{i}^{\ast}(-1,1)s_{j})c_{j}=\sum\limits_{J_{k}}
c_{i}^{\ast}\beta(s_{i}^{\ast}s_{j})c_{j}.
\end{gather*}
For the sum over $J_{1}$ we have
\begin{gather*}
\sum\limits_{J_{1}}c_{i}^{\ast}\beta(s_{i}^{\ast}(-1,1)s_{j})c_{j}=\sum\limits_{J_{1}}
c_{i}^{\ast}\beta(s_{i}^{\ast}(1,-1)(-1,1)(1,-1)s_{j})c_{j}.
\end{gather*}
However the terms
\begin{gather*}
\beta(s_{i}^{\ast}(1,-1)(-1,1)(1,-1)s_{j})   =\sigma\circ\alpha(s_{i}^{\ast}(1,-1)(-1,1)(1,-1)s_{j})
\\
\qquad
=\sigma(\alpha(s_{i}^{\ast})\alpha(-1,1)\alpha(s_{j}))   =\beta(s_{i}^{\ast})\beta(-1,1)\beta(s_{j})
\end{gather*}
for each $i,j\in J_{1}$.
We have
\begin{gather*}
\sum\limits_{J_{1}}c_{i}^{\ast}\beta(s_{i}^{\ast}(-1,1)s_{j})c_{j}=\sum\limits_{J_{1}}
c_{i}^{\ast}\beta(s_{i}^{\ast})\beta(-1,1)\beta(s_{j})c_{j}.
\end{gather*}
while Propositions~\ref{gen square},~\ref{square} imply $\beta(-1,1)$ is in the unit ball of the positive elements of
$C^{\ast}(D_{1},\preceq)$.
Thus the positive matrix
\begin{gather*}
\beta_{k}([s_{i}])^{\ast}\beta(-1,1)\beta_{k}([s_{j}])\leq\beta_{k} ([s_{i}])^{\ast}\beta_{k}([s_{j}])
\end{gather*}
in the appropriate matrix algebra over the $C^*$-algebra~$C$, and
\begin{gather*}
\sum\limits_{J_{1}}c_{i}^{\ast}\beta(s_{i}^{\ast}(1,-1)(-1,1)(1,-1)s_{j})c_{j}
\leq\sum\limits_{J_{1}}c_{i}^{\ast}\beta(s_{i})^{\ast}\beta(s_{j})c_{j}\leq
\sum\limits_{J_{1}}c_{i}^{\ast}\beta(s_{i}^{\ast}s_{j})c_{j}.
\end{gather*}
The latter inequality follows from the Schwarz inequality.
Therefore
\begin{gather*}
\langle (-1,1)x,(-1,1)x\rangle \leq\sum\limits_{J_{1}}c_{i}^{\ast}
\beta(s_{i}^{\ast}s_{j})c_{j}+\sum\limits_{k=2}^{4}\left(\sum\limits_{J_{k}}c_{i}^{\ast}\beta(s_{i}^{\ast}s_{j})c_{j}\right)=\langle x,x\rangle.\tag*{\qed}
\end{gather*}
\renewcommand{\qed}{}
\end{proof}

\begin{Proposition}
Let $\sigma:D_{1}\rightarrow C$ be an order representation of $(D_{1},\preceq)$ to a~$C^*$-algebra~$C$ and
$\beta=\sigma\circ\alpha:D_{1}\rightarrow C$ the order map on $D_{1}$ defining the Hilbert module $\mathcal{E}_{C}$.
The maps
\begin{gather*}
x\rightarrow(-1,1)x
\qquad
\text{and}
\qquad
x\rightarrow(1,-1)x
\end{gather*}
on the dense space $X=\mathbb{C}[D_{1}]\otimes_{\rm alg}C$ are bounded and therefore define positive elements in
$\mathcal{L}(\mathcal{E}_{C})$
of norm~$1$.
\end{Proposition}

\begin{proof}
The lemma above shows these maps are bounded in norm by~$1$, so can be extended to bounded operators on
$\mathcal{E}_{C}$.
They are both clearly adjointable maps on $\mathcal{E}_{C}$, in fact selfadjoint idempotents, and nonzero, so of norm~$1$.
\end{proof}

\subsection[The universal $C^*$-algebra $C^{\ast}(S,\preceq,\mathcal{M})$]
{The universal $\boldsymbol{C^*}$-algebra $\boldsymbol{C^{\ast}(S,\preceq,\mathcal{M})}$}

Without resorting to maps with additional order structure this is where progress stalls.
For a~general element $a\in D_{1}$ we would like to conclude that the basic map $x\rightarrow ax$ on~$X$ is bounded, and
adjointable, but this is unknown for our context.
One problem is that it is not clear that the order map $\beta=\sigma\circ\alpha$ is a~complete order map on $D_{1}$ if
$\sigma$ is only an order representation (Remark~\ref{comp order}).

In order to deal with these dif\/f\/iculties we introduce the universal $C^*$-algebra $C^{\ast}((S,\preceq,\mathcal{M}))$.

\begin{Definition}
Given $(S,\preceq,\mathcal{M})$ a matricially ordered $*$-semi\-group $C^{\ast}(S,\preceq,\mathcal{M})$ is
a~$C^*$-algebra along with a~complete order representation $\iota:S\rightarrow C^{\ast}(S,\preceq,\mathcal{M})$
satisfying the following universal property:

given $\gamma:S\rightarrow B$ a~complete order representation to a~$C^*$-algebra~$B$ there is a~unique $*$-homomorphism
$\pi_{\gamma}=\pi:C^{\ast}(S,\preceq,\mathcal{M})\rightarrow B$ such that $\pi_{\gamma}\circ\iota=\gamma$.
\end{Definition}

If the universal $C^*$-algebra $C^{\ast}(S,\preceq)$ exists then so does the $C^*$-algebra
$C^{\ast}(S,\preceq,\mathcal{M})$, and it is a~quotient of $C^{\ast}(S,\preceq)$.
Remark~\ref{comp order} shows that the universal $C^*$-algebra $C^{\ast}(S,\preceq,\mathcal{M})$ is generally a~proper
quotient of $C^{\ast}(S,\preceq)$.

For an ordered $*$-semi\-group $(S,\preceq)$ and $\beta:S\rightarrow C$ an order representation in a~$C^*$-algebra~$C$
which is a~monomorphism then there is always a~matricial order $(S,\preceq,\mathcal{M}_{\beta})$ on~$S$ so that $\beta$
becomes a~complete order representation.
To see this note that $\beta_{k}:\mathcal{M}_{k}(S)\rightarrow\mathrm{M}_{k}(C)^{\rm sa}$ must also be a~monomorphism, where
for example one could set
\begin{gather*}
\mathcal{M}_{k}(S)=\mathrm{M}_{k}(S)^{\rm sa}=\beta_{k}^{-1}\big(\mathrm{M}_{k}(C)^{\rm sa}\big)
\qquad
\text{for}
\quad
k>1
\end{gather*}
(when $k=1$, $S^{\rm sa}$ is already equal to $\beta^{-1}(C^{\rm sa})$).
Def\/ine the partial order on $\mathcal{M}_{k}(S)$ for $k>1$ to be the pull back partial order $\preceq_{\beta_{k}}$of
$\mathrm{M}_{k}(C)^{\rm sa}$.
That the maps $\iota_{\tau}:\mathcal{M}_{d}(S)\rightarrow\mathcal{M}_{k}(S)$, for $\tau \in\mathcal{P}(d,k)$, are order
maps follows from $\beta_{k}\circ\iota_{\tau}=\iota_{\tau}\circ$ $\beta_{d}$, where the second $\iota_{\tau}$ is the
order map $\mathcal{M}_{d}(C)\rightarrow\mathcal{M}_{k}(C)$.
It only remains to check that this partial order on $\mathcal{M}_{k}(S)$ satisf\/ies
\begin{gather*}
[n_{i}]_{\delta}^{\ast}\left[a_{i,j}\right][n_{i}]_{\delta}\preceq[n_{i}]_{\delta}^{\ast}\left[b_{i,j}\right][n_{j}]_{\delta}
\qquad
\text{for}
\quad
[n_{i}]\in\mathrm{M}_{1,k}(S)
\end{gather*}
whenever $\left[a_{i,j}\right] \preceq\left[b_{i,j}\right] $ in $\mathcal{M}_{k}(S)$.
By def\/inition this follows if
\begin{gather*}
\beta_{k}[(n_{i}^{\ast}a_{ij}n_{j})]\leq\beta_{k}[(n_{i}^{\ast}b_{ij} n_{j})]\qquad\text{in}\quad \mathrm{M}_{k}(C)^{\rm sa}
\end{gather*}
whenever $\left[a_{i,j}\right] \preceq\left[b_{i,j}\right] $ in $\mathcal{M}_{k}(S)$.
However
\begin{gather*}
\beta_{k}[(n_{i}^{\ast}a_{ij}n_{j})]=[\beta(n_{i})^{\ast}\beta(a_{ij})\beta(n_{j})]
\end{gather*}
is the product of matrices $N^{\ast}\beta_{k}[a_{ij}]N$ in $\mathrm{M}_{k}(C)$ where~$N$ is the $k\times k$ diagonal
matrix with entries $\beta(n_{i})$ $(1\leq i\leq k)$ and all other entries $0$.
By def\/inition of the pull back partial order $\beta_{k}[a_{ij}]\preceq\beta_{k}[b_{ij}]$ whenever $\left[a_{i,j}\right]
\preceq\left[b_{i,j}\right] $, so the required condition follows from the properties of the ordering in a~$C^*$-algebra.
If there is an order representation of $(S,\preceq)$ which is a~monomorphism then this is equivalent to the canonical
order representation $\iota:S\rightarrow C^{\ast}(S,\preceq)$ being a~monomorphism.
The above shows that there is a~matricial order $(S,\preceq,\mathcal{M}_{\iota})$ on~$S$ so that $\iota$ becomes
a~complete order representation.
Note that in the procedure to obtain a~matricial order we did not pull back the partial order on the $C^*$-algebra to
the ordered $*$-semi\-group to obtain the ordered $*$-semi\-group $(S,\preceq_{\iota})$.
However, there is no problem doing this since one can replace the partial order $\preceq$ with $\preceq_{\iota}$without
altering the universal $C^*$-algebra for $(S,\preceq)$ (Remark~\ref{orders}).

\begin{Proposition}
For an ordered $*$-semi\-group $(S,\preceq)$ assume that $C^{\ast}(S,\preceq)$ exists and that the canonical order
representation
\begin{gather*}
\iota: \ S\rightarrow C^{\ast}(S,\preceq)
\end{gather*}
is a~monomorphism.
If a~matricial ordering $(S,\preceq,\mathcal{M}_{\iota})$ is defined via $\iota$ then
\begin{gather*}
C^{\ast}(S,\preceq,\mathcal{M}_{\iota})\cong C^{\ast}(S,\preceq).
\end{gather*}
\end{Proposition}

\begin{proof}
It is enough to show that any order representation $\beta:(S,\preceq)\rightarrow C$ is a~complete order representation
of $(S,\preceq,\mathcal{M}_{\iota})$.
By def\/inition $a\preceq b$ in $\mathcal{M}_{k}(S)$ if and only if $\iota_{k}(a)\leq\iota_{k}(b)$ in
$\mathrm{M}_{k}(C^{\ast}(S,\preceq))^{\rm sa}$.
Since $(\pi_{\beta})_{k}$ is a~$*$-homomorphism from $\mathrm{M}_{k}(C^{\ast}(S,\preceq))$ to $\mathrm{M}_{k}(C)$ it is
positive, so $(\pi_{\beta})_{k}\iota_{k}(a)\leq(\pi_{\beta})_{k}\iota_{k}(b)$, i.e., $\beta_{k}(a)\leq\beta_{k}(b)$.
\end{proof}

\begin{Example}
One may apply this observation to the ordered $*$-semi\-group $(\mathbb{N},\preceq)$ since the canonical map
$\iota:(\mathbb{N},\preceq)\rightarrow C^{\ast} (\mathbb{N},\preceq)=C((0,1])$ mapping~$1$ to the positive contraction
$f(t)=t$, $t\in(0,1]$ is an injective order representation.
Therefore the matricial order structure $(\mathbb{N},\preceq,\mathcal{M}_{\iota})$ is def\/ined and
\begin{gather*}
C^{\ast}(\mathbb{N},\preceq,\mathcal{M}_{\iota})=C^{\ast}(\mathbb{N},\preceq).
\end{gather*}
Similarly for a~discrete group~$G$ viewed as a~$*$-semi\-group.
Since $C^{\ast}((G,\preceq))$ is $C^{\ast}(G)$ and~$G$ embeds in $C^{\ast}(G)$ there is a~matricial ordering
$(G,\preceq,\mathcal{M})$, where $\mathcal{M}_{k}(G)$ consists of selfadjoint $k\times k$ matrices
$\mathrm{M}_{k}(G)^{\rm sa}$ with a~partial order obtained by pulling back the partial order on $\mathrm{M}
_{k}(C^{\ast}(G))^{\rm sa}$.
Therefore the universal $C^*$-algebras
\begin{gather*}
C_{b}^{\ast}(G),
\qquad
C^{\ast}(G),
\qquad
C^{\ast}(G,\preceq)
\qquad
\text{and}
\qquad
C^{\ast}(G,\preceq,\mathcal{M})
\end{gather*}
all coincide for a~group~$G$ when~$G$ is viewed as a~$*$-semi\-group.
\end{Example}

\begin{Example}
Consider the $*$-semi\-group~$A$.
Clearly any order representation, in fact any $*$-representation, of~$A$ is also a~complete order representation of this
matricially ordered $*$-semi\-group, so the universal $C^*$-algebra $C^{\ast}(A,\preceq,\mathcal{M})$
is (isomorphic to) $C^{\ast}(A,\mathcal{\preceq})$, and therefore all of the universal $C^*$-algebras
\begin{gather*}
C_{b}^{\ast}(A),
\qquad
C^{\ast}(A),
\qquad
C^{\ast}(A,\preceq)
\qquad
\text{and}
\qquad
C^{\ast}(A,\preceq,\mathcal{M})
\end{gather*}
coincide for the $*$-semi\-group~$A$.
\end{Example}

\subsection[The $C^*$-correspondence $_{C^{\ast}(D_{1},\preceq,\mathcal{M})}\mathcal{E}_{C^{\ast}(D_{1},\preceq,\mathcal{M})}$]
{The $\boldsymbol{C^*}$-correspondence $\boldsymbol{_{C^{\ast}(D_{1},\preceq,\mathcal{M})}\mathcal{E}_{C^{\ast}(D_{1},\preceq,\mathcal{M})}}$}

For the matricial ordered semigroup $(D_{1},\preceq,\mathcal{M})$ let
\begin{gather*}
\iota_{1}: \ D_{1}\rightarrow C^{\ast}(D_{1},\preceq,\mathcal{M})
\end{gather*}
denote the canonical complete order representation.
Set
\begin{gather*}
\beta=\iota_{1}\circ\alpha: \ D_{1}\rightarrow C^{\ast}(D_{1},\preceq,\mathcal{M}).
\end{gather*}
This is a~complete order map, so necessarily an order map, satisfying the Schwarz inequality.
Lemma~\ref{bnd} and the comments preceding it show there is a~Hilbert module
$\mathcal{E}_{C^{\ast}(D_{1},\preceq,\mathcal{M})}$ with dense subspace (a quotient of)
$X=\mathbb{C}[D_{1}]\otimes_{\rm alg}C^{\ast} (D_{1},\preceq,\mathcal{M})$.
The two contractive module maps def\/ined by the idempotents $(-1,1)$ and $(1,-1)$ of $D_{1}$, along with the complete
order structure, yield a~left action of $C^{\ast}(D_{1},\preceq,\mathcal{M})$ on this Hilbert module.

\begin{Theorem}
\label{left action}
There is a~complete order representation $l:D_{1} \rightarrow\mathcal{L}
(\mathcal{E}_{C^{\ast}(D_{1},\preceq,\mathcal{M})})$ of
the matricially ordered $*$-semi\-group $D_{1}$ which yields a~$*$-representation
\begin{gather*}
\phi: \ C^{\ast}(D_{1},\preceq,\mathcal{M})\rightarrow\mathcal{L}(\mathcal{E}_{C^{\ast}(D_{1},\preceq,\mathcal{M})}).
\end{gather*}
This defines a~correspondence $\mathcal{E}$ over the $C^*$-algebra $C^{\ast}(D_{1},\preceq,\mathcal{M})$.
\end{Theorem}

\begin{proof}
Choose $a,b\in D_{1}^{\rm sa}$ with $a\preceq b$ and $[s_{i}]\in\mathrm{M}_{1,k}(D_{1})$.
Then $a^{2}\preceq a$, and since $\beta$ is a~complete order map
\begin{gather*}
\beta_{k}([s_{i}]^{\ast}a^{2}[s_{j}])\leq\beta_{k}([s_{i}]^{\ast}a[s_{j}])\leq\beta_{k}([s_{i}]^{\ast}b[s_{j}])
\end{gather*}
in the $C^*$-algebra $\mathrm{M}_{k}(C^{\ast}(D_{1},\preceq,\mathcal{M}))$.
Since either $b\preceq(-1,1)$, or $b\preceq(1,-1)$ in $(D_{1},\preceq)$ it follows that
$\beta_{k}([s_{i}]^{\ast}b[s_{j}])\leq\beta_{k}([s_{i}]^{\ast}(1,-1)[s_{j}])$ (or
$\beta_{k}([s_{i}]^{\ast}(-1,1)[s_{j}])\leq [\beta(s_{i}^{\ast}s_{j})]$).
For
\begin{gather*}
x=\sum\limits_{i=1}^{k}s_{i}\otimes_{\rm alg}c_{i}\in X=\mathbb{C}[D_{1}]\otimes_{\rm alg}C^{\ast}(D_{1},\preceq,\mathcal{M})
\end{gather*}
we obtain
\begin{gather*}
\langle ax,ax\rangle=\sum\langle c_{i},\beta(s_{i}^{\ast}a^{2}s_{j})c_{j}\rangle
=\langle\overrightarrow{c},\beta_{k}([s_{i}]^{\ast}a^{2}[s_{j}])\overrightarrow{c}\rangle
\leq\langle\overrightarrow{c},\beta_{k}([s_{i}]^{\ast}a[s_{j}])\overrightarrow{c}\rangle
\\
\phantom{\langle ax,ax\rangle}
=\langle x,ax\rangle \leq\langle x,bx\rangle \leq\langle x,x\rangle,
\end{gather*}
where $\overrightarrow{c}=(c_{1},\dots,c_{k})$.
It follows that the map $x\rightarrow ax$ on~$X$ is bounded, so extends to a~bounded map $l(a)$, of norm less than or
equal to~$1$, on $\mathcal{E}_{C^{\ast}(D_{1},\preceq,\mathcal{M})}$.
It is also adjointable, in fact selfadjoint.
The above inequalities show
\begin{gather*}
\langle x,l\left(a\right) x\rangle \leq\langle x,l(b)x\rangle
\qquad
\text{for}
\quad
x\in X/N.
\end{gather*}
Since these are bounded maps this inequality persists for $x\in\mathcal{E}_{C^{\ast}(D_{1},\preceq,\mathcal{M})}$,
which implies~\cite[Lemma~4.1]{l}
that $l(a)\leq l(b)$ in $\mathcal{L}(\mathcal{E}_{C^{\ast}(D_{1},\preceq,\mathcal{M})})$.

For arbitrary $a\in D_{1}$ the element $a^{\ast}a$ is in $D_{1}^{\rm sa}$ so
$l(a^{\ast}a)\in\mathcal{L}(\mathcal{E}_{C^{\ast}(D_{1},\preceq,\mathcal{M})})$.
We have
\begin{gather*}
0\leq\langle ax,ax\rangle \leq\langle x,a^{\ast} ax\rangle \leq\langle x,x\rangle
\end{gather*}
so the map $x\rightarrow ax$ on~$X$ is bounded, and extends to a~bounded map $l(a)$ on
$\mathcal{E}_{C^{\ast}(D_{1},\preceq,\mathcal{M})}$.
This map is clearly adjointable with adjoint $l(a^{\ast})$.
The map~$l$ is therefore an order representation of $(D_{1},\preceq)$
in the $C^*$-algebra $\mathcal{L}(\mathcal{E}_{C^{\ast}(D_{1},\preceq,\mathcal{M})})$.

To show~$l$ is a~complete order map of $D_{1}$ to $\mathcal{L}(\mathcal{E}_{C^{\ast}(D_{1},\preceq,\mathcal{M})})$ it
is enough to show
\begin{gather*}
l_{k}([w_{i}]^{\ast}[w_{j}])\leq
l_{k}([w_{i}]_{\pm}^{\ast}[w_{j}]_{\pm})\qquad\text{in}\quad \mathcal{M}_{k}(\mathcal{L}(\mathcal{E}_{C^{\ast}(D_{1},\preceq,\mathcal{M})}))
\end{gather*}
(so as elements of the $C^*$-algebra $\mathcal{L}(\oplus^{k}\mathcal{E})$) whenever
$[w_{i}]^{\ast}[w_{j}]\in\mathcal{M}_{k}(D_{1})$, with $[w_{i}]\in\mathrm{M}_{1,k}(A)$.
This follows if
\begin{gather*}
\langle \overrightarrow{x},l_{k}([w_{i}]^{\ast}[w_{j}])\overrightarrow {x}\rangle \leq\langle
\overrightarrow{x},l_{k}([w_{i}]_{\pm}^{\ast}[w_{j}]_{\pm})\overrightarrow{x}\rangle
\end{gather*}
for $\overrightarrow{x}=(x_{1},\dots,x_{k})$ with $x_{i}$ in the dense space
$X=\mathbb{C}[D_{1}]\otimes_{\rm alg}C^{\ast}(D_{1},\preceq,\mathcal{M})$ of the module
$\mathcal{E}_{C^{\ast}(D_{1},\preceq,\mathcal{M})}$  \cite[Lemma~4.1]{l}.
Writing
\begin{gather*}
x_{i}=\sum\limits_{h_{i}=1}^{t_{i}}s_{i,h_{i}}\otimes c_{i,h_{i}}
\end{gather*}
with $s_{i,h_{i}}\in D_{1}$, $c_{i,h_{i}}\in C^{\ast}(D_{1},\preceq,\mathcal{M})$, we have
\begin{gather*}
\langle \overrightarrow{x},l_{k}([w_{i}]^{\ast}[w_{j}])\overrightarrow {x}\rangle
=\sum\limits_{i=1}^{k}\sum\limits_{j=1}^{k}\sum\limits_{h_{i}=1}^{t_{i}} \sum\limits_{h_{j}=1}^{t_{j}}\langle
s_{i,h_{i}}\otimes c_{i,h_{i}},l(w_{i}^{\ast}w_{j})s_{j,h_{j}}\otimes c_{j,h_{j}}\rangle.
\end{gather*}
This is equal to the same sum over the terms
\begin{gather*}
\langle c_{i,h_{i}},\beta(s_{i,h_{i}}^{\ast}(w_{i}^{\ast}w_{j})s_{j,h_{j}})c_{j,h_{j}}\rangle
_{C^{\ast}(D_{1},\preceq,\mathcal{M})}.
\end{gather*}

Set
\begin{gather*}
r=\sum\limits_{i=1}^{k}t_{i},\tau=(t_{1},\dots,t_{d})\in\mathcal{P}(k,r),
\\
[s_{i}]=[s_{1,1},\dots,s_{1,t_{1}},s_{2,1},\dots,s_{k,t_{k}}]\in\mathrm{M}_{1,r}(D_{1}),
\end{gather*}
and{\samepage
\begin{gather*}
\overrightarrow{c}=(c_{1,1},\dots,c_{1,t_{1}},c_{2,1},\dots,c_{k,t_{k}})
\end{gather*}
in the Hilbert $C$-module $\oplus_{1}^{r}C$.}

Therefore$\langle \overrightarrow{x},l_{k}([w_{i}]^{\ast}[w_{j}])\overrightarrow{x}\rangle =\langle
\overrightarrow{c},T\overrightarrow{c}\rangle $ where~$T$ is the matrix
\begin{gather*}
[\beta(s_{i,h_{i}}^{\ast}(w_{i}^{\ast}w_{j})s_{j,h_{j}})]=\beta_{r}[(s_{i,h_{i}}^{\ast}(w_{i}^{\ast}w_{j})s_{j,h_{j}})]=\beta_{r}
([s_{i}]_{\delta}^{\ast}\left[w_{i}^{\ast}w_{j}\right]_{\tau} [s_{j}]_{\delta})
\end{gather*}
in $\mathrm{M}_{r}(C)$.
Note $\beta$ is a~complete order map and
\begin{gather*}
\iota_{\tau}:\ \mathcal{M}_{k}(D_{1})\rightarrow\mathcal{M}_{r}(D_{1})
\end{gather*}
is an order map.
Since
\begin{gather*}
\left[w_{i}^{\ast}w_{j}\right] \preceq\left[w_{i\pm}^{\ast}w_{j\pm}\right] \qquad\text{in}\quad \mathcal{M}_{k}(D_{1}),
\end{gather*}
the desired inequality follows and~$l$ is a~complete order map.

The universal property of $C^{\ast}(D_{1},\preceq,\mathcal{M})$ provides a~$*$-representation
\begin{gather*}
(\phi:=)\pi_{l}: \ C^{\ast}(D_{1},\preceq,\mathcal{M})\rightarrow\mathcal{L}(\mathcal{E}_{C^{\ast}(D_{1},\preceq,\mathcal{M})})
\end{gather*}
with $\pi_{l}\circ\iota_{1}=l$.
\end{proof}

\begin{Remark}
\label{corr to C}
For a~complete order representation $\eta:(D_{1},\preceq,\mathcal{M})\rightarrow C$ to a~$C^*$-algebra~$C$ the argument
in Theorem~\ref{left action} also applies to $\beta=\eta\circ\alpha$.
This leads to a~Hilbert module $\mathcal{E}_{C}$ (the closure of $\mathbb{C} [D_{1}]\otimes_{\rm alg}C$) over the
$C^*$-algebra~$C$, along with a~$*$-homomorphism
$\phi_{C}:C^{\ast}(D_{1},\preceq,\mathcal{M})$ $\rightarrow\mathcal{L}(\mathcal{E}_{C})$
def\/ining a~$C^*$-correspondence
$_{C^{\ast}(D_{1},\preceq,\mathcal{M})}\mathcal{E}_{C}$ from $C^{\ast}(D_{1},\preceq,\mathcal{M})$ to~$C$.
\end{Remark}

\subsection[Intrinsic description of $\mathcal{M}_{k}(D_{1})$]{Intrinsic description of $\boldsymbol{\mathcal{M}_{k}(D_{1})}$}

We describe the elements in $\mathcal{M}_{k}(D_{1})$ intrinsically, making use of Propositions~\ref{details} and~\ref{pain}.
This allows one to f\/ind the elements above a~given element in $\mathcal{M}_{k}(D_{1})$, leading to a~complete order
version of the earlier crucial extension result Proposition~\ref{extension 0}.
As in Proposition~\ref{details} the $*$-semi\-group homomorphism $\tau:A\rightarrow\mathbb{Z}$ is a~useful tool.

\begin{Proposition}
If $\mathbf{n}=[w_{i}]^{\ast}[w_{j}]\in\mathcal{M}_{k}(D_{1})$ then $\tau(w_{i}^{\ast})=\tau(w_{j}^{\ast})$ for all $i$, $j$.
The subset
\begin{gather*}
\big\{\tau(w_{i}^{\ast})\,|\,\mathbf{n}=[w_{i}]^{\ast}[w_{j}],~[w_{i}]\in\mathrm{M}_{1,k}(A)~\text{with reduced entries}\big\}
\end{gather*}
of $\mathbb{Z}$ consists of a~single integer $r\leq1$ if and only if $\mathbf{n}$ has a~unique factorization
$[w_{i}]^{\ast}[w_{j}]$.
Otherwise this set is $\{r,r-1\}$ for some $r\leq1$.
\end{Proposition}

\begin{proof}
Write $\mathbf{n}=[n_{ij}]=[w_{i}]^{\ast}[w_{j}]\in\mathcal{M}_{k}(D_{1})$, where $[w_{i}]\in\mathrm{M}_{1,k}(A)$ has
reduced entries.
First note that for all~$i$ and~$j$, $n_{ij}=w_{i}^{\ast}w_{j}$ is an element of $D_{1}$, so of $A_{0}$, and
\begin{gather*}
0=\tau(n_{ij})=\tau(w_{i}^{\ast})+\tau(w_{j})=\tau(w_{i}^{\ast})-\tau(w_{j}).
\end{gather*}
Therefore $\tau(w_{i}^{\ast})=\tau(w_{j}^{\ast})$ for all~$i$ and~$j$.
Remark~\ref{b} shows this set is bounded above by~$1$.
By Proposition~\ref{pain} there is exactly one factorization $[w_{i}]^{\ast}[w_{j}]$ of $\mathbf{n}
\in\mathcal{M}_{k}(D_{1})$ unless $(-1,1)[w_{i}]=[w_{i}]$ or $(1,-1)[w_{i}]=[w_{i}]$, in which case there is a~second
factorization $\mathbf{n} =[v_{i}]^{\ast}[v_{j}]$ with
\begin{gather*}
[v_{i}]=(1)[w_{i}]\quad\text{if}\quad (-1,1)[w_{i}]=[w_{i}],
\qquad\text{or}\qquad
[v_{i}]=(-1)[w_{i}]\quad\text{if}\quad (1,-1)[w_{i}]=[w_{i}].
\end{gather*}
Therefore the set of possible values for $\tau(w_{i}^{\ast})$, where $[w_{i}]\in\mathrm{M}_{1,k}(A)$ has reduced entries
and forms a~factorization $\mathbf{n}=[w_{i}]^{\ast}[w_{j}]$, consists of exactly one or exactly two elements.
The latter occurs if and only if $(-1,1)[w_{i}]=[w_{i}]$ or $(1,-1)[w_{i}]=[w_{i}]$, equivalently whenever the
conditions for f\/inding potential elements of $\mathcal{M}_{k}(D_{1})$ lying above $\mathbf{n}$ in the partial order are
satisf\/ied.
If $(-1,1)[w_{i}]=[w_{i}]$ and $\tau (w_{i}^{\ast})=r$ for some, and therefore all,~$i$ then
$\tau(v_{i}^{\ast})=\tau(w_{i}^{\ast}(-1))=r-1$ for all~$i$, where $[v_{i}]=(1)[w_{i}]$.
\end{proof}

Under the hypothesis of the preceding proposition we further examine the possible values of $\tau(w_{i}^{\ast})$,
resulting in Proposition~\ref{matrix structure}.
The result is crucial for the following extension result, Proposition~\ref{extension}.

\begin{Proposition}
\label{matrix structure}
If $\mathbf{n}=[n_{ij}]=[w_{i}]^{\ast}[w_{j}]\in\mathcal{M}_{k}(D_{1})$ then one of the following holds.
\begin{enumerate}\itemsep=0pt
\item[$1.$] If $\tau(w_{i}^{\ast})=1$ for some $i$
then $[w_{i}]=(-1,1)[w_{i}]$ and there is $[m_{i}]\in\mathrm{M}_{1,k}(D_{1})$, where for each~$i$ either
$(-1,1)m_{i}=m_{i}$ or $m_{i}=(1,-1)$, so that
\begin{gather*}
\mathbf{n}=[m_{i}]^{\ast}(1,-1)[m_{i}].
\end{gather*}
The set
\begin{gather*}
\{\tau(w_{i}^{\ast})\,|\,\mathbf{n}=[w_{i}]^{\ast}[w_{j}],[w_{i}]\in\mathrm{M}_{1,k}(A)\ \text{with reduced entries}\}=\{1,0\}
\end{gather*}

\item[$2.$]
If $\tau(w_{i}^{\ast})=0$  for some $i$
then $[w_{i}]=(-1,1)[w_{i}]$ or $[w_{i}]=(1,-1)[w_{i}]$.
If $[w_{i}]=(1,-1)[w_{i}]$ then $[v_{i}]=(-1)[w_{i}]$ yields a~factorization $[v_{i}]^{\ast}[v_{j}]$ with
$\tau(v_{i}^{\ast})=1$ which is case~$1$.
Otherwise $[w_{i}]=(-1,1)[w_{i}]$ and $[v_{i}]=(1)[w_{i}]$ yields a~factorization $[v_{i}]^{\ast}[v_{j}]$ with
$\tau(v_{i}^{\ast})=-1$.
In this case there is $[a_{i}]\in\mathrm{M}_{1,k}(D_{0})$ with $[a_{i}]^{\ast} [a_{j}]\in\mathcal{M}_{k}(D_{0})$, and
$[m_{i}]\in\mathrm{M}_{1,k}(D_{1})$ so that
\begin{gather*}
\mathbf{n}=[m_{i}]_{\delta}^{\ast}[a_{i}]^{\ast}[a_{j}][m_{j}]_{\delta}
\end{gather*}
$(\delta=(1,\dots,1)\in\mathcal{P}(k,k))$.
The element $[m_{i}]$ is unnecessary if $\mathbf{n}\in\mathcal{M}_{k}(D_{0})$.

\item[$3.$] If $\tau(w_{i}^{\ast})\leq-1$ for all $i$
for any factorization $[w_{i}]^{\ast}[w_{j}]$ of $\mathbf{n.}$ If neither $[w_{i}]=(-1,1)[w_{i}]$ nor
$[w_{i}]=(1,-1)[w_{i}]$ is satisfied, $\mathbf{n}$ is maximal in the partial order on $\mathcal{M}_{k}(D_{1})$.
If one of these is satisfied then there is $[m_{i}]\in\mathrm{M}_{1,k}(D_{1})$ and $[\lambda_{i}]\in\mathrm{M}_{1,k}(A)$
with $[\lambda_{i}]^{\ast}[\lambda_{j}]\in\mathcal{M}_{k}(D_{0})$, $\lambda_{i}{}^{\ast}\lambda_{j}
\in\operatorname*{Irr}(D_{0})\backslash\{(-1,1)\}$ so that
\begin{gather*}
\mathbf{n}=[m_{i}]_{\delta}^{\ast}[\lambda_{i}]^{\ast}[\lambda_{j}][m_{j}]_{\delta}.
\end{gather*}
\end{enumerate}
\end{Proposition}

\begin{proof}
If
$
\tau(w_{i}^{\ast})=1$ for some $i$
the last proposition shows $\tau(w_{i}^{\ast})=1$ for all~$i$.
Therefore $n_{ii}=w_{i}^{\ast}w_{i}$ is of the form $m_{i}^{\ast}(1,-1)m_{i}$, where $m_{i}\in D_{1}$ and
$(-1,1)m_{i}=m_{i}$ if it occurs (Lemma~\ref{sa lemma} and Proposition~\ref{details}).
It follows $(-1)m_{i}$, when reduced in~$A$, must equal $w_{i}$ and so $[w_{i}]=(-1,1)[w_{i}]$.
If $m_{j}$ does not occur for some~$j$ then set $m_{j}=(1,-1)$.
Thus $\mathbf{n}=[m_{i}]^{\ast}(1,-1)[m_{i}]$ with $[m_{i}]\in\mathrm{M}_{1,k}(D_{1})$, where, for each~$j$, either
$(-1,1)m_{j}=m_{j}$ or $m_{j}=(1,-1)$.
The term $(1,-1)$ therefore occurs as the central term of $n_{ii}$ for all $i\mathbf{.}$ Since $w_{i}=(-1)m_{i}$ for
all~$i$, Proposition~\ref{pain} implies that there is an alternative
factorization $\mathbf{n}=[v_{i}]^{\ast}[v_{j}]$
with $[v_{i}]=(1)[w_{i}]$, so then $\tau(v_{i}^{\ast})=\tau(m_{i})=0$.
Thus if $\tau(w_{i}^{\ast})=1$ for some~$i$ then $n_{ii}=w_{i}^{\ast}w_{i}$ must be the minimal (length) factorization.
In conclusion the term $(1,-1)$ occurs as the central term of $n_{ii}$, viewed in $D_{1}$, for some~$i$ if and only if
this occurs for all~$i$.
This in turn is equivalent to
\begin{gather*}
\{\tau(w_{i}^{\ast})\,|\,\mathbf{n}=[w_{i}]^{\ast}[w_{j}],~[w_{i}]\in\mathrm{M}_{1,k}(A)~\text{with reduced entries}\}=\{1,0\}.
\end{gather*}

Next assume that $(1,-1)$ does not occur as the middle term of any $n_{ii}$, so
\begin{gather*}
\tau(w_{i}^{\ast})\leq0\quad\text{for all}\quad i
\end{gather*}
for any factorization $\mathbf{n}=[w_{i}]^{\ast}[w_{j}]\in\mathcal{M}_{k}(D_{1})$, $[w_{i}]\in\mathrm{M}_{1,k}(A)$ with
reduced entries in~$A$.

Consider the case that $\mathbf{n}=[w_{i}]^{\ast}[w_{j}]$ with $\tau (w_{i}^{\ast})=0$ for some~$i$, therefore for
all~$i$.
Then $w_{i}\in D_{1}$ for all~$i$, so by Proposition~\ref{details} each $w_{i}$ must have the form $a_{i}m_{i}$ with
$a_{i}\in D_{0}$ and $m_{i}\in D_{1}$.
If $w_{i}\in D_{0}$ for all~$i$ then $\mathbf{n}=[a_{i}]^{\ast}[a_{j}]\in\mathcal{M}_{k}(D_{0})$.
Otherwise, if $m_{i}$ does not occur for some~$i$, set $m_{i}=(-1,1)$.
It follows that $[w_{i}]=(-1,1)[w_{i}]$ and $\mathbf{n}$ has two factorizations.
Therefore this case is equivalent to
\begin{gather*}
\{\tau(w_{i}^{\ast})\,|\,\mathbf{n}=[w_{i}]^{\ast}[w_{j}],~[w_{i}]\in\mathrm{M}_{1,k}(A)~\text{with reduced entries}\}=\{0,-1\}.
\end{gather*}
In this situation it is possible that $a_{i}=(-1,1)$.
For example consider
\begin{gather*}
[w_{i}]=((-2,3,-2,1),(-1,1))\in\mathrm{M}_{1,2}(A),
\end{gather*}
where $a_{1}=(-2,2)$, $m_{1}=(1,-1)(-1,1)$ and $m_{2}=a_{2}=(-1,1)$.

The remaining case is
\begin{gather*}
\tau(w_{i}^{\ast})\leq-1\quad\text{for any, and therefore all,}~i
\end{gather*}
and for any factorization $\mathbf{n}=[w_{i}]^{\ast}[w_{j}]\in\mathcal{M}_{k}(D_{1})$, $[w_{i}]\in\mathrm{M}_{1,k}(A)$
with reduced entries in~$A$.
Assume either $[w_{i}]=(-1,1)[w_{i}]$ or $[w_{i}]=(1,-1)[w_{i}]$ is satisf\/ied.
Then $\mathbf{n}=[v_{i}]^{\ast}[v_{j}]$ is another factorization, where in the f\/irst case $[v_{i}]=(1)[w_{i}]$, so
$\tau(v_{i}^{\ast})=\tau(w_{i}^{\ast})-1$, or in the second case $[v_{i}]=(-1)[w_{i}]$, so $\tau(v_{i}^{\ast})=\tau
(w_{i}^{\ast})+1$.
By hypothesis $\tau(v_{i}^{\ast})\leq-1$ for all~$i$ also, so $\tau(v_{i}^{\ast})\leq-2$ in the f\/irst case, while in the
second case $\tau(w_{i}^{\ast})\leq-2$.
Thus, if $\tau(w_{i}^{\ast})=-1$ for some (all)~$i$ then $[w_{i}]=(1,-1)[w_{i}]$ is not possible.
It follows $w_{i}$ cannot have the form $(-1,1)m_{i}$ (since then $\tau(w_{i}^{\ast})=0$) or $(1)m_{i}$ for some
$m_{i}\in D_{1}$ (since then $[w_{i}]=(1,-1)[w_{i}]$ and $\tau (w_{i}^{\ast})=-1$).
In particular the central term of $n_{ii}$ cannot be $(-1,1)$ for any~$i$.
It follows that the central term of $n_{ii}$ must be irreducible in $D_{0}$, as otherwise there would be a~factorization
$[w_{i}]^{\ast}[w_{j}]$ with $\tau(w_{i}^{\ast})=0$ for some (so all)~$i$.
Lemma~\ref{fund lemma} then implies that the element $[\lambda_{i}]^{\ast}[\lambda_{j}]\in\mathcal{M}_{k}(D_{0})$, where
$[\lambda_{i}]\in \mathrm{M}_{1,k}(A)$ is used to capture these irreducible middle elements of $n_{ii}$, must have all
its entries in $\operatorname*{Irr}(D_{0})\backslash\{(-1,1)\}$.
\end{proof}

Although these forms are to be viewed as formal matrix structures for the matrices
$\mathbf{n}=[w_{i}]^{\ast}[w_{j}]\in\mathcal{M}_{k}(D_{1})$, they are realizable decompositions whenever the semigroups
are represented in a~$C^*$-algebra.

Given an element $\mathbf{n}$ in $\mathcal{M}_{k}(D_{1})$ we now proceed to use Proposition~\ref{matrix structure} to
make clear which elements $\mathbf{r}$ occur in $\mathcal{M}_{k}(D_{1})$ with $\mathbf{n}\preceq \mathbf{r}$.
This will prove useful in establishing a~version of the previous extension result (Proposition~\ref{extension 0}).

For example, suppose
\begin{gather*}
\mathbf{n}=[m_{i}]^{\ast}(1,-1)[m_{i}]
\end{gather*}
as in case~1 of Proposition~\ref{matrix structure}.
Then the choices for $[w_{i}]_{+}$ are either $m_{i}$ or $(1,-1)m_{i}$ (when viewed as an element of $D_{1}$).
Def\/ine $[t_{i}]\in\mathrm{M}_{1,k}(D_{1})$ to be the element with $t_{i}$ equal to $m_{i}$ or $(1,-1)m_{i}$.
Then $[t_{i}]^{\ast}(1,-1)[t_{j}]=\mathbf{n}$, so $\mathbf{n}\preceq[w_{i}]_{+}^{\ast}[w_{j}]_{+}=[t_{i}]^{\ast}[t_{j}]$.
The other alternative
factorization for such an element $\mathbf{n}$ is $[v_{i}]^{\ast}[v_{j}]$ with
$[v_{i}]=(-1)[w_{i}]$, and $(-1,1)[v_{i}]=[v_{i}]$ in this case.
However $[v_{i}]_{-}^{\ast} [v_{j}]_{-}$ is just $[v_{i}]^{\ast}[v_{j}]$ for any choice of $[v_{i}]_{-}$, so no elements
larger than $\mathbf{n}$ are found through this latter choice of factorization for $\mathbf{n}$.

Suppose (cf.\ case~2 and~3
of Proposition~\ref{matrix structure})
\begin{gather*}
\mathbf{n}=[m_{i}]_{\delta}^{\ast}\mathbf{c}[m_{i}]_{\delta}
\end{gather*}
with $\mathbf{c}\in\mathcal{M}_{k}(D_{0})$.
If $\mathbf{c}$ is not maximal in $\mathcal{M}_{k}(D_{0})$ then $\mathbf{c\prec d}$ in $\mathcal{M}_{k}(D_{0})$
and $\mathbf{n\prec}[m_{i}]_{\delta}^{\ast}\mathbf{d}[m_{i}]_{\delta}$
in $\mathcal{M}_{k}(D_{1})$ by Corollary~\ref{mo}.

If $\mathbf{n}$ is as in case~2
of Proposition~\ref{matrix structure} then
\begin{gather*}
\mathbf{c}=[a_{i}]^{\ast}[a_{j}]\in\mathcal{M}_{k}(D_{0}).
\end{gather*}
If $\mathbf{c}$ is not maximal in $\mathcal{M}_{k}(D_{0})$ then the last paragraph holds, so assume $\mathbf{c}$ is
maximal in $\mathcal{M}_{k} (D_{0})$.
Therefore $\mathbf{c=}[a_{i}]_{+}^{\ast}[a_{j}]_{+}$ for all choices of $[a_{i}]_{+}$ and so for each~$i$, $a_{i}$ must
have the form $\left(-1,r,\dots \right) \in D_{0}$.
By Proposition~\ref{D 1}~$r$ must be~$1$ and $a_{i}=(-1,1)$ for all~$i$.
If $\mathbf{n}=\mathbf{c}$ then $\mathbf{n}$ is maximal.
If $\mathbf{n}\neq\mathbf{c}$ then, using the same approach for $\mathbf{n}$ satisfying case~1,
def\/ine $[t_{i}]\in\mathrm{M}_{1,k}(D_{1})$, where $t_{i}$ is chosen to be $m_{i}$ or $(-1,1)m_{i}$.
Then $n=[t_{i}]^{\ast}(-1,1)[t_{j}]\preceq[t_{i}]^{\ast}[t_{j}]$ yields all possible elements in
$\mathcal{M}_{k}(D_{1})$ (including the extreme choice $[m_{i}]^{\ast}[m_{i}]$) above $\mathbf{n}$.

If $\mathbf{n}$ is as in case~3
and has two factorizations, then
\begin{gather*}
\mathbf{n}=[m_{i}]_{\delta}^{\ast}\mathbf{c}[m_{i}]_{\delta}
\end{gather*}
with $\mathbf{c}\in\mathcal{M}_{k}(D_{0})$ with entries in $\operatorname*{Irr}(D_{0})\backslash\{(-1,1)\}$.
Given any element~$n$ in $D_{0}^{\rm sa}\backslash\{(-1,1)\}$ there is an element~$m$ with $n\prec m$ in $D_{0}^{\rm sa}$,
therefore there is an element $\mathbf{d}$ in $\mathcal{M}_{k}(D_{0})$ with $\mathbf{c}\prec\mathbf{d}$.

It is now possible to establish the version of the previous extension result for $C^{\ast}(D_{j},\preceq)$
(Proposition~\ref{extension 0}) for the universal complete order $C^*$-algebras $C^{\ast}(D_{j},\preceq,\mathcal{M}
_{k})$, $j=0,1$.
The main aspect of the extension result that needs to be addressed is obtaining a~complete order map of
$(D_{1},\preceq,\mathcal{M})$ extending a~given complete order map of $(D_{0},\preceq,\mathcal{M})$.
For this it is suf\/f\/icient to consider the situation where an element $\mathbf{n}\in\mathcal{M}_{k}(D_{1})$ is not
maximal in the partial order, consequently we may restrict attention to those $\mathbf{n}=[w_{i}]^{\ast}[w_{j}]$, where
either $(-1,1)[w_{i}]=[w_{i}]$ or $(1,-1)[w_{i}]=[w_{i}]$ is satisf\/ied.
For the f\/irst two parts of the last proposition this is automatically the case.

\begin{Proposition}
\label{extension}
Given a~complete order representation $\sigma:D_{0} \rightarrow C$ in a~$C^*$-algebra~$C$ there is a~complete order
representation $\rho:D_{1}\rightarrow C$ that extends $\sigma$.
\end{Proposition}

\begin{proof}
Since $\sigma$ is a~complete order representation it is an order representation of $(D_{0},\preceq)$, so we may assume
there exists an order representation $\rho:D_{1}\rightarrow C$ extending $\sigma$ with $\rho(1,-1)=q$ a~nonzero
projection in~$C$ (Proposition~\ref{extension 0}).
It is suf\/f\/icient to show that $\rho$ is a~complete order map, so $\rho_{k}
([w_{i}]^{\ast}[w_{j}])\leq\rho_{k}([w_{i}]_{\pm}^{\ast}[w_{j}]_{\pm})$ whenever
$\mathbf{n}=[w_{i}]^{\ast}[w_{j}]\in\mathcal{M}_{k}(D_{1})$ and $[w_{i}]^{\ast}[w_{j}]\preceq[w_{i}]_{\pm}^{\ast}[w_{j}]_{\pm}$
is a~basic step def\/ining the partial order in $\mathcal{M}_{k}(D_{1})$.
Assume that $\mathbf{n}$ is not maximal with respect to the partial order, since otherwise there is nothing to show.
Consider the possibi\-li\-ties for $\mathbf{n}$ and the possible elements $[w_{i}]_{\pm}^{\ast}[w_{j}]_{\pm}$ larger than
$\mathbf{n}$ described after Proposition~\ref{matrix structure}.
An illustration of two situations suf\/f\/ices.

For example in the f\/irst case
\begin{gather*}
\mathbf{n=}[t_{i}]^{\ast}(1,-1)[t_{j}]\preceq[t_{i}]^{\ast}[t_{j}],
\end{gather*}
where $[t_{i}]\in\mathrm{M}_{1,k}(D_{1})$.
Then
\begin{gather*}
\rho_{k}(\mathbf{n})=[\rho(t_{i})]^{\ast}q[\rho(t_{j})]
\leq[\rho(t_{i})]^{\ast}[\rho(t_{j})]=\rho_{k}([t_{i}]^{\ast}[t_{j}])
\end{gather*}
in $\mathcal{M}_{k}(C)$.

For another situation suppose
\begin{gather*}
\mathbf{n}=[m_{i}]_{\delta}^{\ast}\mathbf{c}[m_{i}]_{\delta}
\qquad\text{with}\quad \mathbf{c\prec d}\quad\text{in}\quad\mathcal{M}_{k}(D_{0}).
\end{gather*}
Then the usual order properties in the $C^*$-algebra $\mathrm{M}_{k}(C)$ yield
\begin{gather*}
\rho_{k}(\mathbf{n})=\rho_{k}([m_{i}]_{\delta}^{\ast})\sigma_{k}
(\mathbf{c})\rho_{k}([m_{i}]_{\delta})\leq\rho_{k}([m_{i}]_{\delta}^{\ast})\sigma_{k}(\mathbf{d})\rho_{k}([m_{i}]_{\delta})
\end{gather*}
since by hypothesis $\sigma_{k}$ is an order map from $\mathcal{M}_{k}(D_{0})$ to $\mathcal{M}_{k}(C)$.
\end{proof}

The universal property implies that there is a~$*$-homomorphism
\begin{gather*}
\pi_{\iota}:\ C^{\ast}(D_{0},\preceq,\mathcal{M})\rightarrow C^{\ast} (D_{1},\preceq,\mathcal{M}),
\end{gather*}
where $\iota:(D_{0},\preceq,\mathcal{M})\rightarrow(D_{1},\preceq,\mathcal{M})$ is the natural inclusion and
\begin{gather*}
\pi_{\iota}\circ\iota_{0}=\iota_{1}\circ\iota
\end{gather*}
with $\iota_{j}:D_{j}\rightarrow C^{\ast}(D_{j},\preceq,\mathcal{M})$ the canonical complete order representations
$(j=0,1)$.
The next result follows using the same argument in Corollary~\ref{injection 0}.

\begin{Corollary}
\label{injection}
The $*$-homomorphism
\begin{gather*}
\pi_{\iota}: \ C^{\ast}(D_{0},\preceq,\mathcal{M})\rightarrow C^{\ast} (D_{1},\preceq,\mathcal{M})
\end{gather*}
is an injection of $C^*$-algebras.
\end{Corollary}

\begin{Remark}
\label{restriction}
Since the complete order map $\alpha$ actually maps $D_{1}$ to the sub-semigroup $D_{0}$, apply Remark~\ref{corr to C}
to the complete order map
\begin{gather*}
\beta_{0}=\iota_{0}\circ\alpha: \ D_{1}\rightarrow C^{\ast}(D_{0},\preceq,\mathcal{M})
\end{gather*}
to def\/ine the Hilbert $C^{\ast}(D_{0},\preceq,\mathcal{M})$ module
$\mathcal{E}_{0}=\mathcal{E}_{C^{\ast}(D_{0},\preceq,\mathcal{M})}$ via the vector space
$X_{0}=\mathbb{C}[D_{1}]\otimes_{\rm alg}C^{\ast}(D_{0},\preceq,\mathcal{M})$.
We obtain a~$C^*$-correspondence
\begin{gather*}
_{C^{\ast}(D_{1},\preceq,\mathcal{M})}\mathcal{E}_{C^{\ast}(D_{0},\preceq,\mathcal{M})}
\end{gather*}
with left action
\begin{gather*}
\phi_{0}: \ C^{\ast}(D_{1},\preceq,\mathcal{M})\rightarrow\mathcal{L}(\mathcal{E}_{C^{\ast}(D_{0},\preceq,\mathcal{M})}).
\end{gather*}
The map $\xi:X_{0}\rightarrow X$ specif\/ied on simple tensors by $s\otimes c\rightarrow s\otimes\pi_{\iota}(c)$
for $s\in D_{1}$, $c\in C^{\ast} (D_{0}$, $\preceq,\mathcal{M})$
yields a~map of Hilbert modules
\begin{gather*}
\xi:\ \mathcal{E}_{C^{\ast}(D_{0},\preceq,\mathcal{M})}\rightarrow \mathcal{E}_{C^{\ast}(D_{1},\preceq,\mathcal{M})}.
\end{gather*}
We have
\begin{gather*}
\langle \xi x,\xi y\rangle =\pi_{\iota}(\langle x,y\rangle)
\qquad
\text{and}
\qquad
\phi(a)\circ\xi=\xi\circ\phi_{0}(a)
\end{gather*}
for $x,y\in\mathcal{E}_{C^{\ast}(D_{0},\preceq,\mathcal{M})}$, $a\in C^{\ast}(D_{1},\preceq,\mathcal{M})$.
Thus $\xi$ is isometric by the previous corollary, and is a~map of $C^*$-correspondences.
In particular
\begin{gather*}
\phi(a)=0
\qquad
\text{implies that}
\quad
\phi_{0}(a)=0.
\end{gather*}
This last observation is crucial in establishing Corollary~\ref{need} below, from which a~specif\/ic isomorphic
representation of $C^{\ast}(D_{1},\preceq,\mathcal{M})$ follows (Corollary~\ref{iso image}).
\end{Remark}

\section[A~Cuntz--Pimsner $C^*$-algebra]{A~Cuntz--Pimsner $\boldsymbol{C^*}$-algebra}
\label{section5}

We consider, for a~particular ideal~$K$ of $C^{\ast}(D_{1},\preceq,\mathcal{M})$, a~relative Cuntz--Pimsner
$C^*$-algebra $\mathcal{O} (K,\mathcal{E})$ associated with the $C^*$-correspondence
\begin{gather*}
_{C^{\ast}(D_{1},\preceq,\mathcal{M})}\mathcal{E}_{C^{\ast}(D_{1},\preceq,\mathcal{M})},
\end{gather*}
and establish an isomorphism of $\mathcal{O}(K,\mathcal{E})$ with the universal $C^*$-algebra $\mathcal{P}$ generated by a~partial isometry.

We brief\/ly sketch some relevant background and refer the reader, for example,
to~\cite{fmr,k,ms, pm} and references therein.
Let $_{B}\mathcal{E}_{B}$ be a~$C^*$-correspondence over a~$C^*$-algebra~$B$.
A~representation $(T,\pi):\mathcal{E}\rightarrow C$ of $_{B}\mathcal{E}_{B}$ in a~$C^*$-algebra~$C$ is
a~$*$-homomorphism $\pi:B\rightarrow C$ along with a~linear map $T:\mathcal{E}\rightarrow C$ which is a~$B$-bimodule
map, and which is a~Hilbert module isometry.
Therefore
\begin{gather*}
T(\phi(b)x) =\pi(b)T(x),T(xb)=T\left(x\right) \pi(b),\qquad\text{and}\qquad
\\
T(x)^{\ast}T(y)=\langle T(x),T(y)\rangle_{C}=\pi(\langle x,y\rangle_{B})
\end{gather*}
for $b\in B$, $x,y\in\mathcal{E}$.
If $(T,\pi):\mathcal{E}\rightarrow C$ is a~representation of $\mathcal{E}$ in a~$C^*$-algebra~$C$ the $C^*$-subalgebra
of~$C$ generated by $T(\mathcal{E})\cup\pi(B)$ is denoted $C^{\ast}(T,\pi)$.

For~$C$ a~$C^*$-algebra, using the identif\/ication of $\mathcal{K}(C)$ with~$C$, a~representation $(T,\pi)$ of
$\mathcal{E}$ in~$C$ yields a~$*$-homomorphism $\Psi_{T}:\mathcal{K}(\mathcal{E})\rightarrow C$ determined by $\theta
_{x,y}\rightarrow T(x)T(y)^{\ast}$.
Denote the ideal $\phi^{-1}(\mathcal{K} (\mathcal{E}))$ of~$B$ by $J(\mathcal{E})$, and say that a~representation
$(T,\pi):\mathcal{E}\rightarrow C$ is coisometric on an ideal~$K$ contained in $J(\mathcal{E})$ if
$\Psi_{T}(\phi(b))=\pi(b)$ for all $b\in K$.
There is a~universal coisometric $C^*$-algebra $\mathcal{O}(K,\mathcal{E})$ called the relative Cuntz--Pimsner algebra
of $\mathcal{E}$ determined by~$K$~\cite{fmr, ms}.
Namely, there is a~representation $(T_{\mathcal{E}},\pi_{\mathcal{E}})$
of $\mathcal{E}$ in $\mathcal{O}(K,\mathcal{E})$ coisometric on~$K$~\cite{fmr}
such that if $(T,\pi):\mathcal{E}\rightarrow C$ is a~representation coisometric on~$K$ then there is a~unique $*$-homomorphism
\begin{gather*}
\rho: \ \mathcal{O}(K,\mathcal{E})\rightarrow C
\qquad
\text{with}
\quad
(T,\pi)=\rho \circ(T_{\mathcal{E}},\pi_{\mathcal{E}}).
\end{gather*}
For the ideal $J_{\mathcal{E}}=\phi^{-1}(\mathcal{K}(\mathcal{E}))\cap (\ker\phi)^{\perp}$ of~$B$ the universal
Cuntz--Pimsner $C^*$-algebra $\mathcal{O}(J_{\mathcal{E}},\mathcal{E})$ determined by this ideal is denoted
$\mathcal{O}_{\mathcal{E}}$~\cite{k}.

Recall the correspondence
$\mathcal{E}=_{C^{\ast}(D_{1},\preceq,\mathcal{M})}\mathcal{E}_{C^{\ast}(D_{1},\preceq,\mathcal{M})}$ of the previous
section with dense subspace (a quotient of) $X=\mathbb{C}[D_{1}]\otimes_{\rm alg}C^{\ast}(D_{1},\preceq,\mathcal{M})$ and $*$-homomorphism
\begin{gather*}
\phi: \ C^{\ast}(D_{1},\preceq,\mathcal{M})\rightarrow\mathcal{L}(\mathcal{E}_{C^{\ast}(D_{1},\preceq,\mathcal{M})})
\end{gather*}
(Lemma~\ref{bnd}, Theorem~\ref{left action}).
Although the semigroup $D_{1}$ is not unital there is still a~crucial distinguished element in the space $\mathcal{E}$,
namely the (class of the) element
\begin{gather*}
e=(1,-1)\otimes\iota_{1}(-1,1)\qquad\text{in}\quad\mathcal{E}.
\end{gather*}
Since $\alpha(1,-1)=(-1,1)$ the $C^{\ast}(D_{1},\preceq,\mathcal{M})$-valued inner product
\begin{gather*}
\langle e,e\rangle =\langle \iota_{1}(-1,1),\iota_{1}
\alpha((1,-1)^{\ast}(1,-1))\iota_{1}(-1,1)\rangle =\iota_{1}(-1,1).
\end{gather*}
Therefore for $(T,\pi):_{C^{\ast}(D_{1},\preceq,\mathcal{M})}\mathcal{E}
_{C^{\ast}(D_{1},\preceq,\mathcal{M})}\rightarrow C$ a~representation we have
\begin{gather*}
T(e)^{\ast}T(e)=\pi(\langle e,e\rangle_{C^{\ast}(D_{1},\preceq,\mathcal{M})})=\pi(\iota_{1}(-1,1)).
\end{gather*}
Since $(-1,1)$ is a~selfadjoint idempotent in the $*$-semi\-group $D_{1}$, $T(e)^{\ast}T(e)$ is a~projection in~$C$ and
$T(e)$ is necessarily a~partial isometry in~$C$.
Compute, for $b\in C^{\ast}(D_{1},\preceq,\mathcal{M})$, that
\begin{gather*}
T(e)^{\ast}\pi(b)T(e)=T(e)^{\ast}T(\phi(b)e)
=\pi(\langle e,b(1,-1)\otimes(-1,1)\rangle_{C^{\ast}(D_{1},\preceq,\mathcal{M})})
\\
\phantom{T(e)^{\ast}\pi(b)T(e)}
=\pi(\sigma((-1,1)\alpha(b)(-1,1)))=\pi(\iota_{1}\alpha(b)).
\end{gather*}
Therefore the partial isometry $T(e)$ implements the complete order map $\beta=\iota_{1}\alpha$ of the ordered
$*$-semi\-group $(D_{1},\preceq)$ in the image $C^*$-algebra $\pi(C^{\ast}(D_{1},\preceq,\mathcal{M}))$.

In particular the $C^*$-subalgebra of~$C$ generated by the partial isometry $T(e)$ must contain the image under
$\pi\circ\iota_{1}$ of the $*$-subsemigroup of $D_{1}$ which is generated by the element $(-1,1)$ and the map $\alpha$.
By Theorem~\ref{gen} this is the $*$-subsemigroup $D_{0}$ of $D_{1}$ and we conclude that $\pi(\iota_{1}(D_{0}))$ is
contained in the $C^*$-algebra generated by $T(e)$.
However, since
\begin{gather*}
\pi(\iota_{1}(1,-1))T(e)=T(\phi(\iota_{1}(1,-1))e)=T(e)
\end{gather*}
the f\/inal projection $T(e)T(e)^{\ast}$ of the partial isometry $T(e)$ is only required to be less than or equal to the
projection $\pi(\iota_{1}(1,-1))$.
Thus the projection $\pi(\iota_{1}(1,-1))$ in $\pi(C^{\ast}(D_{1},\preceq,\mathcal{M}))$ may not be contained in the
$C^*$-algebra generated by $T(e)$.

\begin{Lemma}
\label{equiv}
For $X=\mathbb{C}[D_{1}]\otimes_{\rm alg}C^{\ast}(D_{1},\preceq,\mathcal{M})$ and $k,n,m\in D_{1}$ we have
\begin{gather*}
k\otimes\iota_{1}(m)=k(1,-1)\otimes\iota_{1}((-1,1)m),
\\
k\otimes\iota_{1}(\alpha(n)m)=k(1,-1)n\otimes\iota_{1}(m)
\end{gather*}
in $X/N$.
In particular
\begin{gather*}
k\otimes\iota_{1}(m)=k(1,-1)\otimes\iota_{1}(m)=k\otimes\iota_{1}((-1,1)m).
\end{gather*}
\end{Lemma}

\begin{proof}
Recall Remark~\ref{b}.
By def\/inition the $C^{\ast}(D_{1},\preceq,\mathcal{M})$-valued inner product on~$X$,
\begin{gather*}
\langle k\otimes\iota_{1}(m),a\otimes\iota_{1}(b)\rangle =\langle
\iota_{1}(m),\iota_{1}\alpha(k^{\ast}a)\iota_{1} (b)\rangle
\end{gather*}
while
\begin{gather*}
\langle k(1,-1)\otimes\iota_{1}((-1,1)m),a\otimes\iota_{1} (b)\rangle
\\
\qquad
=\langle \iota_{1}((-1,1)m),\iota_{1}\alpha((1,-1)k^{\ast}a)\iota_{1}(b)\rangle=\iota_{1}(m^{\ast}(-1,1))\iota_{1}\alpha(k^{\ast}a)\iota_{1}(b)
\\
\qquad
=\iota_{1}(m^{\ast})\iota_{1}((-1,1)\alpha(k^{\ast}a))\iota_{1}(b)
=\iota_{1}(m^{\ast})\iota_{1}(\alpha(k^{\ast}a))\iota_{1}(b),
\end{gather*}
which is the f\/irst inner product.
Since the elements $a\otimes\iota_{1}(b)$ span~$X$ the f\/irst equality follows.

For the second equality;
\begin{gather*}
\langle k\otimes\iota_{1}(\alpha(n)m),a\otimes\iota_{1}(b)\rangle =\langle
\iota_{1}(\alpha(n)m),(\iota_{1}\alpha(k^{\ast}a))\iota_{1}(b)\rangle
\\
\qquad
=\langle k(1,-1)n\otimes\iota_{1}(m),a\otimes\iota_{1}(b)\rangle =\langle \iota_{1}(m),
\iota_{1}\alpha(n^{\ast}(1,-1)k^{\ast}a)\iota_{1}(b)\rangle
\\
\qquad
=\langle \iota_{1}(m),\iota_{1}(\alpha(n)^{\ast}\alpha(k^{\ast})\alpha(k^{\ast}a))\iota_{1}(b)\rangle
=\langle \iota_{1} (\alpha(n)m),\iota_{1}(\alpha(k^{\ast}a))\iota_{1}(b)\rangle.
\end{gather*}

The last statements follow by setting $n=(1,-1)$.
\end{proof}

For $(T,\pi):\mathcal{E}\rightarrow C$ a~representation of the $C^*$-correspondence $\mathcal{E}$ over
$C^{\ast}(D_{1},\preceq,\mathcal{M})$, and $k,m\in D_{1}$, the last lemma implies
\begin{gather*}
T(k\otimes\iota_{1}(m))=T(k(1,-1)\otimes\iota_{1}((-1,1)m))
\\
\qquad
=T(l(k)e\iota_{1}(m))=T(\phi(\iota_{1}(k))e\iota_{1}(m))
=\pi(\iota_{1}(k))T(e)\pi(\iota_{1}(m)).
\end{gather*}
Therefore the image space $T(\mathcal{E})$ in~$C$ is contained in
\begin{gather*}
\pi(C^{\ast}(D_{1},\preceq,\mathcal{M}))T(e)\pi(C^{\ast}(D_{1},\preceq,\mathcal{M})).
\end{gather*}

\begin{Lemma}
\label{compact}
If $x,y\in D_{1}$ then
\begin{gather*}
\phi(\iota_{1}(x(1,-1)y^{\ast}))=\theta_{x\otimes\iota_{1}(-1,1),y\otimes \iota_{1}(-1,1)}
\end{gather*}
in $\mathcal{K}(\mathcal{E})$.
In particular $\phi(\iota_{1}(1,-1))=\theta_{e,e}$.
\end{Lemma}

\begin{proof}
For $n,m\in D_{1}$ and $x,y\in D_{1}$,
\begin{gather*}
\langle y\otimes\iota_{1}(-1,1),n\otimes\iota_{1}(m)\rangle =\langle\iota_{1}(-1,1),\iota_{1}\alpha(y^{\ast}n)\iota_{1} (m)\rangle
\\
\qquad
=\iota_{1}((-1,1)\alpha(y^{\ast}n)m)=\iota_{1}(\alpha(y^{\ast}n)m),
\end{gather*}
so
\begin{gather*}
\theta_{x\otimes\iota_{1}(-1,1),y\otimes\iota_{1}(-1,1)}(n\otimes\iota_{1}(m))=(x\otimes\iota_{1}(-1,1))\iota_{1}(\alpha(y^{\ast}n)m)
\\
\qquad
=x\otimes\iota_{1}((-1,1)\alpha(y^{\ast}n)m)=\allowbreak x\otimes\iota_{1}(\alpha(y^{\ast}n)m).
\end{gather*}
The previous lemma shows that this is equal to
\begin{gather*}
x(1,-1)y^{\ast}n\otimes\iota_{1}(m)=\phi(\iota_{1}(x(1,-1)y^{\ast}))(n\otimes\iota_{1}(m))
\end{gather*}
in $X/N$.
Setting $x=y=(1,-1)$ the last statement follows.
\end{proof}

Thus the projection $\iota_{1}(1,-1)$ of $C^{\ast}(D_{1},\preceq,\mathcal{M})$ is in the ideal
$\phi^{-1}(\mathcal{K}(\mathcal{E}))=J(\mathcal{E})$.
Def\/ine the ideal
\begin{gather*}
K=C^{\ast}(D_{1},\preceq,\mathcal{M})\iota_{1}(1,-1)C^{\ast}(D_{1},\preceq,\mathcal{M})
\end{gather*}
generated by $\iota_{1}(1,-1)$.
This is an ideal contained in the ideal $J(\mathcal{E})$.

\begin{Proposition}\sloppy
Let $(T,\pi):\mathcal{E}\rightarrow C$ be a~representation of the $C^*$-correspondence
  $_{C^{\ast}(D_{1},\preceq,\mathcal{M})}\mathcal{E}_{C^{\ast}(D_{1},\preceq,\mathcal{M})}$ in a~$C^*$-algebra~$C$ which
is coisometric on the ideal~$K$ generated by~$\iota_{1}(1,-1)$.
Then the partial isometry $T(e)$ has final projection $\pi(\iota_{1}(1,-1))$.
\end{Proposition}

\begin{proof}
Since $\iota_{1}(1,-1)\in K$ the coisometric hypothesis implies
\begin{gather*}
\pi(\iota_{1}(1,-1))=\Psi_{T}(\phi(\iota_{1}(1,-1))).
\end{gather*}
Lemma~\ref{compact} shows the latter is $\Psi_{T}(\theta_{e,e})=T(e)T(e)^{\ast}$.
\end{proof}

We noted that if $(T,\pi):\mathcal{E}\rightarrow C$ is a~representation of the $C^*$-correspondence
\linebreak $_{C^{\ast}(D_{1},\preceq,\mathcal{M})}\mathcal{E}_{C^{\ast}(D_{1},\preceq,\mathcal{M})}$ then $T(e)$ is a~partial
isometry with initial projection $T(e)^{\ast}T(e)=\pi(\iota_{1}(-1,1))$.
From the last proposition if $(T,\pi)$ is coisometric on the ideal~$K$ then the f\/inal projection
\begin{gather*}
T(e)T(e)^{\ast}=\pi(\iota_{1}(1,-1))\in\pi(C^{\ast}(D_{1},\preceq,\mathcal{M})).
\end{gather*}
Recall the semigroup $D_{1}$ is generated by the idempotent $(1,-1)$ and the map $\alpha$, which is implemented on the
image of $D_{1}$ in $\pi(C^{\ast}(D_{1},\preceq,\mathcal{M}))$ by the partial isometry $T(e)$.
It follows that $\pi(C^{\ast}(D_{1},\preceq,\mathcal{M}))$ is therefore contained in the $C^*$-algebra generated by
$T(e)$.
Combining this with the observations after Lemma~\ref{equiv} yields the following.

\begin{Theorem}
Let
\begin{gather*}
(T,\pi): \ {} _{C^{\ast}(D_{1},\preceq,\mathcal{M})}\mathcal{E}_{C^{\ast} (D_{1},\preceq,\mathcal{M})}\rightarrow C
\end{gather*}
be a~representation of the $C^*$-correspondence $\mathcal{E}$ in a~$C^*$-algebra~$C$ which is coisometric on the
ideal~$K$ of $C^{\ast}(D_{1},\preceq,\mathcal{M})$.
Then the $C^*$-algebra generated by the partial isometry $T(e)$ is equal to the $C^*$-algebra generated by
\begin{gather*}
T(\mathcal{E})\cup\pi(C^{\ast}(D_{1},\preceq,\mathcal{M})).
\end{gather*}
In particular the universal $C^*$-algebra $\mathcal{O}(K,\mathcal{E})$ is generated by the partial isometry
$T_{\mathcal{E}}(e)$, where $(T_{\mathcal{E}},\pi_{\mathcal{E}})$ is the universal representation of $\mathcal{E}$
coisometric on~$K$.
\end{Theorem}

\begin{Corollary}
\label{surj}
With $\mathcal{P}$ the universal $C^*$-algebra generated by a~partial isometry~$v$ there is a~surjective
$*$-homomorphism
\begin{gather*}
\psi: \ \mathcal{P}\rightarrow\mathcal{O}(K,\mathcal{E}),
\end{gather*}
which maps~$v$ to $T_{\mathcal{E}}(e)$.
\end{Corollary}

\begin{Lemma}
\label{coisometric}
Let $(T,\pi):\mathcal{E}\rightarrow C$ be a~representation of a~$C^*$-correspondence $_{B}\mathcal{E}_{B}$ in
a~$C^*$-algebra~$C$.
If $\Psi_{T}(\phi(b))=\pi(b)$ for an element $b\in J(\mathcal{E})$, then $\Psi_{T}(\phi(a))=\pi(a)$ for all~$a$ in the
closed two-sided ideal $BbB$ of~$B$ generated by~$b$.
\end{Lemma}

\begin{proof}
For $c,d\in B$ we have $cbd$ is in the ideal $J(\mathcal{E})=\phi^{-1}(\mathcal{K}(\mathcal{E}))$.
We have
\begin{gather*}
\Psi_{T}(\phi(cbd))=\Psi_{T}(\phi(c)\phi(b)\phi(d))=\pi(c)\Psi_{T}(\phi (b))\pi(d)
\end{gather*}
(e.g., by \cite[Proposition~1.2]{b}), which is $\pi(c)\pi(b)\pi(d)$.
\end{proof}

We are interested in constructing a~particular representation $(S,\pi_{\gamma})$ of the $C^*$-correspondence
$_{C^{\ast}(D_{1},\preceq,\mathcal{M})}\mathcal{E}_{C^{\ast}(D_{1},\preceq,\mathcal{M})}$ in the $C^*$-algebra
$\mathcal{P}$, the universal $C^*$-algebra generated by a~partial isometry $\upsilon$.
First, restrict the $*$-homomorphism $\sigma:A\rightarrow \mathcal{P}$ of the $*$-semi\-group~$A$ (Theorem~\ref{first
iso}) to the $*$-subsemigroup $D_{1}$ of~$A$ to obtain a~$*$-homomorphism $\lambda:D_{1}\rightarrow\mathcal{P}$.
We involve the complete order on $D_{1}$.

\begin{Proposition}
The $*$-homomorphism $\lambda:D_{1}\rightarrow\mathcal{P}$ is a~complete order representation of the completely ordered
$*$-semi\-group $(D_{1},\preceq,\mathcal{M})$ in the $C^*$-algebra $\mathcal{P}$.
\end{Proposition}

\begin{proof}
For $n\in D_{1}$ write $n=w_{\pm}^{\ast}(\pm1,\mp1)w_{\pm}$ a~product taking place in the semigroup~$A$.
We have
\begin{gather*}
\lambda(n)=\sigma(w_{\pm})^{\ast}\sigma(\pm1,\mp1)\sigma(w_{\pm})
\end{gather*}
in the $C^*$-algebra $\mathcal{P}$.
The latter product is less than or equal to the selfadjoint element
\begin{gather*}
\sigma(w_{\pm})^{\ast}\sigma(w_{\pm})=\sigma(w_{\pm}^{\ast}w_{\pm})=\lambda(w_{\pm}^{\ast}w_{\pm})
\end{gather*}
in $\mathcal{P}$.
Since the partial order on $D_{1}^{\rm sa}$ is generated by these relations, $\lambda$ is an order representation.
The complete order property for this map follows similarly.
\end{proof}

The universal property for $C^{\ast}(D_{1},\preceq,\mathcal{M})$ provides a~$*$-homomorphism
\begin{gather*}
\pi_{\lambda}:\ C^{\ast}(D_{1},\preceq,\mathcal{M})\rightarrow\mathcal{P}
\end{gather*}
with $\pi_{\lambda}\circ\iota_{1}=\lambda$.
For $n,m\in D_{1}$ def\/ine~$S$ on the elements $n\otimes\iota_{1}(m)$ of~$X$ by
\begin{gather*}
S(n\otimes\iota_{1}(m))=\lambda(n)v\pi_{\lambda}(\iota_{1}(m))=\pi_{\lambda}(\iota_{1}(n))v\pi_{\lambda}(\iota_{1}(m))
\end{gather*}
and extend linearly to the subspace $\mathbb{C}[D_{1}]\otimes_{\rm alg}\iota_{1}(\mathbb{C}[D_{1}])$ of~$X$.
For $x=s\otimes\iota_{1}(c)$ and $y=t\otimes\iota_{1}(d)$ with $c,d,s,t\in D_{1}$ compute
\begin{gather*}
\pi_{\lambda}(\langle x,y\rangle_{C^{\ast}(D_{1},\preceq,\mathcal{M})})=\pi_{\lambda}\langle
s\otimes\iota_{1}(c),t\otimes \iota_{1}(d)\rangle =\pi_{\lambda}(\iota_{1}(c)^{\ast}\iota_{1}
\alpha(s^{\ast}t)\iota_{1}(d))
\end{gather*}
while
\begin{gather*}
S(x)^{\ast}S(y)=(\pi_{\lambda}\iota_{1}(s)v\pi_{\lambda}(\iota_{1}(c)))^{\ast}\pi_{\lambda}\iota_{1}(t)v\pi_{\lambda}(\iota_{1}(d))
\\
\phantom{S(x)^{\ast}S(y)}
=\pi_{\lambda}(\iota_{1}(c))^{\ast}(v^{\ast}\pi_{\lambda}\iota_{1}(s^{\ast}t)v)\pi_{\lambda}(\iota_{1}(d))=\pi_{\lambda}(\iota_{1}(c))^{\ast}
\pi_{\lambda}(\iota_{1}\alpha(s^{\ast}t))\pi_{\lambda}(\iota_{1}(d)).
\end{gather*}
Therefore
\begin{gather*}
\Vert S(x)\Vert^{2}=\Vert S(x)^{\ast}S(x)\Vert =\Vert \pi_{\lambda}(\langle
x,x\rangle_{C^{\ast} (D_{1},\preceq,\mathcal{M})})\Vert \leq\Vert x\Vert^{2}
\end{gather*}
for $x,y\in\mathbb{C}[D_{1}]\otimes_{\rm alg}\iota_{1}(\mathbb{C}[D_{1}])$, and~$S$ becomes a~well def\/ined linear map on the
quotient of this subspace of~$X$ by~$N$.
Since~$S$ is bounded it may be extended by continuity to a~linear map, also denoted by~$S$, of $\mathcal{E}$ to
$\mathcal{P}$.

We claim that $(S,\pi_{\lambda})$ is a~correspondence representation.
For $k,m,n\in D_{1}$, and $b=\iota_{1}(k)\in C^{\ast}(D_{1},\preceq,\mathcal{M})$, we have
\begin{gather*}
S(\phi(b)(m\otimes\iota_{1}(n)))=S(l(k)(m\otimes\iota_{1}(n)))
\end{gather*}
by Theorem~\ref{left action}.
This equals $S(km\otimes\iota_{1} (n))=\pi_{\lambda}\iota_{1}(km)v\pi_{\lambda}(\iota_{1}(n))$ while
\begin{gather*}
\pi_{\lambda}(b)S(m\otimes\iota_{1}(n))
=\pi_{\lambda}(\iota_{1}\left(k\right))S(m\otimes\iota_{1}(n))
=\pi_{\lambda}(\iota_{1}\left(k\right))\pi_{\lambda}(\iota_{1}(m))v\pi_{\lambda}\iota_{1}(n),
\end{gather*}
and therefore
\begin{gather*}
S(\phi(b)(m\otimes\iota_{1}(n)))=\pi_{\lambda}(b)S(m\otimes\sigma(n)).
\end{gather*}
Check that
\begin{gather*}
S((m\otimes\iota_{1}(n))(b))=S(m\otimes\iota_{1}(nk))
=\pi_{\lambda}\iota_{1}(m)v\pi_{\lambda}\iota_{1}(n)\pi_{\lambda}\iota_{1}(k)=S(m\otimes\iota_{1}(n))\pi_{\lambda}(b).
\end{gather*}
Therefore $S(\phi(b)(x))=\pi_{\lambda}(b)S(x)$ and
\begin{gather*}
S(xb)=S\left(x\right) \pi_{\lambda}(b)
\end{gather*}
for $b\in C^{\ast}(D_{1},\preceq,\mathcal{M})$, $x\in\mathcal{E}$.
Therefore $(S,\pi_{\lambda})$ is a~representation of the correspondence
$_{C^{\ast}(D_{1},\preceq,\mathcal{M})}\mathcal{E}_{C^{\ast}(D_{1},\preceq,\mathcal{M})}$ in the $C^*$-algebra
$\mathcal{P}$.

\begin{Proposition}
There is a~representation
\begin{gather*}
(S,\pi_{\lambda}): \ \mathcal{E}\rightarrow\mathcal{P}
\end{gather*}
of the $C^*$-correspondence $_{C^{\ast}(D_{1},\preceq,\mathcal{M})} \mathcal{E}_{C^{\ast}(D_{1},\preceq,\mathcal{M})}$
which is coisometric on the ideal~$K$ of \linebreak $C^{\ast}(D_{1},\preceq,\mathcal{M})$ generated by $\iota_{1}(1,-1)$.
\end{Proposition}

\begin{proof}
We have just shown that $(S,\pi_{\lambda})$ is a~representation of the $C^*$-correspondence \linebreak
$_{C^{\ast}(D_{1},\preceq,\mathcal{M})}\mathcal{E}_{C^{\ast}(D_{1},\preceq,\mathcal{M})}$.
To show the representation is coisometric on~$K$ it is enough to show
\begin{gather*}
\Psi_{S}(\phi(\iota_{1}(1,-1)))=\pi_{\lambda}(\iota_{1}(1,-1))
\end{gather*}
by Lemma~\ref{coisometric}.
However,
\begin{gather*}
\pi_{\lambda}(\iota_{1}(1,-1))=\lambda((1,-1))=\sigma((1,-1))=vv^{\ast}.
\end{gather*}
We compute
\begin{gather*}
S(e)=S((1,-1)\otimes\iota_{1}(-1,1))=\lambda(1,-1)v\pi_{\lambda}\iota_{1}(-1,1)
\\
\qquad
=\lambda(1,-1)v\lambda(-1,1)=vv^{\ast}vv^{\ast}v=v.
\end{gather*}
By Lemma~\ref{compact} $\phi(\iota_{1}(1,-1))=\theta_{e,e}$, so
\begin{gather*}
\Psi_{S}(\phi(\iota_{1}(1,-1)))=\Psi_{S}(\theta_{e,e})=S(e)S(e)^{\ast}=vv^{\ast}=\pi_{\lambda}(\iota_{1}(1,-1)).\tag*{\qed}
\end{gather*}
\renewcommand{\qed}{}
\end{proof}

The universal property for the relative Cuntz--Pimsner $C^*$-algebra $\mathcal{O}(K,\mathcal{E})$ yields a~unique
$*$-homomorphism $\rho:\mathcal{O}(K,\mathcal{E})\rightarrow\mathcal{P}$ with
$(S,\pi_{\lambda})=\rho\circ(T_{\mathcal{E}},\pi_{\mathcal{E}})$.
Here $(T_{\mathcal{E}},\pi_{\mathcal{E}})$ is the canonical representation of $_{C^{\ast}
(D_{1},\preceq,\mathcal{M})}\mathcal{E}_{C^{\ast}(D_{1},\preceq,\mathcal{M})}$ to $\mathcal{O}(K,\mathcal{E})$.
Since $v=S(e)=\rho\circ T_{\mathcal{E}}(e)$, the image of $\rho$ contains a~generator~$v$ of $\mathcal{P}$, hence $\rho$
is surjective.

\begin{Theorem}
\label{iso}
The $*$-homomorphism $\rho:\mathcal{O}(K,\mathcal{E})\rightarrow \mathcal{P}$ is a~$*$-isomorphism.
\end{Theorem}

\begin{proof}
Since $\rho(\psi(v))=\rho(T_{\mathcal{E}}(e))=v$ by Corollary~\ref{surj} the composition $\rho\circ\psi$ is the
identity map on $\mathcal{P}$, so the surjection $\psi$ is an injection.
\end{proof}

\begin{Example}
If we consider quotient semigroups of $D_{0}$ and $D_{1}$ many aspects of this analysis still hold.
For example, let $S_{1}$ be the unital (unit $u$) two element $*$-semi\-group $\{u,s\}$ with~$s$ also a~selfadjoint
idempotent.
View this as ordered, with the trivial partial order $a\preceq a$ if\/f $a=a$.
With $S_{0}$ the $*$-subsemigroup generated by~$u$, then $C^{\ast}(S_{1},\preceq)=C^{\ast}(S_{1})$ is isomorphic to the
unital $C^*$-algebra $\mathbb{C}^{2}$, and $C^{\ast}(S_{0},\preceq)\cong\mathbb{C}$.
The $C^*$-correspondence constructed above is the $C^*$-algebra $\mathbb{C}^{2}$ viewed as a~right Hilbert module
$\mathcal{E}_{C^{\ast}(S_{1})}=C^{\ast} (S_{1})_{C^{\ast}(S_{1})}$ over itself.
The left action is trivial:
\begin{gather*}
\phi: \ C^{\ast}(S_{1})\rightarrow\mathcal{L}(\mathcal{E}_{C^{\ast}(S_{1})})\cong C^{\ast}(S_{1})
\end{gather*}
with $\phi(a)$ the unit for all $a\in C^{\ast}(S_{1})$.
Set~$K$ to be the ideal, isomorphic to $\mathbb{C}$, generated by~$s$ in $C^{\ast}(S_{1})$.
It follows that the universal Cuntz--Pimsner $C^*$-algebra $\mathcal{O} (K,\mathcal{E})$ is the Toeplitz $C^*$-algebra
generated by an isometry.
\end{Example}

Although not necessary for the development up to this point, it is the case (Corollary~\ref{need}) that the ideal~$K$ of
$C^{\ast}(D_{1},\preceq,\mathcal{M})$ generated by the projection $\iota_{1}(1,-1)$ is actually contained in the ideal
$J_{\mathcal{E}}=\phi^{-1}(\mathcal{K}(\mathcal{E}))\cap(\ker\phi)^{\perp}$.
First recall the Hilbert module $\mathcal{E}_{0}=\mathcal{E}_{C^{\ast}(D_{0},\preceq,\mathcal{M})}$ with left action
\begin{gather*}
\phi_{0}: \ C^{\ast}(D_{1},\preceq,\mathcal{M})\rightarrow\mathcal{L} (\mathcal{E}_{0})
\end{gather*}
described in Remark~\ref{restriction}.
In the following
\begin{gather*}
\pi_{\omega}: \ C^{\ast}(D_{0},\preceq,\mathcal{M})\rightarrow C^{\ast} (D_{1},\preceq,\mathcal{M})
\end{gather*}
is the $*$-homomorphism obtained by applying the universal property to the complete order representation
(Proposition~\ref{0 to 1}) $\iota_{1}\circ \omega:D_{0}\rightarrow C^{\ast}(D_{1},\preceq,\mathcal{M})$.

\begin{Proposition}
With $a\in C^{\ast}(D_{1},\preceq,\mathcal{M})$, and
$
e=(1,-1)\otimes\iota_{0}(-1,1)$
in $\mathcal{E}_{0}=\mathcal{E}_{C^{\ast}(D_{0},\preceq,\mathcal{M})}$ we have
\begin{gather*}
\pi_{\omega}(\langle e,\phi_{0}(a)e\rangle)=\iota_{1} (-1,1)a\iota_{1}(-1,1).
\end{gather*}
\end{Proposition}

\begin{proof}
For $s\in D_{1}$ we have
\begin{gather*}
\langle e,\phi_{0}(\iota_{1}(s))e\rangle =\langle e,s(1,-1)\otimes\iota_{0}(-1,1)\rangle
=\iota_{0}((-1,1)\alpha(s(1,-1))(-1,1)),
\end{gather*}
which is equal to $\iota_{0}(\alpha(s))$ by Remark~\ref{b}.
Applying $\pi_{\omega}$ and using $\pi_{\omega}\circ\iota_{0}=\iota_{1}\circ\omega$ we obtain
\begin{gather*}
\iota_{1}(\omega(\alpha(s)))=\iota_{1}((-1,1)s(-1,1)).
\end{gather*}
Linearity implies that the equality holds for~$a$ in the dense $*$-subalgebra $\mathbb{C}[\iota_{1}(D_{1})]$ of
\linebreak $C^{\ast}(D_{1},\preceq,\mathcal{M})$.
Since both sides are continuous the result follows.
\end{proof}

The proof of the following inclusion follows using Corollary~\ref{injection} via Remark~\ref{restriction}.
It is not clear whether or not this is a~strict inclusion.

\begin{Corollary}
\label{need}
The ideal~$K$ of $C^{\ast}(D_{1},\preceq,\mathcal{M})$ generated by $\iota_{1}(1,-1)$ is contained in the ideal~$J_{\mathcal{E}}$.
\end{Corollary}

\begin{proof}
It is suf\/f\/icient to show $\iota_{1}(1,-1)\in(\ker\phi)^{\perp}$, so to show $a\iota_{1}(1,-1)=0$ whenever \mbox{$\phi(a)=0$}
for $a\in C^{\ast}(D_{1},\preceq,\mathcal{M})$.
If $\phi(a)=0$ then so is $\phi(a^{\ast}a)$ and (Remark~\ref{restriction}) therefore also $\phi_{0}(a^{\ast}a)$.
The previous proposition shows $\iota_{1}(-1,1)a^{\ast}a\iota_{1}(-1,1)$, and therefore $a\iota_{1}(-1,1)$, is zero in
$C^{\ast}(D_{1},\preceq,\mathcal{M})$.
\end{proof}

Observe that we did not have recourse to the gauge-invariant uniqueness theorem for Cuntz--Pimsner algebras to show that $\rho$ is an isomorphism.
However, along with $\rho$ being an isomorphism, Corollary~\ref{need} now enables the eventual use of this theorem
in establishing the next corollaries.
The theorem~\cite{k} states that if $(T,\pi):\mathcal{E}\rightarrow B$ is a~representation of a~correspondence
$\mathcal{E}$ in a~$C^*$-algebra~$B$ which is coisometric on $J_{\mathcal{E}}$ then the induced $*$-homomorphism
$\rho:\mathcal{O}_{\mathcal{E}}\rightarrow C^{\ast}(T,\pi)$ is an isomorphism if and only if $\pi$ is injective and
$(T,\pi)$ admits a~gauge action.
Recall that a~representation $(T,\pi):\mathcal{E}\rightarrow B$ is said to admit a~gauge action if there is
$\gamma:\mathbb{T}\rightarrow\operatorname{Aut} C^{\ast}(T,\pi)$ a~homomorphism with $\gamma_{t}(T(e))=tT(e)$ for all
$e\in\mathcal{E}$ and $\gamma_{t}(\pi(a))=\pi(a)$ for all $a\in A$, $(t\in\mathbb{T})$.
Since the $*$-homomorphism $\Psi_{T}$ is def\/ined by mapping
\begin{gather*}
\theta_{x,y}\rightarrow T(x)T(y)^{\ast}\in C^{\ast}(T,\pi)
\end{gather*}
we necessarily have
\begin{gather*}
\gamma_{t}(\Psi_{T}\left(\theta\right))=\Psi_{T}\left(\theta\right)
\end{gather*}
for $\theta\in\mathcal{K}(\mathcal{E})$ and $t\in\mathbb{T}$.

For~$K$ any ideal in $J(\mathcal{E})$ it follows from the universal property for $\mathcal{O}(K,\mathcal{E})$ that the
universal representation $(T_{\mathcal{E}},\pi_{\mathcal{E}})$ of $\mathcal{E}$ admits a~gauge action, called the
canonical gauge action, on $\mathcal{O}(K,\mathcal{E})$.
We are unable to apply the gauge invariant uniqueness theorem directly to the representation $(S,\pi_{\lambda})$ since
it is coisometric on the ideal~$K$ generated by $\iota_{1}(1,-1)$, which is not necessarily the ideal
$J_{\mathcal{E}}$.

\begin{Corollary}
\label{iso image}
The $*$-representation
\begin{gather*}
\pi_{\lambda}: \ C^{\ast}(D_{1},\preceq,\mathcal{M})\rightarrow\mathcal{P}
\end{gather*}
is injective.
Therefore $C^{\ast}(D_{1},\preceq,\mathcal{M})$ is residually finite.
\end{Corollary}

\begin{proof}
The guage-invariant uniqueness theorem~\cite{k} ensures that the universal $J_{\mathcal{E}}$ covariant representation
$(T,\pi):\mathcal{E}\rightarrow \mathcal{O}(J_{\mathcal{E}},\mathcal{E})=\mathcal{O}_{\mathcal{E}}$ is injective, i.e.,
$\pi$ is injective.
Since $K\subseteq J_{\mathcal{E}}$ by the previous corollary, this representation is also~$K$ covariant, so the
universal property for $\mathcal{O}(K,\mathcal{E})$ implies that there is
\begin{gather*}
\xi: \ \mathcal{O}(K,\mathcal{E})\rightarrow\mathcal{O}(J_{\mathcal{E}},\mathcal{E})
\end{gather*}
a~(surjective) $*$-homomorphism with
\begin{gather*}
\xi\circ(T_{\mathcal{E}},\pi_{\mathcal{E}})=(T,\pi).
\end{gather*}
In particular $\xi\circ\pi_{\mathcal{E}}=\pi$.
Since $\pi$ is injective, so is $\pi_{\mathcal{E}}$.
Now recall (paragraph preceding Theorem~\ref{iso}) $\rho\circ\pi_{\mathcal{E}}=\pi_{\lambda}$.
Since $\rho$ is an isomorphism $\pi_{\lambda}$ is injective.

Therefore $C^{\ast}(D_{1},\preceq,\mathcal{M})$ is isomorphic to a~$C^*$-subalgebra of $\mathcal{P}$.
Since $\mathcal{P}$ is residually f\/inite-dimensional~\cite{bn} so is $C^{\ast}(D_{1},\preceq,\mathcal{M})$
(cf.~\cite{bo}).
\end{proof}

Now $\mathcal{P}$ is not exact~\cite{bn},
and is a~quotient of the Toeplitz--Pimsner algebra of the correspondence
$\mathcal{E}$ over $C^{\ast}(D_{1},\preceq,\mathcal{M})$.
However, the latter is exact if and only if $C^{\ast}(D_{1},\preceq,\mathcal{M})$ is \cite{bo, ds}.

\begin{Corollary}
There are natural non injective quotient maps of the universal and nonexact $C^*$-algebras:
\begin{gather*}
C^{\ast}(D_{1})\overset{\pi_{\preceq}}{\rightarrow}C^{\ast}(D_{1},\preceq)\overset{\pi_{\mathcal{M}}}{\rightarrow}C^{\ast}(D_{1},\preceq,\mathcal{M}).
\end{gather*}
\end{Corollary}

\begin{proof}
We show this for $C^{\ast}(D_{1},\preceq)$ and $C^{\ast}(D_{1},\preceq,\mathcal{M})$ as a~similar argument holds for $C^{\ast}(D_{1})$.

Let $\iota:D_{1}\rightarrow C^{\ast}(D_{1})$, $\iota_{\preceq}:(D_{1},\preceq)\rightarrow C^{\ast}(D_{1},\preceq)$, and
$\iota_{\mathcal{M}}:(D_{1},\preceq,\mathcal{M})\rightarrow C^{\ast}(D_{1},\preceq,\mathcal{M})$ respectively denote the
canonical $*$-representation, order representation, and complete order representation.
Since $\iota_{\mathcal{M}}$ is also an order representation of $(D_{1},\preceq)$, the universal property yields
a~natural quotient map $\pi_{\mathcal{M}}:C^{\ast}(D_{1},\preceq)\rightarrow C^{\ast}(D_{1},\preceq,\mathcal{M})$ with
$\pi_{\mathcal{M}}\circ\iota_{\preceq}=\iota_{\mathcal{M}}$.
The comment preceding the corollary shows $C^{\ast}(D_{1},\preceq,\mathcal{M})$ is not exact, therefore
$C^{\ast}(D_{1},\preceq)$ is not exact.
Similarly, $C^{\ast}(D_{1})$ is not exact.

Recall that any $C^*$-homomorphism must also be a~complete order map.
By considering appropriate diagrams it follows that the map $\pi_{\mathcal{M}}$ is injective, so an isomorphism, if and
only if the order representation $\iota_{\preceq}$ is a~complete order representation.
However, there are order representations of $(D_{1},\preceq)$ that are not complete order representations of
$(D_{1},\preceq,\mathcal{M})$ (cf.\ Remarks~\ref{order} and~\ref{comp order}) so it follows that $\iota_{\preceq}$
cannot be a~complete order representation of $(D_{1},\preceq,\mathcal{M})$ in $C^{\ast}(D_{1},\preceq)$, and
$\pi_{\mathcal{M}}$ is not injective.
\end{proof}

\subsection*{Acknowledgements}

I~am most grateful to the referees for their detailed and helpful commentary regarding the many changes aimed at
improving the readability of my initial submission.
As well, a~referee suggested a~shorter proof of Proposition~\ref{basic irr} and pointed out the potential for a~natural
approach to the pair of $*$-maps $\omega$ and $\beta_{\omega}$ of Section~\ref{section2.1}.
I~am also thankful to the Fields Institute for their hospitality throughout the fall of 2012 during which some of this project was completed.

\LastPageEnding


\begin{thebibliography}{99}
\footnotesize \itemsep=0pt\addcontentsline{toc}{section}{References}

\bibitem{aee}
Abadie B., Eilers S., Exel R., Morita equivalence for crossed products by
  {H}ilbert {$C^*$}-bimodules, \href{http://dx.doi.org/10.1090/S0002-9947-98-02133-3}{\textit{Trans. Amer. Math. Soc.}} \textbf{350}
  (1998), 3043--3054.

\bibitem{a}
Arveson W., Noncommutative dynamics and {$E$}-semigroups, \href{http://dx.doi.org/10.1007/978-0-387-21524-2}{\textit{Springer Monographs
  in Mathematics}}, Springer-Verlag, New York, 2003.

\bibitem{bd}
Barnes B.A., Duncan J., The {B}anach algebra {$l^{1}(S)$}, \href{http://dx.doi.org/10.1016/0022-1236(75)90032-4}{\textit{J.~Funct.
  Anal.}} \textbf{18} (1975), 96--113.

\bibitem{b}
Brenken B., Topological quivers as multiplicity free relations, \textit{Math.
  Scand.} \textbf{106} (2010), 217--242.

\bibitem{bn}
Brenken B., Niu Z., The {$C^*$}-algebra of a partial isometry, \href{http://dx.doi.org/10.1090/S0002-9939-2011-10988-2}{\textit{Proc.
  Amer. Math. Soc.}} \textbf{140} (2012), 199--206.

\bibitem{bo}
Brown N.P., Ozawa N., {$C^*$}-algebras and f\/inite-dimensional approximations,
  \textit{Graduate Studies in Mathe\-ma\-tics}, Vol.~88, Amer. Math.
  Soc., Providence, RI, 2008.

\bibitem{cdp}
Conway J.B., Duncan J., Paterson A.L.T., Monogenic inverse semigroups and their
  {$C^\ast$}-algebras, \href{http://dx.doi.org/10.1017/S030821050002552X}{\textit{Proc. Roy. Soc. Edinburgh Sect.~A}} \textbf{98}
  (1984), 13--24.

\bibitem{dls}
Dales H.G., Lau A.T.M., Strauss D., Banach algebras on semigroups and on their
  compactif\/ications, \href{http://dx.doi.org/10.1090/S0065-9266-10-00595-8}{\textit{Mem. Amer. Math. Soc.}} \textbf{205} (2010),
  vi+165~pages.

\bibitem{dp}
Duncan J., Paterson A.L.T., {$C^\ast$}-algebras of inverse semigroups,
  \href{http://dx.doi.org/10.1017/S0013091500003187}{\textit{Proc. Edinburgh Math. Soc.}} \textbf{28} (1985), 41--58.

\bibitem{ds}
Dykema K.J., Shlyakhtenko D., Exactness of {C}untz--{P}imsner {$C^*$}-algebras,
  \href{http://dx.doi.org/10.1017/S001309159900125X}{\textit{Proc. Edinb. Math. Soc.}} \textbf{44} (2001), 425--444,
  \href{http://arxiv.org/abs/math.OA/9911002}{math.OA/9911002}.

\bibitem{f}
Fell J.M.G., Doran R.S., Representations of {$*$}-algebras, locally compact
  groups, and {B}anach {$*$}-algebraic bundles. {V}ol.~1. Basic representation
  theory of groups and algebras, \textit{Pure and Applied Mathematics}, Vol.~125, Academic Press, Inc., Boston, MA, 1988.

\bibitem{f2}
Fell J.M.G., Doran R.S., Representations of {$*$}-algebras, locally compact
  groups, and {B}anach {$*$}-algebraic bundles. {V}ol.~2. Banach $*$-algebraic
  bundles, induced representations, and the generalized Mackey analysis,
  \textit{Pure and Applied Mathematics}, Vol.~126, Academic Press, Inc.,
  Boston, MA, 1988.

\bibitem{fmr}
Fowler N.J., Muhly P.S., Raeburn I., Representations of {C}untz--{P}imsner
  algebras, \href{http://dx.doi.org/10.1512/iumj.2003.52.2125}{\textit{Indiana Univ. Math.~J.}} \textbf{52} (2003), 569--605.

\bibitem{h}
Howie J.M., An introduction to semigroup theory, Academic Press, London~-- New
  York, 1976.

\bibitem{k}
Katsura T., On {$C^*$}-algebras associated with {$C^*$}-correspondences,
  \href{http://dx.doi.org/10.1016/j.jfa.2004.03.010}{\textit{J.~Funct. Anal.}} \textbf{217} (2004), 366--401,
  \href{http://arxiv.org/abs/math.OA/0309088}{math.OA/0309088}.

\bibitem{l}
Lance E.C., Hilbert {$C^*$}-modules. A~toolkit for operator algebraists,
  \href{http://dx.doi.org/10.1017/CBO9780511526206}{\textit{London Mathematical Society Lecture Note Series}}, Vol.~210, Cambridge
  University Press, Cambridge, 1995.

\bibitem{li}
Li X., Semigroup {$C^*$}-algebras and amenability of semigroups,
  \href{http://dx.doi.org/10.1016/j.jfa.2012.02.020}{\textit{J.~Funct. Anal.}} \textbf{262} (2012), 4302--4340, \href{http://arxiv.org/abs/1105.5539}{arXiv:1105.5539}.

\bibitem{ms}
Muhly P.S., Solel B., Tensor algebras over {$C^*$}-correspondences:
  representations, dilations, and {$C^*$}-envelopes, \href{http://dx.doi.org/10.1006/jfan.1998.3294}{\textit{J.~Funct. Anal.}}
  \textbf{158} (1998), 389--457.

\bibitem{p}
Paulsen V., Completely bounded maps and operator algebras, \textit{Cambridge
  Studies in Advanced Mathematics}, Vol.~78, Cambridge University Press,
  Cambridge, 2002.

\bibitem{pm}
Pimsner M.V., A class of {$C^*$}-algebras generalizing both {C}untz--{K}rieger
  algebras and crossed products by~{${\mathbb Z}$}, in Free Probability Theory
  ({W}aterloo, {ON}, 1995), \textit{Fields Inst. Commun.}, Vol.~12, Amer. Math.
  Soc., Pro\-vi\-dence, RI, 1997, 189--212.

\end{thebibliography}
\end{document}